\documentclass[12pt]{article}
\usepackage{amssymb}
\usepackage{tikz}\usepackage{mathrsfs}
\font\eufm=eufm10 at 14pt\font\eufms=eufm10\font\eufmss=eufm7\newfam\eufam
\textfont\eufam=\eufm\scriptfont\eufam=\eufms\scriptscriptfont\eufam=\eufmss

\def\build#1_#2^#3{\mathrel{\mathop{\kern 0pt#1}\limits_{#2}^{#3}}}
\def\Z{{\bf Z}}\def\R{{\bf R}}\def\C{{\bf C}}

\def\rde{\mathscr}

\def\vc{{\scriptstyle\wedge }}
\def\cqfd{\hfill\vbox{\hrule\hbox{\vrule height6pt depth0pt\hskip 6pt \vrule height6pt}\hrule\relax}}
\def\noi{\noindent}\def\e{{\varepsilon}}
\def\hfl#1#2{\smash{\mathop{\hbox to 12 mm{\rightarrowfill}}\limits^{\scriptstyle#1}_{\scriptstyle#2}}}
\def\vfl#1#2{\llap{$\scriptstyle #1$}\left\downarrow\vbox to 6mm{}\right.\rlap{$\scriptstyle #2$}}
\def\diagram#1{\def\normalbaselines{\baselineskip=0pt\lineskip=10pt\lineskiplimit=1pt} \matrix{#1}}
\def\minus{\setminus }
\begin{document}
\def\zput(#1,#2)#3{\setbox40 =\hbox{#3}\dimen40 = #1pt
\dimen41 =\wd40\divide\dimen41 by 2\advance\dimen40 by -\dimen41
\dimen41 = #2pt\advance\dimen41 by -2.5pt
{\def\unitlength{}\put(\dimen40,\dimen41){\box40}\def\unitlength{1pt}}}

\overfullrule=0pt
\
\vskip 64pt
\centerline{\bf Functoriality of Khovanov homology}
\vskip 12pt
\centerline{Pierre Vogel\footnote{Universit\'e Paris Diderot, Institut de Math\'ematiques de Jussieu-Paris Rive Gauche (UMR 7586),
B\^atiment Sophie Germain, Case 7012, 75205--Paris Cedex 13 --- Email: pierre.vogel@imj-prg.fr}}
\vskip 48pt
\noi{\bf Abstract.} In this paper we prove that every Khovanov homology associated to a Frobenius algebra of rank $2$ can be modified in such a way as to produce a TQFT
on oriented links, that is a monoidal functor from the category of cobordisms of oriented links to the homotopy category of complexes.
\vskip 12pt
\noi{\bf Keywords:} Frobenius algebra, Khovanov homology, cobordism of oriented links, monoidal functor.

\noi{\bf Mathematics Subject Classification (2010):} 13Axx, 13Dxx, 57M25
\vskip 24pt
\noi{\bf Introduction.}
\vskip 12pt
The Khovanov homology was introduced by Khovanov [Kh1] as a categorification of the Jones polynomial. The main ingredient of this homology is a monoidal functor from
the category of cobordisms of oriented $1$-manifolds to a category of modules. But such monoidal functors are characterized by commutative Frobenius algebras [Ko].
Therefore the Khovanov homology can be defined for every commutative Frobenius algebra $R$. The classical Khovanov homology corresponds to the case: 
$R=\Z[\alpha]/(\alpha^2)$ and the Lee-version of this homology corresponds to the case: $R=\Z[\alpha]/(\alpha^2-1)$.

An important problem in this theory is to extend the Khovanov homology to a monoidal functor from the category of cobordisms of oriented links. The first attempt by
Khovanov gave a negative answer because of many problems of signs. Functoriality up to sign was conjectured by Khovanov and proved later by Jacobson [Ja], Bar Natan [BN2]
and Khovanov [Kh3]. This functoriality up to sign was used by Rasmussen [Ra1] to prove a conjecture of Milnor about the slice genus. Strict functoriality for alternative
versions of Khovanov homology were proven by Blanchet [Bl] and Clark, Morrison, Walker [CMW]. In these versions of Khovanov homology, the Kauffman bracket is replaced
with the $sl_2$-polynomial in [Bl] and the $su_2$-polynomial in [CMW].

In this paper we will prove functoriality for the Khovanov homology associated to any Frobenius algebra of rank $2$ (and the classical Kauffman bracket).

Suppose $K$ is a commutative ring and $R$ is a Frobenius $K$-algebra of rank $2$. Denote by $u\mapsto \overline u$ the involution of the extension $K\subset R$ and by
$\delta$ the image of $1$ under the composite map:
$$R\ \hfl{\hbox{\tiny coproduct}}{}\ R\otimes R\ \hfl{\hbox{\tiny product}}{}\ R$$
Denote by ${\rde L}$ the category of cobordisms of oriented links in $\R^3$ and by ${\rde C}_K$ the homotopy category of $K$-complexes. These categories are both 
monoidal. The main result of this paper is the following:
\vskip 12pt
\noi{\bf Theorem A:} {\sl There exists a monoidal functor $\Psi$ from ${\rde L}$ to ${\rde C}_K$ satisfying the following property:

for any diagram $D$ of an oriented link $L$, $\Psi(L)$ is isomorphic to the classical Khovanov complex of $D$.}
\vskip 12pt
Such functors are not unique. The construction produces a lot of functors (called Khovanov functors) satisfying this property. There is a well defined invertible element
in $R$ associated to each Khovanov functor: its weight. Actually if $C$ is an unknotted cobordism of genus $1$ from the unknot diagram (without any crossing) to itself, 
the image of this cobordism under a Khovanov functor of weight $\pi$ is the multiplication by $\delta/\pi$ (resp. $\overline\delta/\overline\pi$) from $R$ to $R$ if the
unknot is oriented clockwise (resp. counterclockwise).
\vskip 12pt
\noi{\bf Theorem B:} {\sl For every invertible element $\pi\in R$ there is a Khovanov functor of weight $\pi$. Moreover two Khovanov functors with the same weight are 
isomorphic.}
\vskip 12pt
Consider a closed oriented surface $S$ in $\R^4$. This surface may be consider as a cobordism from the empty link to itself. Therefore every Khovanov functor $\Psi$
induces an invariant $\Psi(S)\in K$. In the classical case the functor $\Psi$ was well defined up to sign and $\Psi(S)$ was determined by Tanaka [Ta] and Rasmussen [Ra2],
at least for connected surfaces. They prove that $\Psi(S)$ doesn't depend on the embedding $S\subset \R^4$. This result is still true in the general case and $\Psi(S)$
depends only on $S$ and the weight of $\Psi$ and not on the embedding. More precisely we have the following result:
\vskip 12pt
\noi{\bf Theorem C:} {\sl Let $\Psi$ be a Khovanov functor of weight $\pi\in R^*$ and $S$ be a closed oriented surface in $\R^4$. Then we have:
$$\Psi(S)=\prod_i\e(\delta^{p_i}\pi^{1-p_i})$$
where the $p_i$'s are the genus of the components of $S$.}
\vskip 24pt
\noi{\bf 1. Frobenius algebras of rank $2$.}
\vskip 12pt
\noi{\bf 1.1 Definition:} {\sl Let $K$ be a commutative ring. A $K$-algebra of rank $2$ is a $K$-algebra isomorphic to $K[\alpha]/(P)$, where $P\in K[\alpha]$ is a monic
polynomial of degree $2$. 

The element $\alpha$ is called a $K$-generator of $R$ and the involution of the extension $K\subset R$ is called the involution of $R$.}
\vskip 12pt
Let $R$ be a $K$-algebra of rank $2$ and $\alpha$ be a $K$-generator of $R$. The polynomial $P$ is given by:
$$P(\alpha)=\alpha^2-s\alpha+p$$
with $s$ and $p$ in $K$ and the involution of $R$, denoted by $u\mapsto\overline u$, is the identity on $K$ and sends $\alpha$ to $\overline\alpha=s-\alpha$. So we have:
$$s=\alpha+\overline\alpha\hskip 48pt p=\alpha\overline\alpha$$
and, for any $u\in R$, $u+\overline u$ and $u\overline u$ are in $K$.
\vskip 12pt
\noi{\bf 1.2 Proposition: } {\sl Let $K$ be a commutative ring and $R$ be a Frobenius $K$-algebra. Suppose $R$ is a $K$-algebra of rank $2$ and $\alpha$ is a 
$K$-generator of $R$. Then there is a unique invertible element $\omega$ in $R$ such that the coproduct and the counit are defined by:
$$\forall u\in R,\ \ \Delta(u)=u\otimes\omega\alpha-u\overline\alpha\otimes\omega$$
$$\e(\omega)=0\hskip 48pt \e(\omega\alpha)=1$$
$\omega$ is called the twisting element of $R$.

Moreover, if $\omega$ is any invertible element in $R$, the coproduct and the counit defined by the formulae above induce a structure of Frobenius $K$-algebra
on $R$.}
\vskip 12pt
\noi{\bf Proof:} Since $\Delta$ is a $K$-linear map from $R$ to $R\otimes R$, there exist two $K$-linear map $f$ and $g$ from $R$ to $R$ such that:
$$\forall u \in R,\ \ \ \Delta(u)=f(u)\otimes 1+g(u)\otimes\alpha$$
Since $R$ is a Frobenius algebra we have the relation: $\Delta(u)=(u\otimes1)\Delta(1)$ and then:
$$\forall u\in R,\ \ f(u)=uf(1)\ \ \ g(u)=ug(1)$$
We have also the relation: $(\alpha\otimes1-1\otimes\alpha)\Delta(1)=0$ which implies: $f(1)+\overline\alpha g(1)=0$. By setting: $\omega=g(1)$ we get:
$$\Delta(u)=u\omega\otimes\alpha-u\omega\overline\alpha\otimes1=u\otimes\omega\alpha-u\overline\alpha\otimes\omega$$

On the other hand the counit satisfies the relation: $(1\otimes\e)\Delta(u)=u$ which is equivalent to:
$$\e(\omega)=0\hskip 48pt\e(\omega\alpha)=1$$

The last thing to do is to prove that $\omega$ is invertible in $R$. Set: $u=\e(\alpha)-\overline\alpha\e(1)$. It is easy to see the following:
$$\e(\omega u)=\e(1)\hskip 48pt\e(\omega u\alpha)=\e(\alpha)$$
So for any $v\in R$, we have: $\e((\omega u-1)v)=0$. Set: $\omega u-1=a+b\alpha$ with $a$ and $b$ in $K$. Testing this formula with $v=\omega$ and 
$v=\omega\alpha$ implies: $b=a=0$. Therefore $\omega$ is invertible with inverse $u$.

The last part of the proposition is easy to check.\cqfd
\vskip 12pt
\noi{\bf 1.3 Remarks:} Let $R_0$ be the ring $\Z[\alpha,\overline\alpha,a,b,(a+b\alpha)^{-1},(a+b\overline\alpha)^{-1}]$. This ring is equipped with an involution keeping
$a$ and $b$ fixed and exchanging $\alpha$ and $\overline\alpha$. The ring of invariant elements under this involution is:
$$K_0=\Z[\alpha+\overline\alpha,\alpha\overline\alpha,a,b,((a+b\alpha)(a+b\overline\alpha))^{-1}]$$
and $R_0$ is a $K_0$-algebra of rank $2$ with $K_0$-generator $\alpha$. By setting: $\omega=a+b\alpha$, we see that $R_0$ is a Frobenius algebra with twisting element
$\omega$. Moreover $R_0$ is universal in the following sense:

Let $R$ be a Frobenius $K$-algebra of rank $2$ and $\beta$ be a $K$-generator of $R$. Then there exists a unique Frobenius algebra homomorphism from $R_0$ to $R$
sending $\alpha$ to $\beta$.

In particular the endomorphisms of the Frobenius algebra $R_0$ are ring homomorphisms characterized by:
$$\alpha\mapsto \lambda\alpha+\mu\hskip 48pt \overline\alpha\mapsto\lambda\overline\alpha+\mu$$
$$a\mapsto\lambda^{-1}a-\lambda^{-2}\mu b\hskip 48pt b\mapsto \lambda^{-2} b$$
where $(\lambda,\mu)$ is any element of $K_0^*\times K_0$.

Another description of the universal algebra $R_0$ was founded by Khovanov in [Kh4].
\vskip 24pt
\noi{\bf 1.4} From now on $K$ will be a commutative ring, $R$ a Frobenius $K$-algebra of rank $2$ and $\alpha$ a $K$-generator of $R$. The twisting element in $R$
will be denoted by $\omega$. Set:
$$t=\e(1)\hskip24pt \delta=\omega(\alpha-\overline\alpha)\hskip 24pt\theta={\omega\over\overline\omega}$$
It is easy to see that these elements do not depend on the choice of $\alpha$. Moreover we have the following relations:
$$t\delta=1-\theta\hskip 24pt\overline t=t\hskip 24pt \overline \delta=-\theta^{-1}\delta\hskip 24pt\overline\theta=\theta^{-1}$$
Set: $A=\Z[t,\delta,(1-t\delta)^{-1}]$. This ring is a ring with involution and it is contained in the algebra $R_0$. Therefore every Frobenius algebra of rank $2$ is an
$A$-algebra.

By setting: $s=\delta\overline\delta=-\theta^{-1}\delta^2=-(1-t\delta)^{-1}\delta^2$, it is not difficult to see that the ring of elements in $A$ that are fixed by the
involution is $\Z[s,t]$ and $A$ is a $\Z[s,t]$-algebra of rank $2$ with $\Z[s,t]$-generator $\delta$.
\vskip 12pt
\noi{\bf 1.5 Lemma: } {\sl The ring of elements in $R_0$ which are invariant under every endomorphism of $R_0$ is the ring $A$.}
\vskip 12pt
\noi{\bf Proof:} 
Suppose that $a_1,\dots, a_p$ are elements in a commutative ring $\Lambda$. Then denote by $\Lambda<a_1,\dots,a_p>$ the ring obtained by inverting $a_1,\dots, a_p$ in 
$\Lambda$. So we have:
$$R_0=\Z[\alpha,\overline\alpha,a,b]<a+b\alpha,a+b\overline\alpha>$$
We have the following: 
$$\omega=a+b\alpha\hskip 24pt t=\e(1)={-b\over\omega\overline\omega}\hskip 24pt \overline\omega(1-t\delta)=\omega$$
and then
$$R_0=\Z[\alpha,\alpha-\overline\alpha,b,\omega]<\omega,\overline\omega)>=\Z[\alpha,\delta,b,\omega]<\omega,\overline\omega>$$
$$=\Z[\alpha,t,\delta,\omega]<\omega,1-t\delta>=\Bigl(\Z[\alpha,t,\delta]<1-t\delta>\Bigr)[\omega^{\pm}]$$
Consider the endomorphisms of $R_0$ sending $\alpha$ to $\alpha+\mu$, for some $\mu\in K_0$. These endomorphisms keep $t$, $\delta$ and $\omega$ fixed. Then the ring 
$R_1$ of elements in $R_0$ which are invariant under these endomorphisms is:
$$R_1=\Z[t,\delta,\omega]<\omega,1-t\delta>=\Bigl(\Z[t,\delta]<1-t\delta>\Bigr)[\omega^{\pm}]$$
But every endomorphism of $R_0$ keeps $t$ and $\delta$ fixed, and multiplies $\omega$ by any element $\lambda\in K_0^*$. Since $K_0^*$ contains any power of 
$\omega\overline \omega$, the ring $R_2$ of elements in $R_0$ which are invariant under every endomorphism of $R_0$ is:
$$R_2=\Z[t,\delta]<1-t\delta>=A\eqno\cqfd$$
\vskip 12pt
Denote by ${\rde C}$ the category of cobordisms of oriented curves (i.e. closed oriented $1$-dimensional manifolds). The disjoint union induces on ${\rde C}$ a monoidal
structure. TQFT's for oriented surfaces are in one to one correspondance with commutative Frobenius algebras [Ko]. In particular the Frobenius algebra $R$ induces a
monoidal functor $\Phi$ from ${\rde C}$ to the category of $K$-modules. This functor has the following properties: it sends $\emptyset$ to $K$ and $S^1$ to $R$, unit,
counit, product and coproduct of $R$ are the image under $\Phi$ of suitable cobordisms and, for every oriented closed surface $\Sigma$, $\Phi(\Sigma)$ is an element of
$K$. An easy computation gives the following:
\vskip 12pt
\noi{\bf 1.6 Lemma: } {\sl For every $p\geq0$ denote by $\Sigma_p$ an oriented surface of genus $p$. Then we have the
following:
$$\Phi(\Sigma_p)=\e(\delta^p)$$
$$\sum_{p\geq0} x^p \Phi(\Sigma_p)=\e({1\over1-x\delta})={t+x(2-t^2s)\over 1-xts+x^2s}\in K[[x]]$$}
\vskip 12pt
\noi{\bf Remark:} If $R$ is the ring $R_0$ and $S$ is a closed oriented surface, $\Phi(S)$ is an element of $\Z[t,s]$ of degree $\chi(S)$ where the degree in
$\Z[t,s]$ is defined by: $\partial^\circ t=2$, $\partial^\circ s=-4$.
\vskip 12pt
\noi{\bf Remark:} The Khovanov homology is defined by using a Frobenius algebra $R$. The classical Khovanov homology corresponds to the case: 
$R=\Z[\alpha]/(\alpha^2)$ and the Lee version of the Khovanov homology corresponds to: $R=\Z[\alpha]/(\alpha^2-1)$. In both cases, we have:
$$\omega=1\hskip 24pt t=0\hskip 24pt \theta=1\hskip 24pt \delta=2\alpha\hskip 24pt s=-4\alpha^2$$
\vskip 12pt
\noi{\bf 1.7 The category of mixed cobordisms ${\rde C}'$.} The functor $\Phi$ is defined on the category ${\rde C}$, but it is possible to extend it to a bigger 
category ${\rde C}'$.

First of all, if $S$ is a surface, consider the commutative monoid with unit defined by generators and relations. Generators are pairs $(x,a)\in S\times R$ and the 
relations are the following:
$$(x,a)(x,b)\equiv(x,ab)$$
$$(x,1)\equiv 1$$
An element of this monoid will be called a $R$-marking of $S$.

If $C$ and $C'$ are two closed oriented curves, define a mixed cobordism from $C$ to $C'$ as a triple $(S,\Gamma,u)$ where: 

--- $S$ is a compact surface containing a closed curve $\Gamma$ in its interior

--- $S\setminus \Gamma$ is oriented with boundary: $\partial(S\minus\Gamma)=C'-C$

--- when crossing $\Gamma$ the orientation of $S\minus \Gamma$ is changed

--- $u$ is a $R$-marking of $S\minus \Gamma$.

The category ${\rde C}'$ is defined as follows: the objects of ${\rde C}'$ are the objects of ${\rde C}$, that is the oriented curves. A morphism in ${\rde C}'$ from an
oriented curve $C$ to an oriented curve $C'$ is the isomorphism class of a mixed cobordism $(S,\Gamma,u)$ from $C$ to $C'$. So we get a monoidal category ${\rde C}'$
containing ${\rde C}$.
\vskip 12pt
\noi{\bf 1.8 Proposition: } {\sl There is a unique monoidal functor $\Phi'$ from ${\rde C}'$ to the category of $K$-modules satisfying the following:

1) $\Phi'$ is an extension of the functor $\Phi$

2) $\Phi'(S,\Gamma,(x_1,a_1)(x_2,a_2)\dots(x_p,a_p))$ depends only on the isotopy classes of the $x_i$'s in $S\minus\Gamma$

3) $\Phi'(S,\Gamma,(x,a)u)$ is $K$-linear with respect to $a$.

4) if $(S,\Gamma,(x,a)u)$ is a mixed cobordism and $x'$ is obtained by making $x$ go through $\Gamma$, we have: 
$\Phi'(S,\Gamma,(x,a)u)=\Phi'(S,\Gamma,(x',\overline a)u)$

5) $\Phi'(S^1\times[0,1],\emptyset,(x,a))$ is the multiplication by $a$, from $R=\Phi'(S^1)$ to $R$

6) $\Phi'(S^1\times[0,1],S^1\times\{1/2\},1)$ is the map $a\mapsto \delta\overline a$

7) $\Phi'(S,\Gamma,u)$ vanishes if $S$ is not orientable.}
\vskip 12pt
\noi{\bf Proof: } Let $\Phi'$ and $\Phi''$ be two functors satisfying all these conditions. Because of conditions 1) and 5), $\Phi'$ and $\Phi''$ are the same
on mixed cobordisms $(S,\emptyset,u)$. Consider now a mixed cobordism $(S,\Gamma,u)$ from $C$ et $C'$. Let $N$ be a small regular neighborough of $\Gamma$ and $S'$ be
the closure of $S\minus N$. Then the morphism $(S,\Gamma,u)$ is the composite of $(S',\emptyset,u)$ from $C$ to $C'\coprod \partial N$ and $1\times(N,\Gamma,1)$
from  $C'\coprod \partial N$ to $C'$. Therefore, in order to prove that $\Phi'$ and $\Phi''$ are the same on $(S,\Gamma,u)$ it is enough to prove that $\Phi'$ and 
$\Phi''$ are the same on $(N,\Gamma,1)$. But that's a consequence of conditions 6) and 7). So, if the functor $\Phi'$ exists, it is unique.

For the construction of $\Phi'$, it is enough to consider the universal case: $R=R_0$. Consider a mixed cobordism $(S,\emptyset,u)$ from $C_0$ to $C_1$. This cobordism
is a composite of cobordisms on the form $(S',\emptyset,1)$ or $(C\times[0,1],\emptyset,u)$. Because of conditions 1) and 5), there is a unique choice
for the morphism $\Phi'(S,\emptyset,u)$. Moreover properties of $R$ and $\Phi$ imply that this morphism depends only on the cobordism $(S,\emptyset,u)$.

Consider now any mixed cobordism $(S,\Gamma,u)$ from $C$ et $C'$. As before denote by $N$ a small regular neighborough of $\Gamma$ and $S'$ be
the closure of $S\minus N$. Then the morphism $(S,\Gamma,u)$ is the composite of $(S',\emptyset,u)$ from $C$ to $C'\coprod \partial N$ and $1\times(N,\Gamma,1)$
from  $C'\coprod \partial N$ to $C'$. So to define $\Phi'(S,\Gamma,u)$, it's enough to define $\Phi'(N,\Gamma,1):\Phi'(\partial N)\rightarrow K$. Let $\Gamma_i$ be the 
component of $\Gamma$ and $N_i$ be the component of $N$ containing $\Gamma_i$. We must set:
$$\Phi'(N,\Gamma,1)=\otimes_i\Phi'(N_i,\Gamma_i,1)$$
If $N_i$ is a M\" obius band, we have to set: $\Phi'(N_i,\Gamma_i,1)=0$. If not, $N_i$ is a band $\Gamma_i\times[-1,1]$ and we define the map $\Phi'(N_i,\Gamma_i,1)$
from $R\otimes R$ to $K$ by: $a\otimes b\mapsto a\overline b+b\overline a$. 

So we get a monoidal functor from ${\rde C}'$ to the category of $K$-modules. It is not difficult to check all the conditions except the last one. 

Suppose that
$(S,\Gamma,u)$ is a mixed cobordism from $C$ to $C'$, where $S$ is nonorientable. Since $S$ isn't orientable, there exists a loop $\gamma$ starting at some point
$x\in S\minus\Gamma$ and intersecting transversally $\Gamma$ an odd number of times. If $u$ is a product of $(x_i,a_i)$, we may suppose that $\gamma$ doesn't meet any of
the $x_i$'s. Let $D$ be a small disk in $S\minus(\Gamma\cup\gamma)$ near $x$ and $S'$ be the closure of $S\minus D$.
Then $(S,\Gamma,u)$ is the composite of $(S',\Gamma,u)$ from $C$ to $C'\coprod S^1$ and $1\times(D,\emptyset,1)$ from $C'\coprod S^1$ to $C'$ and it is enough to
prove that $\Phi'(S',\Gamma,u)$ vanishes.

Because of conditions 3) and 4) we have:
$$\Phi'(S',\Gamma,(x,\alpha)u)=\Phi'(S',\Gamma,(x,\overline\alpha)u)\ \ \ \Longrightarrow\ \ \ \Phi'(S',\Gamma,(x,\alpha-\overline\alpha)u)=0$$
But $\Phi'(S',\Gamma,u)$ is a morphism from $\Phi'(C)$ to $\Phi'(C'\coprod S^1)=\Phi'(C')\otimes R$. So we have:
$$\bigl(1\otimes(\alpha-\overline\alpha)\bigr)\Phi'(S',\Gamma,u)(v)=0$$
for any $v\in \Phi'(C)$. Since $\alpha-\overline\alpha$ is not a zero divisor in $R=R_0$, $\Phi'(S',\Gamma,u)(v)$ vanishes for every
$v\in\Phi'(C)$ and $\Phi'(S',\Gamma,u)$ is the zero map.\cqfd
\vskip 24pt
\noi{\bf 2. Khovanov complexes of diagrams.}
\vskip 24pt
In this section we'll construct many Khovanov complexes associated to link diagrams (oriented or not). All these complexes are graded (or bigraded) differential
$K$-modules (or $R$-modules).
\vskip 12pt
\noi{\bf 2.1 Notations: } Let $X$ be a finite set. Denote by $\Lambda(X)$ the maximal exterior power of the $\Z$-module freely generated by $X$. This module is a
free $\Z$-module of rank $1$.

Suppose $X$ is a graded set. The grading if defined by a map $e$ from $X$ to $\Z$ and we get a $\Z$-grading on $\Lambda(X)$ by the rule:
$$\partial^\circ(x_1\vc x_2\vc \dots\vc x_p)=\sum_ie(x_i)$$
This graded module will be denoted by $\Lambda^e(X)$. If $e$ is the map $x\mapsto 1$ (resp. $x\mapsto -1$), $\Lambda^e(X)$ will be denoted by $\Lambda^+(X)$ (resp.
$\Lambda^-(X))$.
\vskip 12pt
\noi{\bf 2.2 The Khovanov complex $kh(D)$.}
\vskip 12pt
Let $D$ be a link diagram and $X$ be the set of crossings of $D$. Denote by $\widehat X$ the set of maps $s: X\rightarrow \{\pm1\}$. Such a map is called a state on $D$.

If $s\in\widehat X$ is a state, we can modify $D$ near each crossing $x$ by the rule:
$$\begin{tikzpicture}[scale=1/10] \draw (0,0) -- (10,10) ; \draw (10,0) -- (5.4,4.6) ; \draw (4.6,5.4) -- (0,10) ; \draw (20,5) node {$\mapsto$} ;
\draw (30,0)--(33,3);\draw (37,7)--(40,10);\draw (40,0)--(37,3);\draw (30,10)--(33,7);\draw (33,3) arc (-45:45:2.828);\draw (37,7) arc (135:225:2.828);
\draw (70,5) node {if $s(x)=1\ \ $} ; \end{tikzpicture}$$
$$\begin{tikzpicture}[scale=1/10] \draw (0,0) -- (10,10) ; \draw (10,0) -- (5.4,4.6) ; \draw (4.6,5.4) -- (0,10) ; \draw (20,5) node {$\mapsto$} ;
\draw (30,0)--(33,3);\draw (37,7)--(40,10);\draw (40,0)--(37,3);\draw (30,10)--(33,7);\draw (33,3) arc (135:45:2.828);\draw (37,7) arc (-45:-135:2.828);
\draw (70,5) node {if $s(x)=-1$} ;\end{tikzpicture}$$
So we get a new diagram $D_s$ called the $s$-resolution of $D$. This diagram is a curve embedded in the plane and there is a unique compact $K\subset\R^2$ such 
that $D_s$ is the boundary of $K$. Since $K$ is oriented by the plane, $D_s$ is an oriented curve. Denote by $X_s$ the set $s^{-1}(-1)$. So for every state $s$ we have an
oriented curve $D_s$ and a set $X_s$. 

Consider the following graded $K$-module:
$$E=\build\oplus_s^{}\Lambda^-(X_s)\otimes\Phi(D_s)$$

Let $x$ be a crossing of $D$. By sending $x$ to $-1$ and all the other crossings to $1$ we get a state $e_x$. If $s$ is a state denote by $s'$ the state $e_x s$. The
manifold $D_{s'}$ is obtained from $D_s$ by a surgery along a curve $\gamma_x$ near $x$. So we get an oriented cobordism from $D_s$ to $D_{s'}$ and, via the functor
$\Phi$, a map from $\Phi(D_s)$ to $\Phi(D_{s'})$ still denoted by $\gamma_x$. Using this map we have a map $x\otimes\gamma_x$ from $\Lambda^-(X_s)\otimes\Phi(D_s)$ to 
$\Lambda^-(X_{s'})\otimes\Phi(D_{s'})$ defined by:
$$u\otimes v\mapsto x\vc u\otimes\gamma_x(v)$$
Notice that this map is trivial on $\Lambda^-(X_s)\otimes\Phi(D_s)$ if $x$ belongs to $X_s$. 

It is easy to see that the map $d=\sum_x x\otimes\gamma_x$ is a differential on $E$ of degree $-1$. So we get a complex $(E,d)$ denoted by $kh(D)$ (or $kh(D,R)$).
\vskip 12pt
\noi{\bf 2.3 Remark:} If $R=\Z[\alpha]/(\alpha^2)$ and $D$ is oriented, the classical Khovanov complex of $D$ is essentially isomorphic to some suspension of $kh(D)$.
\vskip 12pt
\noi{\bf 2.4 The operators $T_p$.}
\vskip 12pt
If $D$ be a link diagram, a point in $D$ which is not a crossing will be called a regular point in $D$.
\vskip 12pt
Let $D$ be a link diagram and $p$ be a regular point in $D$. Let $a$ be an element of $R$. If $s$ is a state on $D$, there is a unique component $C$ of $D_s$ 
containing $p$. Denote by $D'$ the complement of $C$ in $D_s$. The multiplication by $a\otimes1$ in $\Phi(D_s)=\Phi(C)\otimes\Phi(D')=R\otimes\Phi(D')$
is an endomorphism $f$ of $\Phi(D_s)$. So the map:
$$1\otimes f: u\otimes v\mapsto u\otimes f(v)$$
is an endomorphism of $kh(D)$. This endomorphism will be denoted by $T_p(a)$.
\vskip 12pt
\noi{\bf 2.5 Proposition: } {\sl Let $D$ be a link diagram. Then, for every regular point $p\in D$ and every $a\in R$, 
$T_p(a)$ is a morphism of complexes of degree $0$ from $kh(D)$ to itself. Moreover these operators satisfy the following properties:

--- The operators $T_p(a)$ commute.

--- The map $a\mapsto T_p(a)$ is a $K$-algebra homomorphism from $R$ to End$(kh(D))$.

--- Let $p$ and $q$ be two regular points in $D$. Suppose that these points are the endpoints of a path in $D$ going through exactly one crossing. Then, for every
$a\in R$ the two operators $T_p(a)$ and $T_q(\overline a)$ are homotopic.}
\vskip 12pt
\noi{\bf Proof: } All these properties are easy to check except the last one.

Denote by $x$ the crossing between $p$ and $q$.
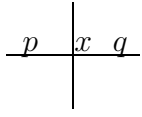
$$\begin{picture}(50,40) \put(0,20){\line(1,0){50}}\put(25,0){\line(0,1){40}}\zput(9,24){$p$}\zput(43,24){$q$}\zput(29,24){$x$}\end{picture}$$
Denote by $D^+$ (resp. $D^-$) the diagram obtained from $D$ by a positive (resp. negative) resolution at $x$. By setting $U=kh(D^+)$ and $V=kh(D^-)$, we have:
$$kh(D)=1\otimes U\oplus x\otimes V$$
Denote by $\gamma$ the surgery homomorphism corresponding to the surgery along the path $\gamma_x$. It is a morphism from $U$ to $V$ and from $V$ to $U$. Denote by $k$
the map from $kh(D)$ to itself defined by:
$$k(1\otimes u)=0\hskip 48pt k(x\otimes v)=1\otimes\gamma(v)$$
We can check that the corresponding homotopy $d(k)=d\circ k+k\circ d$ is the map $1\otimes\gamma^2$.

Set: $T=T_p\otimes T_q$. This operator is a map from $R\otimes R$ to the algebra of endomorphisms of $kh(D)$. It is easy to check the following:
$$d(k)=T(\Delta(1))=T_p(\omega)T_q(\alpha)-T_p(\omega\overline\alpha)=T_q(\omega)T_p(\alpha)-T_q(\omega\overline\alpha)=T_q(\omega)(T_p(\alpha)-T_q(\overline\alpha))$$
and $T_p(\alpha)$ is homotopic to $T_q(\overline\alpha)$ because $\omega$ is invertible. Let $a=u+v\alpha$ be any element in $R$ with $u$ and $v$ in $K$. If $\sim$ is the
homotopy relation, we have:
$$T_p(a)=u+vT_p(\alpha)\sim u+vT_q(\overline\alpha)=T_q(\overline a)\eqno{\cqfd}$$
\vskip 12pt
\noi{\bf 2.6 The Khovanov complex $kh(D,p)$.}
\vskip 12pt
Let $D$ be a link diagram and $p$ be a regular point in $D$. Such a pair $(D,p)$ will be called a pointed diagram. The operator $T_p$ induces an action of $R$ on the 
complex $kh(D)$. Using this action $kh(D)$ becomes a graded differential $R$-module denoted by $kh(D,p)$ (or $kh(D,p,R)$). It is easy to see that $kh(D,p)$ is free over
$R$.

The algebra $R$ may be big: the transcendance degree of $R_0$ is $4$. Nevertheless the complex $kh(D,p)$ can be reduced to a smaller complex.

Consider a graded commutative ring $\Lambda$. Denote by ${\rde M}_{**}(\Lambda)$ the class of bigraded differential free $\Lambda$-modules $C$ satisfying the following:
$$\partial^\circ a=n,\ \ \partial^\circ u=(p,q)\ \ \Longrightarrow\ \ \partial^\circ(au)=(p,q+n),\ \ \partial^\circ(du)=(p-1,q-1)$$
for every $a\in \Lambda$ and $u\in C$. The first component of this degree is called the homological degree and the second the $q$-degree.

Consider the ring $A=\Z[t,\delta,(1-t\delta)^{-1}]$. This ring is graded by the rule:
$$\partial^\circ t=2\hskip 24pt \partial^\circ \delta=-2$$
Moreover $R$ is an $A$-algebra.
\vskip 12pt
\noi{\bf 2.7 Proposition: } {\sl There is a correspondance associating to every pointed diagram $(D,p)$ a complex $kh'(D,p)\in{\rde M}_{**}(A)$ and
an isomorphism:
$$kh(D,p)\build\longrightarrow_{}^\sim R\build\otimes_A^{}kh'(D,p)$$
compatible with the degree in $kh(D,p)$ and the homological degree in $kh'(D,p)$.}
\vskip 12pt
\noi{\bf Proof: } Let $(D,p)$ be a pointed diagram and $X$ be the set of crossings of $D$. For any state $s: X\rightarrow\{\pm1\}$, denote by $C_s$ the set of 
connected components of $D_s$. This set is pointed by the component $c_0$ of $D_s$ containing $p$. Denote by $U_s$ the set of maps $\lambda: C_s\rightarrow\{\pm1\}$
sending $c_0$ to $1$. Finally denote by $U$ the set of pairs $(s,\lambda)$, where $s$ is a state and $\lambda$ is an element of $U_s$.

Because of the universality of $R_0$, we may as well suppose that $R$ is the algebra $R_0$. In this case, we have a degree in $R$:
$$\partial^\circ\alpha=\partial^\circ\overline\alpha=-2\hskip 24pt\partial^\circ\omega=\partial^\circ\overline\omega=0\ \ \Longrightarrow \ \ \partial^\circ t=2
\ \ \ \partial^\circ \delta=-2$$

For every state $s$, set:
$$N_s=\build\oplus_{\lambda\in U_s}^{}R e(s,\lambda)\hskip 48pt M_s=\Lambda^-(X_s)\otimes N_s$$
where the elements $e(s,\lambda)$ are formal vectors in one to one correspondance to the elements of $U$. We set also:
$$M=\build\oplus_s^{} M_s$$
We put a degree on $N_s$ by:
$$\partial^\circ a e(s,\lambda)=\partial^\circ a+\sum_{c\in C_s} \lambda(c)$$
and a bidegree on $M$ by:
$$\partial^\circ(u\otimes v)=(\partial^\circ u,\partial^\circ v)$$

Let $s$ be a state. Take a numbering of $C_s$: $C_s=\{c_0,c_1,\dots,c_{q-1}\}$ in such a way that $c_0$ contains the point $p$. Every $\lambda\in U_s$ is on the form:
$c_i\mapsto \lambda_i$ with $\lambda_0=1$.

Consider elements $a_i\in R$ for $0\leq i<q$. For every $\lambda\in U_s$ we set:
$$a_i(\lambda)=\left\{\matrix{a_i&\hbox{if}\ \lambda_i=1\cr \e(a_i)&\hbox{if}\ \lambda_i=-1\cr}\right.$$
Then we get a map $\varphi_s$ from $\Phi(D_s)=R^{\otimes C_s}$ to $N_s$ defined by:
$$\build\otimes_i^{} a_i\mapsto \sum_{\lambda\in U_s}\Bigl(\prod_i a_i(\lambda)\Bigr)e(s,\lambda)$$
It is easy to check that this map is $R$-linear and bijective. Its inverse is given by:
$$\varphi_s^{-1}(e(s,\lambda))=\prod_{0\leq i<q} b_i$$
with:
$$b_i=\left\{\matrix{\omega_i\omega_0^{-1}&\hbox{if}\ \lambda_i=1\cr \omega_i(\alpha_i-\alpha_0)&\hbox{if}\ \lambda_i=-1\cr}\right.$$
and: $u_i=1^{\otimes i}\otimes u\otimes 1^{\otimes(q-i-1)}$ for $u=\omega$ or $u=\alpha$.

The isomorphisms $\varphi_s$ induce an isomorphism $\varphi$ from $kh(D,p)$ to $M$. Via this isomorphism, the differential on $kh(D,p)$ induces a differential $d'$
on $M$. A straightforward computation shows that the bidegree of $d'$ is $(-1,-1)$ and the entries of the matrix associated to $d'$ are in 
$\{0,\pm1,\pm\delta,\pm\theta^{-1},\pm\theta^{-1}t,\pm\theta^{-1}\delta\}$. That implies the result with:
$$kh'(D,p)=\build\oplus_s^{} \Lambda^-(X_s)\otimes\Bigl(\build\oplus_{\lambda\in U_s}^{}A e(s,\lambda)\Bigr)\eqno{\cqfd}$$
\vskip 12pt
Using technics in [Kh4] it is possible to get a stronger reduction. Let $\beta$ be the element $\alpha-\overline\alpha=\delta/\omega\in R_0$. Every Frobenius
algebra of rank $2$ is a $\Z[\beta]$-algebra. Moreover the map $t\mapsto 0$ and $\delta\mapsto \beta$ induces a ring homomorphism $A\rightarrow\Z[\beta]$ and $\Z[\beta]$
is an $A$-algebra.
\vskip 12pt
\noi{\bf 2.8 Proposition: } {\sl For each pointed diagram $(D,p)$, denote by $kh''(D,p)$ the complex $\Z[\beta]\otimes_A kh'(D,p)\in{\rde M}_{**}(\Z[\beta])$. 

Then, for each pointed diagram $(D,p)$, there is an isomorphism of complexes:
$$kh(D,p)\build\longrightarrow_{}^\sim R\build\otimes_{\Z[\beta]}^{}kh''(D,p)$$
compatible with the degree in $kh(D,p)$ and the homological degree in $kh''(D,p)$.}
\vskip 12pt
\noi{\bf Proof: } Let $R'$ be the algebra $R$ equipped with the following coproduct and counit:
$$\Delta'(u)=\Delta(u/\omega)\hskip 24pt \e'(u)=\e(\omega u)$$
It is easy to see that $R'$ is a Frobenius algebra with generator $\alpha$ and twisting element $1$. Then $R'$ is an $A$-algebra and $t$ and $\delta$ in $A$ are sent
to $0$ and $\beta$ in $R$. So the $A$-algebra structure of $R'$ is actually an $\Z[\beta]$-algebra structure.

Let $(D,p)$ be a pointed diagram. For every state $s$, $D_s$ is oriented. In particular $D_1$ (corresponding to the state $x\mapsto 1$) is the oriented boundary of a
unique compact $K$ in the plane. Consider a point $q$ near $p$ in the interior of $K$. 

Let $s$ be a state. If $c$ is a component of $D_s$, $c$ is the oriented boundary of a unique compact $K_s(c)$
in $S^2$. So we define an integer $f(s,c)$ by:
$$f(s,c)=\left\{\matrix{\chi(K_s(c)\cap K)+1-\chi(K)&\hbox{if}\ q\in K_s(c)\cr \chi(K_s(c)\cap K)&\hbox{otherwise}\cr}\right.$$
where $\chi$ is the Euler characteristic.

It is easy to see the following:

Let $s$ be a state, and $x$ be a crossing of $D$ with $s(x)=1$. Denote by $s'$ the state $s$ modified at $x$ ($s'=s e_x$). If the surgery operator $\gamma_x$ connects 
two components $c$ and $c'$ of $D_s$ into a component $c''$ of $D_{s'}$, we have:
$$f(s',c'')=f(s,c)+f(s,c')$$
If the surgery operator disconnects a component $c$ of $D_s$ into two components $c'$ and $c''$ in $D_{s'}$, we have:
$$f(s,c)+1=f(s',c')+f(s',c'')$$

So, following [Kh4], we get an isomorphism between the complexes $kh(D,p,R)$ and $kh(D,p,R')$. This isomorphism is defined as follows:
$$u\otimes\Bigl(\build\otimes_c^{} v_c\Bigr)\mapsto u\otimes\Bigl(\build\otimes_c^{}{v_c\over\omega^{f(s,c)}}\Bigr)$$
for $u\in\Lambda^-(X_s)$ and $v_c\in R$. Thus we have:
$$kh(D,p,R)\simeq kh(D,p,R')\simeq R'\build\otimes_A^{}kh'(D,p)$$
$$\simeq  R'\build\otimes_{\Z[\beta]}^{}\Z[\beta]\build\otimes_A^{}kh'(D,p)=R'\build\otimes_{\Z[\beta]}^{}kh''(D,p)\eqno{\cqfd}$$
The only problem with this new reduction is the fact that $\beta$ is not necessarily stable under the endomorphisms of $R$.

Another (more serious) problem is the fact that these reductions do not induce any canonical reduction for $kh(D)$. 

\vskip 12pt
\noi{\bf 2.9 The operators $\widehat T_p$.}
\vskip 12pt
If a link diagram $D$ is oriented, it is possible to modify the operators $T_p$ in the following way:

Suppose $p$ is a regular point in $D$. Let $q$ be a point in a neighborough of $p$ and on the left hand side of $D$ and $n$ be the winding number of $D$ about $q$.
For every $a\in R$, define $\widehat T_p(a)$ as the operator $T_p(a)$ if $n$ is even and the operator $T_p(\overline a)$ if $n$ is odd. So $\widehat T_p$ acts on
$kh(D)$. It is easy to see that all these operators commute and that the homotopy class of $\widehat T_p(a)$ depends only on $a$ and the component of $D$ containing $p$.
\vskip 12pt
\noi{\bf 2.10 Khovanov complexes of oriented diagrams.} 
\vskip 12pt
Let $D$ be an oriented link diagram. Since $D$ is oriented, each crossing of $D$ has a sign. Denote by $X_-$ the set of negative crossings of $D$. So we define:
$$KH(D)=kh(D)\otimes\Lambda^+(X_-)$$
This is a graded differential $K$-module. If needed this complex will be denoted by $KH(D,R)$.

The main result of this paper is to prove that the correspondance $D\mapsto KH(D)$ comes from a monoidal functor from the category of cobordisms
of oriented links in $\R^3$ to the homotopy category of $K$-complexes. 

If $R=\Z[\alpha]/\alpha^2$, $KH(D)$ is essentially isomorphic to the classical Khovanov complex. But the isomorphism between these two complexes is not canonical. It is
canonical only up to sign.

If $p$ is a regular point in $D$, the operators $\widehat T_p(a)$ induce a structure of $R$-complex on $KH(D)$. This complex will be denoted by $KH(D,p)$ (or 
$KH(D,p,R)$).

We have also the complex $KH'(D,p)$: In the case: $\widehat T_p(a)=T_p(a)$ the complex $KH'(D,p)$ is the complex $kh'(D,p)\otimes\Lambda^+(X_-)$. In the case: 
$\widehat T_p(a)=T_p(\overline a)$, $KH'(D,p)$ is the complex $kh'(D',p)\otimes\Lambda^+(X_-)$ where the $D'$ is the 
diagram $D$ with the opposite orientation. In any case $KH'(D,p)$ is a complex in ${\rde M}_{**}(A)$ and $KH(D,p)$ is isomorphic to $R\otimes_A KH'(D,p)$.

As before we have also the complex $KH''(D,p)=\Z[\beta]\otimes_A KH'(D,p)$ in ${\rde M}_{**}(\Z[\beta])$.
\vskip 12pt
\noi{\bf 2.11 Proposition: } {\sl Let $D$ be an oriented link diagram. Consider a circle $C$ contained in a half plane disjoint from $D$ and oriented clockwise. Denote
by $(D^\circ,p)$ the union $D\cup C$ pointed by a point $p$ in $C$. Then we have canonical isomorphisms of $R$-complexes:
$$R\build\otimes_K^{} KH(D)\simeq R\build\otimes_A^{} KH'(D^\circ,p)\simeq R\build\otimes_{\Z[\beta]}^{} KH''(D^\circ,p)$$}
\vskip 12pt
\noi{\bf Proof:} The $R$-complex $R\otimes KH(D)$ is obviously isomorphic to $KH(D^\circ,p,R)$. The result follows.\cqfd
\vskip 12pt
\noi{\bf 2.12 The modules $E(D)$ and  $E(D,p)$}
\vskip 12pt
Consider an oriented link diagram $D$. Denote by $C$ the set of components of $D$, by $X$ the set of crossings of $D$ and by $e: X\rightarrow\{\pm\}$ the map sending
each crossing to its sign. For each map $\sigma:C\rightarrow\{\pm\}$ and each sign $e$ denote by $D^e_\sigma$ the subdiagram of $D$ where $\sigma$ is equal to $e$.
Denote also by $D(\sigma)$ the diagram $D$ where the orientation on $D^-_\sigma$ is changed and by $Y(\sigma)\subset X$ the set of crossings between $D^+_\sigma$
and $D^-_\sigma$.

For any ring $B$ denote by $E(D,B)$ the following graded $B$-module:
$$E(D,B)=\build\oplus_{\sigma:C\rightarrow\{\pm\}}^{}\Lambda^{-e}(Y(\sigma))\otimes B v(\sigma)$$
where the $v(\sigma)$'s are formal vectors in one to one correspondance with the $\sigma$'s.

Let $p$ be a regular point in $D$. Let $e_0$ be the sign with is equal to $+$ if and only if the operators $T_p$ and $\widehat T_p$ are the same on $KH(D)$.
The set $C$ is pointed by the component $c_0$ containing $p$. Denote by $\widehat C$ the set of maps $\sigma: C\rightarrow \{\pm\}$ sending $c_0$ to $e_0$.
So, for any ring $B$ denote by $E(D,p,B)$ the following graded $B$-module:
$$E(D,p,B)=\build\oplus_{\sigma\in\widehat C}^{}\Lambda^{-e}(Y(\sigma))\otimes B v(\sigma)$$

If $B$ is a graded ring, $E(D,B)$ and $E(D,p,B)$ is bigraded by the rule:
$$\partial^\circ (u\otimes bv(\sigma))=(\partial^\circ u,\partial^\circ b+b_\sigma)$$
where $b_\sigma$ the number of components of the oriented resolution of $D(\sigma)$. By equipping $E(D,B)$ and $E(D,p,B)$ with the zero differential, these modules
become $B$-complexes or complexes in ${\rde M}_{**}(B)$.
\vskip 12pt
\noi{\bf 2.13 Lemma: } {\sl Let $(D,p)$ be a pointed oriented diagram. Then there are caconical morphisms of complexes:
$$\varphi(R): KH(D,p,R)\longrightarrow E(D,p,R)$$
$$\varphi': KH'(D,p)\longrightarrow E(D,p,A)$$
$$\varphi'': KH''(D,p)\longrightarrow E(D,p,\Z[\beta])$$
of degree $0$ or $(0,0)$.}
\vskip 12pt
\noi{\bf Proof:} Let $\sigma$ be an element of $\widehat C$. The oriented resolution of $D(\sigma)$ (i.e. the only resolution compatible with the orientation) is the
diagram $D_s$ where $s$ is the state $x\mapsto e(x)\sigma(c)\sigma(c')$ and $c$ and $c'$ the components of $D$ containing $x$.

Let $Z(\sigma)$ be the complement of $Y(\sigma)$ in $X$. For each subset $H$ of $X$ we denote by $H_+$ (resp. $H_-$) the set of positive crossings (resp. negative
crossings) in $H$. So we have:
$$X_s=Y(\sigma)_+\cup Z(\sigma)_-\hskip 24pt X_-=Y(\sigma)_-\cup Z(\sigma)_-$$

For any $u\in R$ and any sign $e$ define the element $u^{(e)}$ by:
$$u^{(+)}=u\hskip 24pt u^{(-)}=\overline u$$

The set $C_s$ of components of $D_s$ has $q=b_\sigma$ elements. Take a numbering of $C_s$: $C_s=(c_0(s),c_1(s),\dots,c_{q-1}(s))$, such that $c_0(s)$ contains the point
$p$. For each $i<q$, denote by $d_i$ the winding number of $D_s$ about a point to the left of $c_i(s)$ and by $a_i$ the sign $(-1)^{d_i}$. 

We have a morphism $f_\sigma$ from $KH(D)$ to $\Lambda^{-e}(Y(\sigma))\otimes R$:

This morphism is trivial on $\Lambda^-(X_{s'})\otimes\Phi(D_{s'})\otimes\Lambda^+(X_-)$ for $s'\not=s$. For $s'=s$, we have:
$$\Lambda^-(X_s)\otimes\Phi(D_s)\otimes\Lambda^+(X_-)\simeq \Lambda^-(Y(\sigma)_+)\otimes \Lambda^-(Z(\sigma)_-)\otimes R^{\otimes q}\otimes\Lambda^+(Y(\sigma)_-)\otimes 
\Lambda^+(Z(\sigma)_-)$$
$$\simeq \Lambda^{-e}(Y(\sigma))\otimes R^{\otimes q}$$
and $\Lambda^-(X_s)\otimes\Phi(D_s)\otimes\Lambda^+(X_-)$ is canonically isomorphic to $\Lambda^{-e}(Y(\sigma))\otimes R^{\otimes q}$. So the map $f_\sigma$
is on $\Lambda^-(X_s)\otimes\Phi(D_s)\otimes\Lambda^+(X_-)$ the following:
$$\Lambda^-(X_s)\otimes\Phi(D_s)\otimes\Lambda^+(X_-)\build\longrightarrow_{}^\sim\Lambda^{-e}(Y(\sigma))\otimes R^{\otimes q}\build\longrightarrow_{}^{1\otimes g} 
\Lambda^{-e}(Y(\sigma))\otimes R$$
where $g$ is the map:
$$b_0\otimes b_1\otimes\dots\otimes b_{q-1}\mapsto \prod_i b_i^{(a_i)}$$
More precisely, $f_\sigma$ is the map:
$$u_1\otimes u_2\otimes u_3\otimes u_4\otimes u_5\ \ \mapsto\ \ (-1)^{|u_2||u_4|}<u_2,u_5> u_1\vc u_4\otimes g'(u_3)$$
for every $(u_1,u_2,u_3,u_4,u_5)\in\Lambda^-(Y(\sigma)_+)\times\Lambda^-(Z(\sigma)_-)\times\Phi(D_s)\times\Lambda^+(Y(\sigma)_-)\times\Lambda^+(Z(\sigma)_-)$,
where $|u|$ is the degree of $u$, $g'$ is the map:
$$\Phi(D_s)\build\longrightarrow_^\sim R^{\otimes q}\build\longrightarrow_{}^g R$$
and $<?,?>$ is the isomorphism $\Lambda^-(Z(\sigma)_-)\otimes\Lambda^+(Z(\sigma)_-)\build\longrightarrow_^\sim \Z$ defined by:
$$<x_1\vc x_2\vc\dots\vc x_n,x_1\vc x_2\vc\dots\vc x_n>=(-1)^{n(n-1)/2}$$

It is easy to see that this map is a morphism of degree $0$ from $KH(D)$ to $\Lambda^{-e}(Y(\sigma))\otimes R$ inducing a morphism $\varphi_\sigma(R)$ from $KH(D,p,R)$ to 
$\Lambda^{-e}(Y(\sigma))\otimes Rv(\sigma)$. 

On the level of complexes $KH'$, we get a morphism $\varphi_\sigma$ from $KH'(D,p)$ to $\Lambda^{-e}(Y(\sigma))\otimes Av(\sigma)$. This morphism sends $e(s,\lambda)$ to
$0$ if there exists an $i$ with $a_i=e_0$ and $\lambda_i=-$. If it is not the case, the morphism sends $e(s,\lambda)$ to $\theta^{-c}(-\delta)^d v(D')$, where $d$ is the
number of $i$ such that: $\lambda_i=-$ and $c$ the number of $i$ such that: $a_i=-e_0$. Because of the degree of $v(\sigma)$, the bidegree of $\varphi_\sigma$ is $(0,0)$.

By taking the sum of all these morphisms we get the desired morphisms $\varphi(R)$, $\varphi'$ and then $\varphi''$.\cqfd
\vskip 12pt
\noi{\bf 2.14 Remark:} Suppose $\delta$ is invertible in $R$. Then the map $\varphi(R)$ has a section. By using the notations above this section is a morphism of
complexes defined by:
$$\Lambda^{-e}(Y(\sigma))\otimes Rv(\sigma)\build\longrightarrow_{}^\sim \Lambda^-(X_s)\otimes Rv(\sigma)\otimes\Lambda^+(X_-)\build\longrightarrow_{}^{1\otimes g\otimes 1}
\Lambda^-(X_s)\otimes\Phi(D_s)\otimes\Lambda^+(X_-)$$
where $g$ is the map:
$$u v(\sigma)\mapsto u_0\prod_{0\leq i<j<q}\Bigl(\delta_0^{-1}\bigl(\omega_i^{(a_i)}\alpha_j^{(a_j)}-(\omega_i\overline\alpha_i)^{(a_i)}\bigr)\Bigr)$$
and $u_i$ is, for every $u\in R$, the element $1^{\otimes i}\otimes u\otimes1^{\otimes(q-i-1)}$.
\vskip 12pt
\noi{\bf 2.15 Theorem:} {\sl Let $(D,p)$ be an pointed oriented link diagram. Suppose that the underlying graph of $D$ is connected. Then the morphisms:
$$\varphi(R): KH(D,p,R)\longrightarrow E(D,p,R)$$
$$\varphi': KH'(D,p)\longrightarrow E(D,p,A)$$
$$\varphi'': KH''(D,p)\longrightarrow E(D,p,\Z[\beta])$$
are surjective and the action of $\beta$ is homotopically nilpotent on their kernels.}
\vskip 12pt
\noi{\bf Proof:} Because the underlying graph of $D$ is connected, the map sending $\sigma\in\widehat C$ to the state $s: x\mapsto e(x)\sigma(c)\sigma(c')$ is injective.
Therefore all maps $\varphi$ are surjective, in particular $\varphi''$.

The last thing to do is to prove that the action of $\beta$ on the kernel $U$ of $\varphi''$ is homotopically nilpotent.
But that is equivalent to the fact that $\Z[\beta^\pm]\otimes_{\Z[\beta]} U$ is acyclic. Let $B$ be the ring $Z[\beta^\pm]$ and $R_1$ be the algebra $B\times B$.
Using the diagonal map $B\rightarrow B\times B$, $R_1$ is a Frobenius $B$-algebra of rank $2$ with generator 
$(1,0)$. In this algebra, the involution is: $(u,v)\mapsto (v,u)$, the twisting element is $1$ and the counit is: $(u,v)\mapsto u-v$. 

It is clear that $B\otimes U$ is acyclic if and only if $R_1\otimes U$ is acyclic. Therefore it is enough to prove that 
$\varphi(R_1)$ induces an isomorphism in homology from $KH(D,p,R_1)$ to $E(D,p,R_1)$. But in $R_1$ we have two orthogonal idempotents:
$$\pi_+=(1,0)\hskip 24pt \pi_-=(0,1)$$
So we can use the Karoubi completion method of Bar-Natan and Morrison [BM] to prove that $\varphi(R_1)$ is a homology equivalence.\cqfd
\vskip 12pt
\noi{\bf 2.16 Remark:} If the underlying graph of $D$ is not connected the result is still true but the kernels have to be replaced by the homotopy kernels. In 
particular the morphisms $\varphi$ are homotopy equivalences if $\delta$ is invertible in $R$.
\vskip 24pt
\noi{\bf 3. Elementary moves.}
\vskip 24pt
Consider an  elementary move (i.e. a Reidemeister move or a surgery move) $f:D\rightarrow D'$ transforming a link diagram $D$ into a link diagram $D'$. If the 
correspondance $D\mapsto kh(D)$ comes from a functor on the category of cobordisms of links, the move $f$ would have to induce a morphism $f_*$ from
$kh(D)$ to $kh(D')$. So the first thing to do is to associate to every elementary move $f:D\rightarrow D'$ a morphism from $kh(D)$ to $kh(D')$. 
\vskip 12pt
\noi{\bf 3.1 Proposition: } {\sl There is a correspondance associating to every elementary move $f:D\rightarrow D'$ a morphism $f^0: kh(D)\rightarrow kh(D')$ such that:

if $f$ is a Reidemeister move with inverse move $g$, $f^0$ is a homotopy equivalence and $g^0$ is a homotopy inverse of $f$.}
\vskip 12pt
\noi{\bf Proof: } In order to construct the correspondance $f\mapsto f^0$, we have to consider all types of elementary moves.
\vskip 12pt
\noi{\bf Surgery moves:} There is three kinds of surgery moves: surgery moves of index $0$, $1$ or $2$. A surgery move $f:D\rightarrow D'$ of index $0$ transforms $D$
by adding a circle bounding a disk in the plane which is disjoint from $D$ and the inverse move $g:D'\rightarrow D$ is a surgery move of index $2$. 
$$\begin{tikzpicture}[scale=1/8] \draw (5,5) node {$f:$} ; \draw [very thick] (10,0)--(30,0)--(30,10)--(10,10)--(10,0) ; \draw [very thick] (20,0)--(20,10) ;
\draw (25,5) circle (2) ;\end{tikzpicture}$$

Let $D$ be a link diagram. Consider a path $\gamma$ embedded in the plane which doesn't meet any crossing point of $D$ and intersects $D$ in its boundary. Such a path
is called a surgery path of $D$.

A surgery move $f:D\rightarrow D'$ of index $1$ is a modification $D\mapsto D'$ by surgery along some surgery path $\gamma$ of $D$. 
$$\begin{tikzpicture}[scale=1/8] \draw (-5,5) node {$f:$} ; \draw [very thick] (0,0)--(20,0)--(20,10)--(0,10)--(0,0) ; \draw [very thick] (10,0)--(10,10) ;
\draw (2,2) arc (-45:45:4.242);\draw (8,8) arc (135:225:4.242); \draw (18,2) arc (45:135:4.242);\draw (12,8) arc (-135:-45:4.242);
\draw [very thin] (3.242,5)--(6.758,5);\draw (5,4) node {$\scriptstyle\gamma$};\end{tikzpicture}$$
Suppose $f:D\rightarrow D'$ is a surgery move of index $0$ and $g:D' \rightarrow D$ its inverse move. Then we have: $kh(D')=kh(D)\otimes R$ and the map $f^0$ and $g^0$
are defined by:
$$\forall (u,a)\in kh(D)\times R,\ \ f^0(u)=u\otimes 1,\hskip 24pt g^0(u\otimes a)=u\e(a)$$

Suppose $f:D\rightarrow D'$ is a surgery move of index $1$. This surgery move is defined by a surgery path $\gamma$. For every state $s$ the surgery induces, via
the functor $\Phi$, a map $\gamma:\Phi(D_s)\rightarrow \Phi(D'_s)$ and these maps induce a map $f^0$ from $kh(D)$ to $kh(D')$. It is easy to see that $f^0$ is a 
morphism of complexes.

\vskip 12pt
\noi{\bf Reidemeister moves of type I:} There is four kinds of Reidemeister moves of type I depending of two signs. The first sign is $+$ (or $1$) if the move creates a
new crossing and $-$ (or $-1$) if it removes one crossing. The second one is the sign of this crossing (for any orientation of the diagram). We say that $f$ is a 
Reidemeister move of type I$_{e'}^e$ if $f$ is a Reidemeister move where $e$ is the first sign and $e'$ the second one. It is clear that the inverse move of a 
Reidemeister move of type I$_{e'}^e$ is a Reidemeister move of type I$_{e'}^{-e}$.

Consider a Reidemeister move $f: D\rightarrow D'$ of type I$_e^+$ and denote by $g$ its inverse move (of type I$_e^-$). Let $x$ be the created crossing and $U$ be the
complex $kh(D)$. Let $p$ be a point in $D$ which is near $x$. So we have:
$$kh(D')=\left\{\matrix{U\otimes R\oplus x\otimes U&\hbox{if}\ e=+\cr U\oplus x\otimes U\otimes R&\hbox{if}\ e=-\cr}\right.$$

Suppose $e=+$. Then the move is the following:
$$\begin{tikzpicture}[scale=1/8] \draw [very thick] (-7,-3)--(17,-3)--(17,13)--(-7,13)--(-7,-3); \draw [very thick] (5,-3)--(5,13); \draw (-12,5) node {$f:$} ;
\draw (-1,0)--(-1,10); \draw [domain=0:8.2] plot ({12-2*cos(36*\x)},1/13*\x^3-15/13*\x^2+63/13*\x) ; \draw (9,5) node {$x$}; \draw (-3,5) node {$p$};
\draw [domain=8.65:10] plot ({12-2*cos(36*\x)},1/13*\x^3-15/13*\x^2+63/13*\x) ;\end{tikzpicture}$$
and  we set:
$$\forall (u,a)\in U\times R,\ \ f^0(u)=u\otimes \omega\alpha-T_p(\alpha)u\otimes\omega,\hskip 24pt g^0(x\otimes u)=0,\ \ g^0(u\otimes a)=u\e(a)$$
Suppose $e=-$. Then the move is the following:
$$\begin{tikzpicture}[scale=1/8] \draw [very thick] (-7,-3)--(17,-3)--(17,13)--(-7,13)--(-7,-3); \draw [very thick] (5,-3)--(5,13); \draw (-12,5) node {$f:$} ;
\draw (-1,0)--(-1,10); \draw [domain=0:1.3] plot ({12-2*cos(36*\x)},1/13*\x^3-15/13*\x^2+63/13*\x) ; \draw (9,5) node {$x$}; \draw (-3,5) node {$p$};
\draw [domain=1.8:10] plot ({12-2*cos(36*\x)},1/13*\x^3-15/13*\x^2+63/13*\x) ;\end{tikzpicture}$$
and  we set:
$$\forall (u,a)\in U\times R,\ \ f^0(u)=x\otimes u\otimes 1\hskip 24pt g^0(u)=0,\ \ g^0(x\otimes u\otimes a)=T_p(\overline a)u$$
It is not difficult to see that $f^0$ and $g^0$ are morphisms of complexes and that $g^0$ is a left inverse of $f^0$. Moreover the cokernel of $f^0$ is the mapping cone
of an isomorphism. So $f^0$ is a homotopy equivalence and $g^0$ is a homotopy inverse of $f^0$.
\vskip 12pt
\noi{\bf Reidemeister moves of type II:} There is two kinds of Reidemeister moves of type II: Reidemeister moves of type II$^+$ that create two new crossings and 
their inverses of type II$^-$.

Consider a Reidemeister move $f: D\rightarrow D'$ of type II$^+$ and its inverse move $g:D'\rightarrow D$. Let $D''$ be the diagram obtained from $D$ by a surgery along
a path joining the two branches of $D$ modified by $f$. 
$$\begin{tikzpicture}[scale=1/8]  \draw [very thick] (-5,-3)--(23,-3)--(23,13)--(-5,13)--(-5,-3);\draw [very thick] (9,-3)--(9,13);\draw (-2,5) node {$p$};
\draw (0,0)--(0,10);\draw (4,0,)--(4,10); \draw [domain=0:10] plot ({14+3*sin(18*\x)},\x) ; \draw (13.3,7.5) node {$x$};\draw (13.3,2.5) node {$y$};
\draw [domain=0:1.6] plot ({18-3*sin(18*\x)},\x) ;\draw [domain=3:7] plot ({18-3*sin(18*\x)},\x) ;\draw [domain=8.4:10] plot ({18-3*sin(18*\x)},\x) ;
\draw (-10,5) node {$f:$};\draw (40,5) node {$D'':$}; \draw (45,0) arc (180:0:2);\draw (45,10) arc (-180:0:2);
\end{tikzpicture}$$
Denote by $U$ and $V$ the complexes $kh(D)$ and $kh(D'')$. Then we have:
$$kh(D')=V\oplus x\otimes U\oplus y\otimes V\otimes R\oplus x\vc y\otimes V$$
where $x$ and $y$ are the created crossings in $D'$. Surgery moves from $D$ to $D''$ and from $D''$ to $D$ induce morphisms $\gamma:U\rightarrow V$ and
$\gamma: V\rightarrow U$. The morphisms $f^0$ and $g^0$ are defined by:
$$\forall u\in U,\ \ f^0(u)=x\otimes u+y\otimes\gamma(T(\omega^{-1})u)\otimes\omega, \ \ \ g^0(x\otimes u)=u$$
$$\forall (v,a)\in V\times R,\ \ g^0(v)=g^0(x\vc y\otimes v)=0,\ \ \ g^0(y\otimes v\otimes a)=-\e(a)\gamma(v)$$
where $T$ is the operator $T_p$ for some $p$ in one of the two branches of $D$.

We can see that $g^0$ doesn't depend on $p$ and $f^0$ and $g^0$ are morphisms of complexes. Moreover $g^0$ is a left inverse of $f^0$ and the cokernel of $f^0$ is the
mapping cone of an isomorphism. So $f^0$ is a homotopy equivalence and $g^0$ is a homotopy inverse of $f^0$.
\vskip 12pt
\noi{\bf Reidemeister moves of type III:} There is two kinds of Reidemeister moves of type III depending on a sign. If $f:D\rightarrow D'$ is a Reidemeister move of
type III, three crossings of $D$ are modified. These crossings are the vertices of a triangle $\Theta$ which is oriented by the orientation of the plane. We number the 
branches starting with the top branch and ending with the bottom one. So we get a numbering of the edges of the triangle. We say that the move is of type III$_+$ if
this numbering is compatible with the orientation of $\Theta$ and III$_-$ if it is not the case.

Let's denote by $x$ the intersection of the top branch and the bottom branch. The other vertices of $\Theta$ are denoted by $y$ and $z$ in such a way that the numbering 
$(x,y,z)$ is compatible with the orientation of the triangle. We proceed similarly for the diagram $D'$ and the move is the following:
$$\begin{tikzpicture}[scale=1/4] \draw (-2.736,-.5)--(2.736,-.5);\draw (-1.5,-2.5)--(1,2.5);\draw (1.5,-2.5)--(-1,2.5);
\draw (2+3.264,.5)--(2+8.736,.5);\draw (2+7.5,2.5)--(2+5,-2.5);\draw (2+4.5,2.5)--(2+7,-2.5);\draw (-7,0) node {$f:$};
\draw [very thick] (-4,-4)--(12,-4)--(12,4)--(-4,4)--(-4,-4);\draw [very thick] (4,-4)--(4,4);\draw (0,2) node {$x$};\draw (8,-2) node {$x$};
\draw (-1.8,-1.5) node {$y$};\draw (9.8,1.5) node {$y$};\draw (1.8,-1.5) node {$z$};\draw (6.2,1.5) node {$z$};
\end{tikzpicture}$$

By resolution of $D$ near the triangle we get five new diagrams:
$$\begin{tikzpicture}[scale=1/4] \draw (3,0) arc (-90:-210:1.732);\draw (-3,0) arc (-90:30:1.732);\draw (1.5,-2.598) arc (30:150:1.732);\draw (-5.5,0) node {$D_0:$};
\draw (16-.4,0) arc (90:210:1.732);\draw (10-.4,0) arc (90:-30:1.732);\draw (14.5-.4,2.598) arc (-30:-150:1.732);\draw (7.5-.4,0) node {$D_1:$};
\draw (23-.8,0)--(29-.8,0);\draw (27.5-.8,-2.598) arc (30:150:1.732);\draw (27.5-.8,2.598) arc (-30:-150:1.732);\draw (20.5-.8,0) node {$D_x:$};
\draw (40.5-1.2,-2.598)--(37.5-1.2,2.598);\draw (36-1.2,0) arc (90:-30:1.732);\draw (42-1.2,0) arc (-90:-210:1.732);\draw (33.5-1.2,0) node {$D_y:$};
\draw (53.5-1.6,2.598)--(50.5-1.6,-2.598);\draw (49-1.6,0) arc (-90:30:1.732);\draw (55-1.6,0) arc (90:210:1.732);\draw (46.5-1.6,0) node {$D_z:$};
\end{tikzpicture}$$
and, for $D'$, we get: $D'_0=D_1$, $D'_1=D_0$, $D'_x=D_x$, $D'_y=D_y$ and $D'_z=D_z$.

With these new diagrams we get five new complexes: $U=kh(D_0)$, $V=kh(D_1)$, $A=kh(D_x)$, $B=kh(D_y)$, $C=kh(D_z)$. 

Consider six points $p_i, i\in \Z/6$, sitted as follows in $D$:
$$\begin{tikzpicture}[scale=1/4] \draw (-2.736,-.5)--(2.736,-.5);\draw (-1.5,-2.5)--(1,2.5);\draw (1.5,-2.5)--(-1,2.5);
\draw (-3.6,-.5) node {$\scriptstyle p_2$};\draw (3.6,-.5) node {$\scriptstyle p_5$};\draw (-1.2,2.9) node {$\scriptstyle p_1$};
\draw (1.2,2.9) node {$\scriptstyle p_6$};\draw (-1.7,-2.9) node {$\scriptstyle p_3$};\draw (1.7,-2.9) node {$\scriptstyle p_4$};
\end{tikzpicture}$$
All these points are in a circle $\Gamma$ which bounds a disk $\Delta$. Since the move modifies the diagram $D$ only in $\Delta$, the points $p_i$ are also in
diagrams $D'$, $D_0$, $D_1$, $D_x$, $D_y$ and $D_z$.  

The arc of $\Gamma$ with endpoints $p_i$ and $p_{i-1}$ will be denoted by $\gamma_i$. This arc induces a surgery operator still denoted by $\gamma_i$. We'll denote also
by $T_i$ the operator $T_{p_i}$. All these operators are well defined on the complexes $U$, $V$, $A$, $B$ and $C$.

Suppose $f$ is a move of type III$_+$. Then the move is the following:
$$\begin{tikzpicture}[scale=1/4] \draw (-2.736,-.5)--(-.7,-.5);\draw  (-.3,-.5)--(2.736,-.5);\draw (-1.5,-2.5)--(1,2.5);
\draw (1.5,-2.5)--(.65,-.8);\draw (.35,-.2)--(.15,.2);\draw (-.15,.8)--(-1,2.5);
\draw (8+2.736,.5)--(8.7,.5);\draw  (8.3,.5)--(8-2.736,.5);\draw (9.5,2.5)--(7,-2.5);
\draw (8-1.5,2.5)--(8-0.65,.8);\draw (8-0.35,.2)--(8-0.15,-.2);\draw (8.15,-.8)--(9,-2.5);\draw (-7,0) node {$f:$};
\draw [very thick] (-4,-4)--(12,-4)--(12,4)--(-4,4)--(-4,-4);\draw [very thick] (4,-4)--(4,4);\draw (0,2) node {$x$};\draw (8,-2) node {$x$};
\draw (-1.8,-1.5) node {$y$};\draw (9.8,1.5) node {$y$};\draw (1.8,-1.5) node {$z$};\draw (6.2,1.5) node {$z$};
\end{tikzpicture}$$
And we have:
$$kh(D)=V\oplus x\otimes V\otimes R\oplus y\otimes C\oplus z\otimes B\oplus y\vc z\otimes U\oplus x\vc z\otimes V\oplus x\vc y\otimes V\oplus x\vc y\vc z\otimes A$$
$$kh(D')=U\oplus x\otimes U\otimes R\oplus y\otimes C\oplus z\otimes B\oplus y\vc z\otimes V\oplus x\vc z\otimes U\oplus x\vc y\otimes U\oplus x\vc y\vc z\otimes A$$
In $kh(D)$ and $kh(D')$ the differential $d$ is on the form $d=d_1+d_2$, where $d_1$ corresponds to the differentials of the complexes $U,V,A,B,C$. The map $d_2$ is 
defined on $kh(D)$ by the following (for any $r\in R$, $u\in U$, $v\in V$, $a\in A$, $b\in B$, $c\in C$):
$$d_2(v)=x\otimes \gamma_1(v)+y\otimes \gamma_2(v)+z\otimes \gamma_6(v)$$
$$d_2(x\otimes v\otimes r)=-x\vc y\otimes T_3(r)v-x\vc z\otimes T_5(r)v$$
$$d_2(y\otimes c)=x\vc y\otimes\gamma_1(c)-y\vc z\otimes\gamma_4(c)$$
$$d_2(z\otimes b)=x\vc z\otimes\gamma_1(b)+y\vc z\otimes\gamma_4(b)$$
$$d_2(y\vc z\otimes u)=x\vc y\vc z\otimes\gamma_1(u)$$
$$d_2(x\vc z\otimes v)=-x\vc y\vc z\otimes\gamma_4(v)$$
$$d_2(x\vc y\otimes v)=x\vc y\vc z\otimes\gamma_4(v)$$
$$d_2(x\vc y\vc z\otimes a)=0$$
and, by exchanging $U$ and $V$ and replacing $\gamma_i$ and $T_i$ by $\gamma_{i+3}$ and $T_{i+3}$, we get the differential $d_2$ on $kh(D')$.

So we define the map $f^0$ by the following:
$$f^0(v)=0$$
$$f^0(x\otimes v\otimes r)=\e(r\omega)\Bigl(x\otimes T_1(\omega^{-1})\gamma_2\gamma_6(v)\otimes1+y\otimes T_1(\omega^{-1})\gamma_2(v)+z\otimes 
T_1(\omega^{-1})\gamma_6(v)\Bigr)$$
$$f^0(y\otimes c)=-x\otimes\gamma_4(c)\otimes1-y\otimes c$$
$$f^0(z\otimes b)=-z\otimes b$$
$$f^0(x\vc y\otimes v)=y\vc z\otimes v$$
$$f^0(x\vc z\otimes v)=-y\vc z\otimes v$$
$$f^0(y\vc z\otimes u)=-x\vc z\otimes u$$
$$f^0(x\vc y\vc z\otimes a)=x\vc y\vc z\otimes a$$
We can check that $f^0$ is a morphism of complexes of degree $0$. Let $W$ be the $\Z$-module freely generated by $\lambda$, $d(\lambda)$, $\mu$ and $d(\mu)$. By setting:
$$\partial^\circ\lambda=0\hskip 24pt\partial^\circ\mu=-1$$
$W$ becomes an acyclic graded differential $\Z$-module. We have maps $\Psi: W\otimes V\rightarrow kh(D)$ and $\Psi': W\otimes U\rightarrow kh(D')$ defined by:
$$\Psi(\lambda\otimes v)=v\hskip 24pt\Psi(d(\lambda)\otimes v)=d_2(v)$$
$$\Psi(\mu\otimes v)=x\otimes v\otimes1\hskip24pt\Psi(d(\mu)\otimes v)=d_2(x\otimes v\otimes1)$$
$$\Psi'(\lambda\otimes u)=u\hskip 24pt\Psi'(d(\lambda)\otimes u)=d_2(u)$$
$$\Psi'(\mu\otimes u)=x\otimes u\otimes1\hskip24pt\Psi'(d(\mu)\otimes u)=d_2(x\otimes u\otimes1)$$
It is easy to see that $\Psi$ ans $\Psi'$ are injective morphisms of complexes of degree $0$. Since $W$ is acyclic, the images of $\Psi$ and $\Psi'$ are acyclic too. 
Moreover the image of $\Psi$ is killed by $f^0$ and the image of $\Psi'$ is therefore killed by $g^0$, where $g$ is the inverse move of $f$. Modulo these images we have:
$$kh(D)\equiv y\otimes C\oplus z\otimes B\oplus y\vc z\otimes U\oplus x\vc z\otimes V\oplus x\vc y\vc z\otimes A$$
$$kh(D')\equiv y\otimes C\oplus z\otimes B\oplus y\vc z\otimes V\oplus x\vc z\otimes U\oplus x\vc y\vc z\otimes A$$
and maps $f^0$ and $g^0$ are on these quotients:
$$f^0:\ \ \left\{\matrix{y\otimes c\mapsto -y\otimes c\cr z\otimes b\mapsto -z\otimes b\cr y\vc z\otimes u\mapsto -x\vc z\otimes u\cr x\vc z\otimes v\mapsto -y\vc z
\otimes v\cr x\vc y\vc z\otimes a\mapsto x\vc y\vc z\otimes a\cr}\right.\hskip 48pt
g^0:\ \ \left\{\matrix{y\otimes c\mapsto -y\otimes c\cr z\otimes b\mapsto -z\otimes b\cr y\vc z\otimes v\mapsto -x\vc z\otimes v\cr x\vc z\otimes u\mapsto -y\vc z
\otimes u\cr x\vc y\vc z\otimes a\mapsto x\vc y\vc z\otimes a\cr}\right.$$
Therefore $f^0$ and $g^0$ are homotopy equivalences and $g^0$ is a homotopy inverse of $f^0$.

Suppose now the move is of type III$_-$. The move is the following:
$$\begin{tikzpicture}[scale=1/4] \draw (2.736,-.5)--(.7,-.5);\draw  (.3,-.5)--(-2.736,-.5);\draw (1.5,-2.5)--(-1,2.5);
\draw (-1.5,-2.5)--(-.65,-.8);\draw (-.35,-.2)--(-.15,.2);\draw (.15,.8)--(1,2.5);
\draw (8-2.736,.5)--(8-.7,.5);\draw  (8-.3,.5)--(8+2.736,.5);\draw (8-1.5,2.5)--(8+1,-2.5);
\draw (8+1.5,2.5)--(8+0.65,.8);\draw (8+0.35,.2)--(8+0.15,-.2);\draw (8-.15,-.8)--(8-1,-2.5);\draw (-7,0) node {$f:$};
\draw [very thick] (-4,-4)--(12,-4)--(12,4)--(-4,4)--(-4,-4);\draw [very thick] (4,-4)--(4,4);\draw (0,2) node {$x$};\draw (8,-2) node {$x$};
\draw (-1.8,-1.5) node {$y$};\draw (9.8,1.5) node {$y$};\draw (1.8,-1.5) node {$z$};\draw (6.2,1.5) node {$z$};
\end{tikzpicture}$$
And we have:
$$kh(D)=A\oplus x\otimes U\oplus y\otimes V\oplus z\otimes V\oplus y\vc z\otimes V\otimes R\oplus x\vc z\otimes C\oplus x\vc y\otimes B\oplus x\vc y\vc z\otimes V$$
$$kh(D')=A\oplus x\otimes V\oplus y\otimes U\oplus z\otimes U\oplus y\vc z\otimes U\otimes R\oplus x\vc z\otimes C\oplus x\vc y\otimes B\oplus x\vc y\vc z\otimes U$$
As before the differentials of $kh(D)$ and $kh(D')$ are on the form $d=d_1+d_2$, where $d_1$ corresponds to the differentials of the complexes $U,V,A,B,C$. 
The map $d_2$ is defined on $kh(D)$ by the following:
$$d_2(a)=x\otimes \gamma_2(a)+y\otimes \gamma_5(a)+z\otimes \gamma_5(a)$$
$$d_2(x\otimes u)=-x\vc y\otimes \gamma_3(u)-x\vc z\otimes \gamma_5(u)$$
$$d_2(y\otimes v)=x\vc y\otimes \gamma_6(v)-y\vc z\otimes \gamma_5(v)$$
$$d_2(z\otimes v)=x\vc z\otimes \gamma_2(v)+y\vc z\otimes \gamma_3(v)$$
$$d_2(y\vc z\otimes v\otimes r)=x\vc y\vc z\otimes T_1(r)v$$
$$d_2(x\vc y\otimes b)=x\vc y\vc z\otimes \gamma_1(b)$$
$$d_2(x\vc z\otimes c)=-x\vc y\vc z\otimes \gamma_1(c)$$
$$d_2(x\vc y\vc z\otimes v)=0$$
and, as before, we get the differential $d_2$ on $kh(D')$ by exchanging $U$ and $V$ and replacing $\gamma_i$ and $T_i$ by $\gamma_{i+3}$ and $T_{i+3}$.
So we define the map $f^0$ by the following:
$$f^0(a)=a$$
$$f^0(x\otimes u)=(y+z)\otimes u$$
$$f^0(y\otimes v)=x\otimes v$$
$$f^0(z\otimes v)=0$$
$$f^0(y\vc z\otimes v\otimes r)=\e(r)\bigl(x\vc z\otimes\gamma_2(v)+y\vc z\otimes T_4(\omega^{-1})\gamma_4\gamma_2(v)\otimes \omega\bigr)$$
$$f^0(x\vc y\otimes b)=-x\vc y\otimes b+y\vc z\otimes T_4(\omega^{-1})\gamma_4(b)\otimes \omega$$
$$f^0(x\vc z\otimes c)=-x\vc z\otimes c-y\vc z\otimes T_4(\omega^{-1})\gamma_4(c)\otimes \omega$$
$$f^0(x\vc y\vc z\otimes v)=0$$
We can check that $f^0$ is a morphism of complexes of degree $0$. We have maps $\Psi: W\otimes V\rightarrow kh(D)$ and $\Psi': W\otimes U\rightarrow kh(D')$:
$$\Psi(\lambda\otimes v)=z\otimes v\hskip 48pt \Psi(d(\lambda)\otimes v)=d_2(z\otimes v)$$
$$\Psi(\mu\otimes v)=y\vc z\otimes v\otimes\omega\hskip 48pt \Psi(d(\mu)\otimes v)=d_2(y\vc z\otimes v\otimes\omega)$$
$$\Psi'(\lambda\otimes u)=z\otimes u\hskip 48pt \Psi'(d(\lambda)\otimes u)=d_2(z\otimes u)$$
$$\Psi'(\mu\otimes u)=y\vc z\otimes u\otimes\omega\hskip 48pt \Psi'(d(\mu)\otimes u)=d_2(y\vc z\otimes u\otimes\omega)$$
These two maps are injective morphisms of complexes of degree $-1$ and their images are acyclic subcomplexes killed by $f^0$ and $g^0$. Modulo these images we have:
$$kh(D)\equiv A\oplus x\otimes U\oplus y\otimes V\oplus x\vc y\otimes B\oplus x\vc z\otimes C$$
$$kh(D')\equiv A\oplus x\otimes V\oplus y\otimes U\oplus x\vc y\otimes B\oplus x\vc z\otimes C$$
and maps $f^0$ and $g^0$ are:
$$f^0:\ \ \left\{\matrix{a\mapsto a\cr x\otimes u\mapsto y\otimes u\cr y\otimes v\mapsto x\otimes v\cr x\vc y\otimes b\mapsto-x\vc y\otimes b\cr x\vc z\otimes c
\mapsto -x\vc z\otimes c\cr}\right. \hskip 48pt
g^0:\ \ \left\{\matrix{a\mapsto a\cr x\otimes v\mapsto y\otimes v\cr y\otimes u\mapsto x\otimes u\cr x\vc y\otimes b\mapsto-x\vc y\otimes b\cr x\vc z\otimes c
\mapsto -x\vc z\otimes c\cr}\right.$$
Therefore $f^0$ and $g^0$ are homotopy equivalences and $g^0$ is a homotopy inverse of $f^0$.\cqfd
\vskip 12pt
\noi{\bf 3.2 Remark:} For every elementary move $f$, the map $f^0$ is a morphism of degree $0$ except if $f$ is a Reidemeister move of type I$_-^e$ or II$^e$ (with 
$e=\pm$). In these cases $f^0$ is a morphism of degree $-e$.
\vskip 12pt
\noi{\bf 3.3 Elementary moves for oriented link diagrams.} Consider an oriented link diagram $D$. Denote by $w$ the function sending each point $x$ outside of $D$ to the
winding number of $D$ about $x$. The sign $(-1)^{w(x)}$ will be called the $D$-sign of $x$. If $x$ is a crossing of $D$, the function $w$ takes three
values $n-1,n,n+1$ near $x$. The sign $(-1)^n$ will also be called the $D$-sign of $x$.

In order to describe the orientation of a curve we'll use the following convention: suppose $a$ is a sign. Then the figure
$$\begin{tikzpicture} \draw [>=latex,->] (0,0)--(1,0);\draw (.9,.2) node {$a$};\end{tikzpicture}$$
means that the curve is oriented by the arrow if $a=+$ and by the opposite orientation if $a=-$.

Let $f:D\rightarrow D'$ be a Reidemeister move of type I$_e^+$ between oriented diagrams and $g$ be its inverse move. Let $a$ be the winding number of the created
loop and $h$ be the $D'$-sign of the created crossing. We say that $f$ is a Reidemeister move of type I$^+(e,a,h)$ and that $g$ is a Reidemeister move
of type I$^-(e,a,h)$. In this case the diagram $D'$ is oriented as follows:
$$
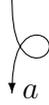
\begin{tikzpicture}[scale=1/8] \draw [domain=0:10] plot ({12-2*cos(36*\x)},1/13*\x^3-15/13*\x^2+63/13*\x);\draw [->,>=latex] (10,0.1)--(10,0); \draw (12,.2) node {$a$};
\end{tikzpicture}$$

Let $f:D\rightarrow D'$ be a Reidemeister move of type II$^+$ between oriented diagrams and $g$ be its inverse move. Let $\Gamma$ be the created loop. This loop is
oriented by the orientation of the plane. Let $a$ (resp. $b$) be the sign which is equal to $+$ if and only if the orientation of the top branch (resp. the bottom branch)
agrees with the orientation of $\Gamma$. Let $h$ be the commun $D'$-sign of the created crossings. We say that $f$ is a Reidemeister move of 
type II$^+(a,b,h)$ and that $g$ is a Reidemeister move of type II$^-(a,b,h)$. In this case the diagram $D'$ is oriented as follows:
$$
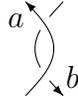
\begin{tikzpicture}[scale=1/8] \draw [domain=0:10] plot ({14+3*sin(18*\x)},\x) ;\draw [domain=0:1.6] plot ({18-3*sin(18*\x)},\x) ;
\draw [domain=3:7] plot ({18-3*sin(18*\x)},\x) ;\draw [domain=8.4:10] plot ({18-3*sin(18*\x)},\x) ;
\draw [->,>=latex] (14.1,9.9)--(14,10);\draw [->,>=latex] (17.9,.1)--(18,0); \draw (13,8) node {$a$};\draw (19,2) node {$b$};
\end{tikzpicture}$$
For technical reasons we'll say also that $f$ is a Reidemeister move of type II$^e(e',a,b,h)$ if $f$ is a Reidemeister move of type II$^e(a,b,h)$ (resp. II$^e(b,a,h)$)
if $e'=+$ (resp. $e'=-$).

Let $f:D\rightarrow D'$ be a Reidemeister move of type III$_e$ between oriented diagrams. Let $\Theta$ be the triangle in $D$ modified by $f$. This triangle is oriented
by the orientation of the plane. Let $a$ (resp. $b$, $c$) be the sign which is equal to $+$ if and only if the orientation of the top branch (resp. the middle branch, 
the bottom branch) agrees with the orientation of the triangle. Let $h$ be the $D$-sign of the center of $\Theta$. We say that $f$ is a 
Reidemeister move of type III$(e,a,b,c,h)$. The diagram $D$ is oriented as follows:
$$
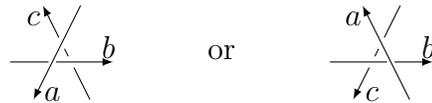
\begin{tikzpicture}[scale=1/4] \draw (-2.736,-.5)--(-.7,-.5);\draw [->,>=latex] (-.3,-.5)--(2.736,-.5);\draw (2.5,0.3) node {$b$};\draw (-1.5,1.8) node {$c$};
\draw (1.5,-2.5)--(.65,-.8);\draw (.35,-.2)--(.15,.2);\draw [->,>=latex] (-.15,.8)--(-1,2.5);\draw [->,>=latex] (1,2.5)--(-1.5,-2.5);\draw (-.5,-2.2) node {$a$};
\draw [->,>=latex] (17.7,-.5)--(19.736,-.5);\draw  (17.3,-.5)--(14.264,-.5);\draw [->,>=latex] (18.5,-2.5)--(16,2.5);\draw (19.5,0.3) node {$b$};
\draw [->,>=latex] (16.35,-.8)--(15.5,-2.5);\draw (16.65,-.2)--(16.85,.2);\draw (17.15,.8)--(18,2.5);\draw (16.5,-2.2) node {$c$};\draw (15.5,1.8) node {$a$};
\draw (8.5,0) node {or};
\end{tikzpicture}$$
depending if $e=+$ or $e=-$.

It is easy to see that the inverse move of $f$ is a Reidemeister move of type III$(e,-a,-b,-c,-h)$.

Let $f:D\rightarrow D'$ be a surgery move of index $0$ between oriented diagrams and $g$ be its inverse move. Let $a$ be the winding number of the created circle and 
$h$ be the $D$-sign of a point in this circle. We say that $f$ is a surgery move of type $(0,a,h)$ and that $g$ is a surgery move of type $(2,a,h)$.

Let $f:D\rightarrow D'$ be a surgery move of index $1$ between oriented diagrams. Let $\gamma$ be a path inducing this surgery. Let $a$ be the sign which is equal to $+$ 
if and only if $\gamma$ is on the left of $D$ and $h$ be the $D$-sign of the middle point in $\gamma$. We say that $f$ is a surgery move of type $(1,a,h)$. It is
easy to see that the inverse move of $f$ is a surgery move of type  $(1,-a,-h)$.
\vskip 12pt
\noi{\bf 3.4 The correspondance $f\mapsto f^1$.}
\vskip 12pt
Consider an elementary move $f:D\rightarrow D'$ where $D$ and $D'$ are oriented link diagrams. This move is called an elementary move of oriented diagrams if $D$ and
$D'$ have the same orientation outside the modification area. 

Here we'll associate to each elementary move $f:D\rightarrow D'$ of oriented diagrams a morphism $f^1$ from $KH(D)$ to $KH(D')$.

Denote by $X_-$ and $X'_-$ the set of negative crossings of $D$ and $D'$.

Suppose the move is a Reidemeister move of type I$_-^+$ or II$^+$. In this case $X'_-$ is the union of $X_-$ and one crossing $x$ and the map $f^1$ is defined by:
$$\forall (u,v)\in kh(D)\times\Lambda^+(X_-), \ \ f^1(u\otimes v)=(-1)^{|u|}f^0(u)\otimes x\vc v$$
where $|u|$ is the degree of $u$.

Suppose the move is a Reidemeister move of type I$_-^-$ or II$^-$. In this case $X_-$ is the union of $X'_-$ and one crossing $x$ and the map $f^1$ is defined by:
$$\forall (u,v)\in kh(D)\times\Lambda^+(X'_-),\ \ f^1(u\otimes x\vc v)=-(-1)^{|u|}f^0(u)\otimes v$$

In all other cases, $X_-$ and $X'_-$ are the same and the map $f^1$ is defined by:
$$\forall (u,v)\in kh(D)\times\Lambda^+(X_-), \ \ f^1(u\otimes v)=f^0(u)\otimes v$$

It is not difficult to see that each map $f^1$ is a morphism of complexes of degree $0$. Moreover is $f$ is a Reidemeister move with inverse move $g$, $g^1$ is a 
homotopy inverse of $f^1$. 

Notice that each map $f^0$ and each map $f^1$ are invariant under any endomorphism of $R$.
\vskip 12pt
\noi{\bf 3.5 The correspondance $f\mapsto f^{\rde K}$.}
\vskip 12pt
Consider elements $A(e,a,h),X(a,h),Y(a,h),Z(a,h)$ in $R^*$, elements $B(a,b,h)$ in $(R\otimes R)^*$ and elements $C(e,a,b,c)$ in $(R\otimes R\otimes R)^*$, depending
on signs $e,a,b,c,h$. Such a data ${\rde K}=(A,B,C,X,Y,Z)$ will be called a Khovanov data. Using this data we'll construct a correspondance associating to each elementary
move $f:D\rightarrow D'$ a morphism $f^{\rde K}:KH(D)\rightarrow KH(D')$.

Let $f:D\rightarrow D'$ be a Reidemeister move of type I$^+(e,a,h)$ and $g$ be its inverse move. Let $p\in D$ be a point in the modified branch of $D$. The morphisms 
$f^{\rde K}$ and $g^{\rde K}$ are defined by:
$$f^{\rde K}=f^1\circ\widehat T_p(A(e,a,h))\hskip 48pt g^{\rde K}=\widehat T_p(A(e,a,h)^{-1})\circ g^1$$
 
Let $f:D\rightarrow D'$ be a Reidemeister move of type II$^+(a,b,h)$ and $g$ be its inverse move. Let $p$ be a point in the top branch and $q$ be a point in the
bottom branch. The morphisms $f^{\rde K}$ and $g^{\rde K}$ are defined by:
$$f^{\rde K}=f^1\circ\widehat T(B(a,b,h))\hskip 48pt g^{\rde K}=\widehat T(B(a,b,h)^{-1})\circ g^1$$
where $\widehat T$ is the map: $u\otimes v\mapsto \widehat T_p(u)\widehat T_q(v)$.

Let $f:D\rightarrow D'$ be a Reidemeister move of type III$(e,a,b,c,h)$ and $g$ be its inverse move. The move $g$ is a Reidemeister move of type III$(e,-a,-b,-c,-h)$.
Let $B_1$ (resp. $B_2$, $B_3$) be the top branch (resp. the middle branch, the bottom branch) near the modified triangle in $D$. For $i=1,2,3$, denote by $p_i$ a
point in $B_i$ and denote by $\widehat T$ the map: $u\otimes v\otimes w\mapsto \widehat T_{p_1}(u)\widehat T_{p_2}(v)\widehat T_{p_3}(w)$. Suppose that $D(e,a,b,c,h)$ are
elements in $(R\otimes R\otimes R)^*$ then we can set: $\widehat f=f^1\circ \widehat T(D(e,a,b,c,h))$. But the condition: $\widehat f\circ\widehat g\sim$Id is equivalent
to the condition:
$$D(e,a,b,c,h)D(e,-a,-b,-c,-h)=1$$
and that's equivalent to the fact that $D(e,a,b,c,h)$ is on the form $D'(e,ah,bh,ch)^h$. So we define $f^{\rde K}$ by:
$$f^{\rde K}=f^1\circ\widehat T(C(e,ah,bh,ch)^h)$$

Let $f:D\rightarrow D'$ be a surgery move of type $(0,a,h)$ and $g$ be its inverse move. Let $p$ be a point in the created circle. Then we define $f^{\rde K}$ and 
$g^{\rde K}$ by:
$$f^{\rde K}=\widehat T_p(X(a,h))\circ f^1\hskip 48pt g^{\rde K}=g^1\circ \widehat T_p(Z(a,h))$$

Let $f:D\rightarrow D'$ be a surgery move of type $(1,a,h)$. Let $p$ be a point in the boundary of the surgery path. Then $f^{\rde K}$ is defined by:
$$f^{\rde K}=f^1\circ\widehat T_p(Y(a,h))$$

It is easy to see that, for every elementary move $f$, $f^{\rde K}$ is a morphism of degree $0$ well defined up to homotopy. Moreover, if $g$ is the inverse 
move of a Reidemeister move $f$, $g^{\rde K}$ is a homotopy inverse of $f^{\rde K}$.

\vskip 24pt
\noi{\bf 4. Movie moves.}
\vskip 24pt
In all this section ${\rde K}=(A,B,C,X,Y,Z)$ is a given Khovanov data. For technical reasons we define the elements $B^e(a,b,h)$ by:
$$B(a,b,h)=u\otimes v\ \ \Longrightarrow\ \ \ B^0(a,b,h)=uv\ \ \ B^+(a,b,h)=u\otimes v\ \ \ B^-(b,a,h)=v\otimes u$$
\vskip 12pt

Let ${\rde D}$ be the set of link diagrams or the set of oriented link diagrams and $A$ be a correspondance associating to each $D\in{\rde D}$ a $K$-complex $A(D)$.

In ${\rde D}$ we have elementary moves and every elementary move $f:D\rightarrow D'$ has an inverse move: $\overline f:D'\rightarrow D$.

Consider diagrams $D_0$, $D_1$, \dots, $D_n$ in ${\rde D}$ and elementary moves $f_i:D_{i-1}\rightarrow D_i$. Such a sequence $\varphi=(D_0,D_1,\dots D_n)$ (or 
$\varphi=(f_1,f_2,\dots f_n)$) will be called a movie sequence from $D_0$ to $D_n$. If $D=D_0$ and $D'=D_n$ are the diagrams of two links $L$ and $L'$, the movie 
sequence $\varphi$ induces a cobordism from $L$ to $L'$ and this cobordism is oriented in the oriented case. By replacing every elementary moves by its inverse, we
get a new movie sequence: $\overline\varphi=(\overline f_n, \dots, \overline f_1)$ from $D'$ to $D$ and the cobordism associated to $\overline\varphi$ is isotopic
to the opposite of the cobordism associated to $\varphi$.

Suppose all the elementary moves of a movie sequence $\varphi$ are Reidemeister moves. Then the movie sequence induces an isotopy from $L$ to $L'$. Suppose also
$D$ is equal to $D'$ and the isotopy from $L$ to $L'=L$ is isotopic to the identity. In this case the movie sequence $\varphi$ will be called a movie move of type I.

Consider $\varphi=(D_0,D_1,\dots D_p)$ and $\psi=(D'_0,D'_1,\dots,D'_q)$ two movie sequences from $D_0=D'_0$ to $D_p=D'_q$. The diagrams $D_0$ and $D_p$ are
the diagrams of two links $L$ and $L'$ and the two movie sequences induce two cobordisms from $L$ to $L'$. Suppose these two cobordisms are isotopic.
Then the pair $(\varphi,\psi)$ will be called a movie move of type II.

Consider a correspondance $f\mapsto \widehat f$ associating to every elementary move $f:D\rightarrow D'$ in ${\rde D}$ a morphism $\widehat f$ from $A(D)$ to $A(D')$
and suppose that this correspondance has the following property:  for every Reidemeister move $f$ with inverse move $g$, $\widehat g$ is a homotopy inverse of 
$\widehat f$.

Let $\varphi=(f_1,f_2,\dots,f_n)$ be a movie sequence from $D$ to $D'$. Then we have a map $\widehat\varphi=\widehat f_n\circ\dots\circ\widehat f_1$ from $A(D)$ to 
$A(D')$.

Suppose this correspondance is coming from a functor from the category of cobordisms of links to the homotopy category of $K$-complexes. In this case, for each movie 
move $\varphi$ of type I, the morphism $\widehat\varphi$ is homotopic to the identity and for each movie move $(\varphi,\psi)$ of type II, the two morphisms 
$\widehat\varphi$ and $\widehat\psi$ are homotopic.

Let $\varphi$ be a movie move of type I. Then the morphism $\widehat{\overline\varphi}$ is a homotopy inverse of $\widehat\varphi$ and the condition ''$\widehat\varphi$
is homotopic to the identity'' implies the same condition for the morphism $\widehat{\overline\varphi}$.

But if $(\varphi,\psi)$ is a movie move of type II, we get two conditions: $\widehat\varphi$ and $\widehat\psi$ are homotopic and $\widehat{\overline\varphi}$ and
$\widehat{\overline\psi}$ are also homotopic.

So a movie move of type I induces one condition on the correspondance and a movie move of type II induces two conditions.

We'll apply this construction in three cases:

--- ${\rde D}$ is the set of link diagrams, $A(D)$ is the Khovanov complex $kh(D)$ and the correspondance is: $f\mapsto f^0$.

--- ${\rde D}$ is the set of oriented link diagrams, $A(D)$ is the Khovanov complex $KH(D)$ and the correspondance is: $f\mapsto f^1$.

--- ${\rde D}$ is the set of oriented link diagrams, $A(D)$ is the Khovanov complex $KH(D)$ and the correspondance is: $f\mapsto f^{\rde K}$.

\vskip 12pt
\noi{\bf 4.0 Notation: } From now on we'll use the following notation: for every element $u\in R$ and every sign $e$ we set:
$$u^{(e)}=\left\{\matrix{u&\hbox{if}\ e=+\cr\overline u&\hbox{if}\ e=-\cr}\right.$$
\vskip 12pt
\noi{\bf 4.1 Movie moves of type MVM$_0$.} Consider a link diagram $D$ and two elementary moves $f$ and $g$ modifying $D$ on disjoint areas. Let $D_1$ be the diagram $D$
modified by $f$ and $D_2$ be the diagram $D$ modified by $g$. The move $f$ (resp. $g$) induces a well defined move $f'$ (resp. $g'$) on $D_2$ (resp. $D_1$) and we have a
commutative diagram:
$$\diagram{D&\hfl{f}{}&D_1\cr \vfl{g}{}&&\vfl{g'}{}\cr D_2&\hfl{f'}{}&D'\cr}$$
So we get two movie sequences: $\varphi=(f,g')$ and $\psi=(g,f')$ and a movie move $(\varphi,\psi)$ of type II. This movie move will be called a movie move of type 
MVM$_0$.
\vskip 12pt
\noi{\bf 4.1.a Lemma:} {\sl Let $(\varphi,\psi)$ be a movie move of type MVM$_0$ associated to elementary moves $f$ and $g$. Then we have:
$$\varphi^0=(-1)^{pq}\psi^0$$
where $p$ and $q$ are the degrees of the map $f^0$ and $g^0$.

In the oriented case we have:
$$\varphi^1=\psi^1\hskip 48pt \varphi^{\rde K}=\psi^{\rde K}\hskip 48pt \overline\varphi^{\rde K}=\overline\psi^{\rde K}$$}
\vskip 12pt
\noi{\bf Proof:} It is easy to see that we have some commutativity in the graded sense and that implies the first formula and therefore the other ones.

For example, consider the case where $f$ and $g$ are Reidemeister moves of type I$^+(-,a,h)$ and I$^+(-,b,k)$. Let $x$ (resp. $y$)  be the crossing created by $f$ (resp.
$g$) and $p$ (resp. $q$) be a regular point in $D$ near $x$ (resp. near $y$). Let $X_-$ be the set of negative crossings of the diagram $D$. So, for every $u\in kh(D)$
and every $v\in\Lambda^+(X_-)$, we have:
$$\varphi^0(u)=f^0(y\otimes u\otimes 1)=x\vc y\otimes u\otimes 1\otimes 1$$
$$\psi^0(u)=g^0(x\otimes u\otimes 1)=y\vc x\otimes u\otimes 1\otimes 1\ \ \Longrightarrow\ \ \psi^0=-\varphi^0$$
$$\varphi^1(u\otimes v)=f^1(g'^0(u)\otimes y\vc v)=\varphi^0(u)\otimes x\vc y\vc v$$
$$\psi^1(u\otimes v)=g^1(f'^0(u)\otimes x\vc v)=\psi^0(u)\otimes y\vc x\vc v\ \ \Longrightarrow\ \ \psi^1=\varphi^1$$
$$\varphi^{\rde K}=\varphi^1\circ \widehat T_p(A(-,a,h))\widehat T_q(A(-,b,k))=\psi^1\circ \widehat T_q(A(-,b,k))\widehat T_p(A(-,a,h))=\psi^{\rde K}$$
The last formula $\overline\varphi^{\rde K}=\overline\psi^{\rde K}$ is just the formula $\varphi^{\rde K}=\psi^{\rde K}$ applied for Reidemeister moves $\overline g$ and 
$\overline f$.\cqfd

\vskip 12pt
\noi{\bf 4.2 Movie moves of type MVM$_1$.} Consider a link diagram $D$ and a regular point $p$ in $D$. Let $e$ be a sign. We may apply a Reidemeister 
move of type I$_e^+$ on $D$ near $p$, then a new Reidemeister move of type I$_{-e}^+$ in the created loop. By applying a Reidemeister move of type II$^-$ we may remove
the two created crossings and recover the diagram $D$. So we have a movie move of type I. This move will be called a movie move of type MVM$_1(e)$ near $p$.
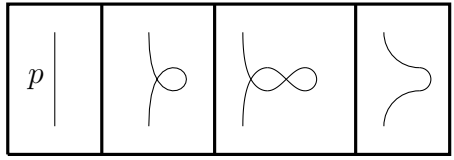
$$\begin{tikzpicture}[scale=1/8] \draw (0,0) -- (0,10) ; \draw (0,5) node[left] {$p$} ;
\draw [domain=0:10] plot ({12-2*cos(36*\x)},1/13*\x^3-15/13*\x^2+63/13*\x) ;
\draw [domain=0:3] plot ({22-2*cos(36*\x)},1/13*\x^3-15/13*\x^2+63/13*\x) ;
\draw [domain=7:10] plot ({22-2*cos(36*\x)},1/13*\x^3-15/13*\x^2+63/13*\x) ;
\draw [domain=0:4] plot ({22+2*cos(72)+\x},{5+16/13*cos(45*\x)}) ; \draw [domain=0:4] plot ({22+2*cos(72)+\x},{5-16/13*cos(45*\x)}) ;
\draw ({22+2*cos(72)+4},5-16/13) arc (-90:90:16/13) ;
\draw (35,0) arc (180:90:49/13) ; \draw (35,10) arc (-180:-90:49/13) ; \draw (35+49/13,49/13) arc (-90:90:16/13) ; 
\draw [very thick] (-5,-3) -- (42,-3) -- (42,13) -- (-5,13) -- (-5,-3) ; \draw [very thick] (5,-3) -- (5,13) ; \draw [very thick] (17,-3) -- (17,13) ; 
\draw [very thick] (32,-3) -- (32,13) ; \end{tikzpicture}$$
\vskip 12pt
\noi{\bf 4.2.a Lemma:} {\sl Let $\varphi$ be a movie move of type MVM$_1(e)$ near a point $p\in D$. Then we have:
$$\varphi^0=-eT_p({\omega\over\omega^{(e)}})$$

In the oriented case, suppose the first move is a Reidemeister move of type I$^+(e,a,h)$. Then we have:
$$\varphi^1=-eT_p({\omega\over\omega^{(e)}})=-e\widehat T_p({\omega^{(ah)}\over\omega^{(eah)}})$$
$$\varphi^{\rde K}=-e\widehat T_p({\omega^{(ah)}\over\omega^{(eah)}}{A(e,a,h)A(-e,-a,h)\over B^0(a,a,h)})$$}
\vskip 12pt
\noi{\bf Proof: } A straightforward computation gives the following:
$$e=+\ \ \Longrightarrow\ \ \varphi^0=-\hbox{Id}$$
$$e=-\ \ \Longrightarrow\ \ \varphi^0=T_p({\omega\over\overline\omega})$$
So we get the first formula. The second formula: $\varphi^1=-e T_p(\omega/\omega^{(e)})$ is easy to deduce.

Since the first move is a Reidemeister move of type I$^+(e,a,h)$, the orientation of $D$ is given by the following:
$$\begin{tikzpicture}[scale=1/8] \draw [->,>=latex] (0,10) -- (0,0) ; \draw (0,5) node[left] {$p$} ; \draw (1.7,.4) node {$a$};  \end{tikzpicture}$$
So the second move is a Reidemeister move of type I$^+(-e,-a,h)$ and the last one is a Reidemeister move of type II$^-(a,a,h)$. The result follows.\cqfd
\vskip 12pt
\noi{\bf 4.3 Movie moves of type MVM$_2$.} Let $D$ be a link diagram and $x$ be a crossing of $D$. We may modify $D$ by a Reidemeister move of type II$^+$ which creates
two new crossings $y$ and $z$ and then remove $x$ and $y$ by a Reidemeister move of type II$^-$. 
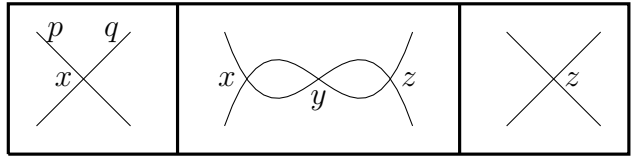
$$\begin{tikzpicture}[scale=1/8]  \draw (0,0) -- (10,10) ; \draw (0,10) -- (10,0) ; \draw (5,5) node[left] {$x$} ;
\draw (2,10) node {$p$} ; \draw (8,10) node {$q$} ;
\draw [domain=0:4] plot (20+5*\x,3/2*\x^3-9*\x^2+29/2*\x) ; \draw [domain=0:4] plot (20+5*\x,-3/2*\x^3+9*\x^2-29/2*\x+10) ;
\draw (22.362374,5) node[left] {$x$} ; \draw (30,5) node[below] {$y$} ; \draw (37.637626,5) node[right] {$z$} ;
\draw (50,0) -- (60,10) ; \draw (50,10) -- (60,0) ; \draw (55,5) node[right] {$z$} ;
\draw [very thick] (-3,-3) -- (63,-3) -- (63,13) -- (-3,13) -- (-3,-3) ; \draw [very thick] (15,-3) -- (15,13) ; \draw [very thick] (45,-3) -- (45,13) ;
\end{tikzpicture}$$
So we get a movie move. We say that this move is a movie move of type MVM$_2(e)$ near $(p,q)$ where $e=+$ (resp. $e=-$) if the branch containing $p$ is over (resp. under)
the branch containing $q$.
\vskip 12pt
\noi{\bf 4.3.a Lemma:} {\sl Let $e$ be a sign and $\varphi$ be a movie move of type MVM$_2(e)$ near $(p,q)$. Then we have:
$$\varphi^0=eT_p(\omega)T_q(\omega^{-1})$$

In the oriented case, suppose the first move is a Reidemeister move of type II$^+(e,a,b,h)$. Then we have:
$$\varphi^1=-abT_p(\omega)T_q(\omega^{-1})=-ab\widehat T(\omega^{(-bh)}\otimes{1\over\omega^{(ah)}})$$
$$\varphi^{\rde K}=-ab\widehat T\Bigl((\omega^{(-bh)}\otimes{1\over\omega^{(ah)}}){B^e(a,b,h)\over B^e(-a,-b,h)}\Bigr)$$
where $\widehat T$ is the map $u\otimes v\mapsto\widehat T_p(u)\widehat T_q(v)$.}
\vskip 12pt
\noi{\bf Proof: } A straightforward computation gives the first formula. Let $f$ and $g$ be the Reidemeister moves in $\varphi$ and $X_-$ be the set of negative 
crossings of $D$.

So it is easy to see that the orientation of the second diagram is given by:
$$
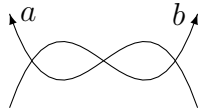
\begin{tikzpicture}[scale=1/8] \draw [domain=0:4] plot (20+5*\x,3/2*\x^3-9*\x^2+29/2*\x) ; \draw [domain=0:4] plot (20+5*\x,-3/2*\x^3+9*\x^2-29/2*\x+10) ;
\draw [->,>=latex] (39.9,9.71)--(40,10); \draw (38,10) node {$b$};\draw [->,>=latex] (20.1,9.71)--(20,10);\draw (22,10) node {$a$};\end{tikzpicture}$$
and the two moves are Reidemeister moves of type II$^+(e,a,b,h)$ and II$^-(e,-a,-b,h)$. Moreover the signs of the crossings $x, y, z$ are: $-eab,eab,-eab$.

If $eab=+$, $x$ is in $X_-$ and we have, for every $u\in kh(D)$ and every $x\vc v\in\Lambda^+(X_-)$:
$$\varphi^1(u\otimes x\vc v)=(-1)^{|u|}g^1(f^0(u)\otimes z\vc x\vc v)=-(-1)^{|u|}g^1(f^0(u)\otimes x\vc z\vc v)$$
$$=(-1)^{|f^0(u)|}(-1)^{|u|}g^0(u)\otimes z\vc v=-\varphi^0(u)\otimes z\vc v$$
and that implies:
$$\varphi^1=-abT_p(\omega)T_q(\omega^{-1})$$

If $eab=-$, $x$ is not in $X_-$ and we have, for every $u\in kh(D)$ and every $v\in\Lambda^+(X_-)$:
$$\varphi^1(u\otimes v)=(-1)^{|u|}g^1(f^0(u)\otimes y\vc v)=-(-1)^{|f^0(u)|}(-1)^{|u|}g^0(u)\otimes v=\varphi^0(u)\otimes v$$
and that implies again:
$$\varphi^1=-abT_p(\omega)T_q(\omega^{-1})$$

On the other hand, the $D$-sign of a point in the plane between $p$ and $q$ is $abh$. So we get:
$$\varphi^1=-ab\widehat T_p(\omega^{(-bh)})\widehat T_q({1\over\omega^{(ah)}})=-ab\widehat T(\omega^{(-bh)}\otimes{1\over\omega^{(ah)}})$$
and the last formula follows easily.\cqfd
\vskip 12pt
\noi{\bf 4.4 Movie moves of type MVM$_3$.} Let $D$ be a link diagram and $p$ and $q$ be two regular points in $D$. Suppose that these two points may
be joined by a path outside of $D$. Then we have a movie move described by the following figure:
$$
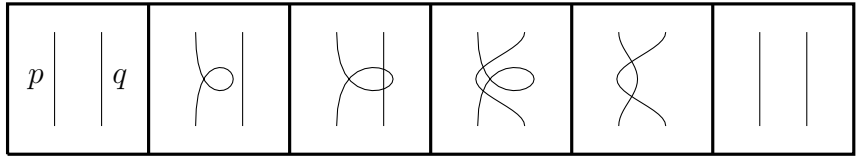
\begin{tikzpicture}[scale=1/8] \draw (0,5) node[left] {$p$} ; \draw (0,0) --(0,10) ; \draw (5,0) -- (5,10) ; \draw (5,5) node[right] {$q$} ;
\draw [domain=0:10] plot ({17-2*cos(36*\x)},1/13*\x^3-15/13*\x^2+63/13*\x) ; \draw (20,0) -- (20,10) ;
\draw [domain=0:10] plot ({33-3*cos(36*\x)},1/13*\x^3-15/13*\x^2+63/13*\x) ; \draw (35,0) -- (35,10) ;
\draw [domain=0:10] plot ({48-3*cos(36*\x)},1/13*\x^3-15/13*\x^2+63/13*\x) ; \draw [domain=0:10] plot ({50-2.6*(1-cos(36*\x)},\x) ;
\draw [domain=0:10] plot ({60+1-cos(36*\x)},\x) ;\draw [domain=0:10] plot ({65-2.6*(1-cos(36*\x)},\x) ;
\draw (75,0) -- (75,10) ; \draw (80,0) -- (80,10) ;
\draw [very thick] (-5,-3) -- (85,-3) -- (85,13) -- (-5,13) -- (-5,-3) ; \draw [very thick] ((10,-3) -- (10,13) ; \draw [very thick] ((25,-3) -- (25,13) ; 
\draw [very thick] (40,-3) -- (40,13) ; \draw [very thick] (55,-3) -- (55,13) ; \draw [very thick] (70,-3) -- (70,13) ;\end{tikzpicture}$$
This movie move depends on two sign $e$ and $e'$. The sign $e$ is the sign of the first Reidemeister move of type I and $e'$ is $+$ if and only if the loop goes
over the branch containing $q$. We say that this move is a movie move of type MVM$_3(e,e')$ near $(p,q)$.
\vskip 12pt
\noi{\bf 4.4.a Lemma: } {\sl Let $e$ and $e'$ be two signs. Let $p$ and $q$ be two points in a link diagram $D$ and $\varphi$ be a movie move of type MVM$_3(e,e')$ near
$(p,q)$. Then we have:
$$\varphi^0=T_p({\omega\over\omega^{(e)}})$$

In the oriented case, let $a$, $b$, $h$ be the signs such that the first move of $\varphi$ is a Reidemeister move of type I$^+(e,a,h)$ and $a$ is equal to $b$ if and
only if the two branches are oriented in the same way. Then we have:
$$\varphi^1=e T_p({\omega\over\omega^{(e)}})=e\widehat T_p({\omega^{(ah)}\over\omega^{(eah)}})$$
$$\varphi^{\rde K}=e\widehat T\Bigl(({\omega^{(ah)}\over\omega^{(eah)}}\otimes 1)({A(e,a,h)\over A(e,a,-h)}\otimes 1){B^{e'}(a,b,-abh)\over B^{e'}(-a,b,abh)}C\Bigr)$$
where $\widehat T$ is the map $u\otimes v\mapsto\widehat T_p(u)\widehat T_q(v)$ and $C$ is defined by:
$$e'=+,\ \ C(e,-ah,-ah,bh)=u\otimes v\otimes w\ \ \Longrightarrow\ \ C=(uv)^{-h}\otimes w^{-h}$$
$$e'=-,\ \ C(e,bh,-ah,-ah)=u\otimes v\otimes w,\ \ \Longrightarrow\ \ C=(wv)^{-h}\otimes u^{-h}$$.}
\vskip 12pt
\noi{\bf Proof: } A straightforward computation gives the first formula.

The diagram $D$ is oriented as follows:
$$
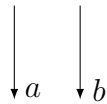
\begin{tikzpicture}[scale=1/8] \draw [->,>=latex] (0,10)--(0,0) ; \draw [->,>=latex] (7,10)--(7,0) ;\draw (2,1) node {$a$};\draw (9,1) node {$b$};\end{tikzpicture}$$
Denote by $(f_1,f_2,f_3,f_4,f_5)$ the sequence $\varphi$. Then the types of $f_1$, $f_2$, $f_3$, $f_4$, $f_5$  are:
$$\hbox{I}^+(e,a,h),\ \hbox{II}^+(a,b,-abh),\ \hbox{III}(e,a,a,-b,-h),\ \hbox{I}^-(e,a,-h),\ \hbox{II}^-(-a,b,abh)$$
if $e'=+$ and:
$$\hbox{I}^+(e,a,h),\ \hbox{II}^+(b,a,-abh),\ \hbox{III}(e,-b,a,a,-h),\ \hbox{I}^-(e,a,-h),\ \hbox{II}^-(b,-a,abh)$$
if $e'=-$.

In the case: $e=+$, only $f_2$ and $f_5$ modify the set of negative crossings. Then, during the move $\varphi$, only one negative crossing appears and then disappears.
So we get:
$$\varphi^1=T_p({\omega\over\omega^{(e)}})$$
In the case: $e=-$, $f_1$ creates a negative crossing $x$, $f_2$ creates a negative crossing $y$, $f_4$ removes $x$ and $f_5$ removes $y$. This operation produces a sign
and we have:
$$\varphi^1=-T_p({\omega\over\omega^{(e)}})$$
So we get the second formula. The other formulae are easy to check.\cqfd

\vskip 12pt
\noi{\bf 4.5 Movie moves of type MVM$_4$.} Consider a link diagrams and two consecutive crossings $x$ and $y$. Then we have a movie move described by the following 
figure:
$$
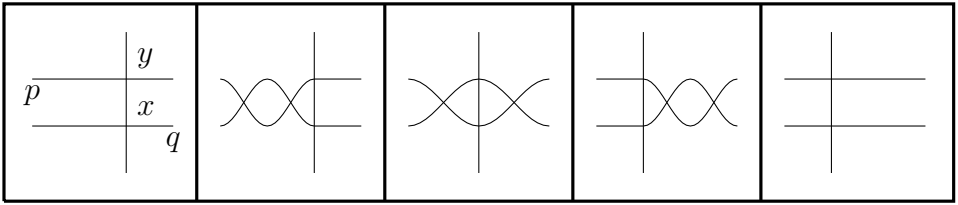
\begin{tikzpicture}[scale=1/8] \draw (0,5) -- (15,5) ; \draw (0,10) -- (15,10) ; \draw (10,0) -- (10,15) ; \draw (10,5) node[above right] {$x$} ;
\draw (10,10) node[above right] {$y$} ;
\draw (0,10.5) node[below] {$p$} ; \draw (15,5.5) node[below] {$q$} ;
\draw [domain=0:10] plot (20+\x,{15/2-5/2*cos(36*\x)}) ; \draw [domain=0:10] plot (20+\x,{15/2+5/2*cos(36*\x)}) ; \draw (30,0) -- (30,15) ;
\draw (30,5) -- (35,5) ; \draw (30,10) -- (35,10) ;
\draw [domain=0:15] plot(40+\x,{15/2-5/2*cos(24*\x)}) ; \draw [domain=0:15] plot(40+\x,{15/2+5/2*cos(24*\x)}) ; \draw (95/2,0) -- (95/2,15) ;
\draw [domain=0:10] plot (65+\x,{15/2-5/2*cos(36*\x)}) ; \draw [domain=0:10] plot (65+\x,{15/2+5/2*cos(36*\x)}) ; \draw (65,0) -- (65,15) ;
\draw (60,5) -- (65,5) ; \draw (60,10) -- (65,10) ; \draw (80,5) -- (95,5) ; \draw (80,10) -- (95,10) ; \draw (85,0) -- (85,15) ; 
\draw [very thick] (-3,-3) -- (98,-3) -- (98,18) -- (-3,18) -- (-3,-3) ; \draw [very thick] (17.5,-3) -- (17.5,18) ; \draw [very thick] (37.5,-3) -- (37.5,18) ;
\draw [very thick] (57.5,-3) -- (57.5,18) ; \draw [very thick] (77.5,-3) -- (77.5,18) ; \end{tikzpicture}$$
Up to replacing this movie move by its inverse we may as well suppose that the horizontal branch containing $y$ goes over the other horizontal branch after the first 
Reidemeister move. So this move depends only on an element $e\in\{0,+,-\}$ where $e=+$ (resp. $e=0$, $e=-$) if the vertical branch is over (resp. between, under) the 
two other branches. We say that this move is a movie move of type MVM$_4(e)$ near $(p,q)$.
\vskip 12pt
\noi{\bf 4.5.a Lemma: } {\sl Let $p$ and $q$ be two points in a link diagram $D$. Let $\varphi$ be a movie move of type MVM$_4(e)$ near $(p,q)$. Then the morphism 
$\varphi^0$ is homotopic to the operator:
$$(1-2e^2)T_p(\omega^{-1})T_q(\omega)$$

In the oriented case, suppose that the first move in $\varphi$ is a Reidemeister move of type II$^+(a,b,h)$ and that $c=+$ (resp. $c=-$) if the 
vertical branch is oriented downward (resp. upward). Then we have:
$$\varphi^1\sim (1-2e^2)T_p(\omega^{-1})T_q(\omega)=(1-2e^2)\widehat T_p({1\over\omega^{(bh)}})\widehat T_q(\omega^{(-ah)})$$
$$\varphi^{\rde K}\sim (1-2e^2)\widehat T\Bigl(({1\over\omega^{(bh)}}\otimes\omega^{(-ah)}\otimes 1)({B(a,b,h)\over B(a,b,-h)}\otimes 1)C^{-abh}\Bigr)$$
where  $\widehat T$ is the map $u\otimes v\otimes w\mapsto\widehat T_p(u)\widehat T_q(v)\widehat T_r(w)$ (for any $r$ in the vertical branch)  and $C$ is defined by:
$$e=-,\ \ C(+,bh,ah,abch)C(-,-bh,-ah,abch)=u\otimes v\otimes w\ \ \Longrightarrow\ \ \ C=u\otimes v\otimes w$$
$$e=0,\ \ C(-,bh,abch,ah)C(+,-bh,abch,-ah)=u\otimes v\otimes w\ \ \Longrightarrow\ \ \ C=u\otimes w\otimes v$$
$$e=+,\ \ C(+,abch,bh,ah)C(-,abch,-bh,-ah)=u\otimes v\otimes w\ \ \Longrightarrow\ \ \ C=w\otimes u\otimes v$$} 
\vskip 12pt
\noi{\bf Proof: } The union of the branches is a graph with six points in its boundary. Take a counterclockwise numbering of these points begining with the point
$p$. So we have: $p_1=p$, $p_4=q$. The operator $T_{p_i}$ will be denoted by $T_i$.

There are four complexes $U$, $V$, $W$ and $H$ such that:
$$kh(D)=1\otimes U\oplus x\otimes V\oplus y\otimes W\oplus x\vc y\otimes H$$

Using small paths near $x$ or $y$, we get surgery maps $\gamma_x$ and $\gamma_y$. These maps are morphisms of complexes. The morphism $\gamma_x$ is a map from $U$ to 
$V$, from $V$ to $U$, from $W$ to $H$ and from $H$ to $W$. The morphism $\gamma_y$ is a map from $U$ to $W$, from $W$ to $U$, from $V$ to $H$ and from $H$ to $V$. A 
straightforward computation shows the following:

If $\varphi$ is a movie move of type MVM$_4(-)$, we have (for $u\in U$, $v\in V$, $w\in W$, $h\in H$):
$$\varphi^0(1\otimes u)=-1\otimes u$$
$$\varphi^0(x\otimes v)=-x\otimes T_1(\overline\omega)T_3(\omega^{-1})(v)+\e(\omega^2)y\otimes T_1(\omega^{-2})\gamma_x\gamma_y(v)$$
$$\varphi^0(y\otimes w)=y\otimes\Bigl(-1+T_1(\omega^{-1})T_3(\overline\omega)-T_1(\omega^{-1})T_3(\omega)\Bigr)(w)-\e(\omega^2)x\otimes T_1(\omega^{-2})\gamma_x\gamma_y(w)$$
$$\varphi^0(x\vc y\otimes h)=-x\vc y\otimes h$$

If $\varphi$ is a movie move of type MVM$_4(0)$, we have:
$$\varphi^0(1\otimes u)=1\otimes T_3(\omega)T_5(\overline\omega^{-1})(u)$$
$$\varphi^0(x\otimes v)=x\otimes T_1(\omega)T_5(\overline\omega^{-1})(v)-\e(1)y\otimes \gamma_x\gamma_y(v)$$
$$\varphi^0(y\otimes w)=y\otimes w$$
$$\varphi^0(x\vc y\otimes h)=x\vc y\otimes T_3(\omega^{-1})T_5(\omega)(h)$$

If $\varphi$ is a movie move of type MVM$_4(+)$, we have:
$$\varphi^0(1\otimes u)=-1\otimes T_1(\overline\omega)T_3(\omega^{-1})(u)$$
$$\varphi^0(x\otimes v)=-x\otimes T_1(\omega^{-1})T_5(\omega)(v)$$
$$\varphi^0(y\otimes w)=-y\otimes T_1(\overline\omega)T_3(\omega^{-1})(w)+\e(\omega^{2})x\otimes T_1(\omega^{-2})\gamma_x\gamma_y(w)$$
$$\varphi^0(x\vc y\otimes h)=-x\vc y\otimes h$$

Consider the maps $k_x$ and $k_y$ from $kh(D)$ to $kh(D)$ that vanish on $1\otimes U\oplus x\otimes V\oplus y\oplus W$ and satisfy the following:
$$k_x(x\vc y\otimes h)=x\otimes \gamma_y(h)\hskip 24pt k_y(x\vc y\otimes h)=y\otimes \gamma_x(h)$$
for every $h\in H$. These maps induce the following homotopies:
$$h_x=d(k_x)=d\circ k_x+k_x\circ d\hskip 24pt h_y=d(h_y)=d\circ k_y+k_y\circ d$$
and we have:
$$h_x(1\otimes u)=0\hskip 48pt h_y(1\otimes u)=0$$
$$h_x(x\otimes v)=-x\otimes \gamma_y^2(v)\hskip 48pt h_y(x\otimes v)=-y\otimes \gamma_x\gamma_y(v)$$
$$h_x(y\otimes w)=x\otimes \gamma_x\gamma_y(w)\hskip 48pt h_y(y\otimes w)=y\otimes \gamma_x^2(w)$$
$$h_x(x\vc y\otimes h)=-x\vc y\otimes \gamma_y^2(h)\hskip 48pt h_y(x\vc y\otimes h)=x\vc y\otimes \gamma_x^2(h)$$

Moreover one checks that $\gamma_x^2$ (resp. $\gamma_y^2$) is on $W$ and $H$ (resp. on $V$ and $H$) the map $T_2(\omega)T_4(\alpha)-T_2(\omega\overline\alpha)$ (resp.
the map $T_1(\omega)T_5(\alpha)-T_1(\omega\overline\alpha)$). So we get:
$$\varphi^0=\e(1)T_1(\overline\omega)T_3(\omega^{-1})h_x-\e(\omega^2)T_1(\omega^{-2})h_y-T_2(\omega^{-1})T_5(\omega)$$
in the first case,
$$\varphi^0=\e(1)h_y+T_4(\omega)T_5(\overline\omega^{-1})$$
in the second case and:
$$\varphi^0=\e(\omega^2)T_1(\omega^{-2})h_x-T_1(\overline\omega)T_2(\omega^{-1})$$
in the last case. If we denote by $\sim$ the homotopy relation, we get:
$$\varphi^0\sim -T_2(\omega^{-1})T_5(\omega)\sim -T_4(\overline\omega^{-1})T_1(\overline\omega)=-T_1(\omega^{-1})T_4(\omega)$$
in the first case,
$$\varphi^0\sim T_4(\omega)T_5(\overline\omega^{-1})\sim T_1(\omega^{-1})T_4(\omega)$$
in the second case and:
$$\varphi^0\sim-T_1(\overline\omega)T_2(\omega^{-1})\sim-T_1(\overline\omega)T_4(\overline\omega^{-1})=-T_1(\omega^{-1})T_4(\omega)$$
in the last case. So we get the first formula.

The movie move creates one negative crossing and then remove it. Therefore we get the same formula for $\varphi^1$.

The second diagram $D'$ is oriented as follows:
$$
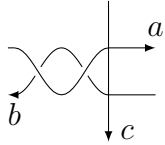
\begin{tikzpicture}[scale=1/8] \draw [domain=0:10] plot (20+\x,{15/2+5/2*cos(36*\x)});\draw (20,10)--(19.3,10);
\draw [domain=0:2.2] plot (20+\x,{15/2-5/2*cos(36*\x)});\draw [->,>=latex] (20,5)--(19.3,5); \draw (20,3) node {$b$};\draw (35,12) node {$a$};\draw (32,1) node {$c$};
\draw [domain=2.8:7.2] plot (20+\x,{15/2-5/2*cos(36*\x)});\draw [domain=7.8:10] plot (20+\x,{15/2-5/2*cos(36*\x)});
\draw [->,>=latex] (30,15)--(30,0);\draw (30,5)--(35,5);\draw [->,>=latex] (30,10)--(35,10);\end{tikzpicture}$$

Denote by $(f_1,f_2,f_3,f_4)$ the sequence $\varphi$. Then the types of $f_1$, $f_2$, $f_3$, $f_4$ are:
$$\hbox{II}^+(a,b,h),\ \hbox{III}(+,-a,-b,-c,-abh),\ \hbox{III}(-,a,b,-c,-abh),\ \hbox{II}^-(a,b,h)$$
if $e=-$,
$$\hbox{II}^+(a,b,h),\ \hbox{III}(-,-a,-c,-b,-abh),\ \hbox{III}(+,a,-c,b,-abh),\ \hbox{II}^-(a,b,h)$$
if $e=0$ and
$$\hbox{II}^+(a,b,h),\ \hbox{III}(+,-c,-a,-b,-abh),\ \hbox{III}(-,-c,a,b,-abh),\ \hbox{II}^-(a,b,h)$$
if $e=+$. The result follows.\cqfd

\vskip 12pt
\noi{\bf 4.6 Movie moves of type MVM$_5$.} Consider a link in $\R^3$ that projects on the plane by an immersion and suppose that every crossing is simple except one
and the image of the link looks like the following near the singular crossing:
$$
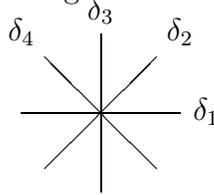
\begin{picture}(60,60) \put(30,0){\line(0,1){60}}\put(0,30){\line(1,0){60}}\put(8.79,8.79){\line(1,1){42.42}}\put(51.21,8.79){\line(-1,1){42.42}}
\zput(70,30){$\delta_1$}\zput(60,60){$\delta_2$}\zput(30,68){$\delta_3$}\zput(0,60){$\delta_4$}\end{picture}$$

If we move the link a little bit, we get a diagram like:
$$
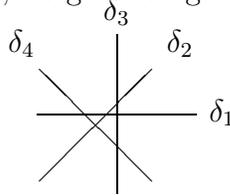
\begin{picture}(60,60) \put(30,0){\line(0,1){60}}\put(0,30){\line(1,0){60}}\put(0.79,4.79){\line(1,1){42.42}}\put(43.21,4.79){\line(-1,1){42.42}}
\zput(70,30){$\delta_1$}\zput(54,56){$\delta_2$}\zput(30,68){$\delta_3$}\zput(-6,56){$\delta_4$}\end{picture}$$
For $\{i,j\}\subset \{1,2,3,4\}$, denote by $x_{ij}$ the intersection of $\delta_i$ and $\delta_j$. One checks the following fact: if, during an small isotopy,
$x_{24}-x_{13}$ makes one counterclockwise turn around $0$, $x_{14}-x_{23}$ and $x_{34}-x_{12}$ make one clockwise turn around $0$. Therefore any small isotopy such that
$x_{ij}-x_{kl}$ makes one turn around $0$ is equivalent to an isotopy where $\delta_3$ and $\delta_4$ are fixed and $\delta_1$ and $\delta_2$ are moving by translations.
During this isotopy, we get diagrams $D_z$ where $z=x_{12}$ and $z$ makes one counterclockwise turn around $x=x_{34}$. If $z$ is not in one of the initial 
$\delta_i$'s, $D_z$ is a link diagram. So we get a movie move which is the composite of eight Reidemeister moves of type III. 
$$
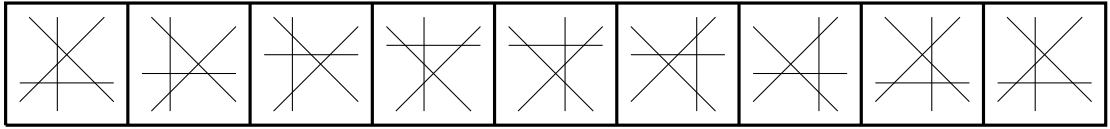
\begin{tikzpicture}[scale=1/8] \draw [very thick] (-1.5,-1.5)--(115.5,-1.5)--(115.5,11.5)--(-1.5,11.5)--(-1.5,-1.5);\draw [very thick] (11.5,-1.5)--(11.5,11.5);
\draw [very thick] (24.5,-1.5)--(24.5,11.5);\draw [very thick] (37.5,-1.5)--(37.5,11.5);\draw [very thick] (50.5,-1.5)--(50.5,11.5);
\draw [very thick] (63.5,-1.5)--(63.5,11.5);\draw [very thick] (76.5,-1.5)--(76.5,11.5);\draw [very thick] (89.5,-1.5)--(89.5,11.5);
\draw [very thick] (102.5,-1.5)--(102.5,11.5);
\draw (0,3)--(10,3);\draw (4,0)--(4,10);\draw (0,1)--(9,10);\draw (10,1)--(1,10);
\draw (13+0,4)--(13+10,4);\draw (13+3,0)--(13+3,10);\draw (13+1,0)--(13+10,9);\draw (13+10,1)--(13+1,10);
\draw (26+0,6)--(26+10,6);\draw (26+3,0)--(26+3,10);\draw (26+1,0)--(26+10,9);\draw (26+10,1)--(26+1,10);
\draw (39+0,7)--(39+10,7);\draw (39+4,0)--(39+4,10);\draw (39+1,0)--(39+10,9);\draw (39+9,0)--(39+0,9);
\draw (52+0,7)--(52+10,7);\draw (52+6,0)--(52+6,10);\draw (52+1,0)--(52+10,9);\draw (52+9,0)--(52+0,9);
\draw (65+0,6)--(65+10,6);\draw (65+7,0)--(65+7,10);\draw (65+0,1)--(65+9,10);\draw (65+9,0)--(65+0,9);
\draw (78+0,4)--(78+10,4);\draw (78+7,0)--(78+7,10);\draw (78+0,1)--(78+9,10);\draw (78+9,0)--(78+0,9);
\draw (91+0,3)--(91+10,3);\draw (91+6,0)--(91+6,10);\draw (91+0,1)--(91+9,10);\draw (91+10,1)--(91+1,10);
\draw (104+0,3)--(104+10,3);\draw (104+4,0)--(104+4,10);\draw (104+0,1)--(104+9,10);\draw (104+10,1)--(104+1,10);
\end{tikzpicture}$$
This move depends essentially on the heights $c_1,c_2,c_3,c_4$ of the lines $\delta_1,\delta_2,\delta_3,\delta_4$. We say that this move is a movie move of type 
MVM$_5(c_1,c_2,c_3,c_4)$.
\vskip 12pt
\noi{\bf 4.6.a Lemma: } {\sl Let $\varphi$ be a movie move of type MVM$_5$. Then $\varphi^0$ is homotopic to the identity.

In the oriented case, $\varphi^1$ is also homotopic to the identity.}
\vskip 12pt
\noi{\bf Proof:} Consider a movie move $\varphi$ of type MVM$_5(c_1,c_2,c_3,c_4)$. The morphism $\varphi^0$ acts on $kh(D)$ for some link diagram $D$. The modification
area of $\varphi$ involves six crossings and, to describe $\varphi^0$, we have to decompose $kh(D)$ into a direct sum of $64$ complexes. So an explicit description of
$\varphi^0$ is almost impossible, at least by hand. But it is possible to compute $\varphi^0$ using a program on a computer and we check that $\varphi^0$ is an 
idempotent (i.e. $\varphi^0\circ\varphi^0=\varphi^0$) when the sequence $(c_i)$ is increasing or decreasing. This property is may be related to the following
fact which is easy to check: all the Reidemeister moves of $\varphi$ are of type III$_+$ if $(c_i)$ is increasing and of type III$_-$ if $(c_i)$ is decreasing.

Since $\varphi^0$ is a homotopy equivalence it has a homotopy inverse $\psi$. So when $(c_1,c_2,c_3,c_4)$ is increasing or decreasing, we have:
$$\varphi^0\sim \psi\circ\varphi^0\circ\varphi^0=\psi\circ\varphi^0\sim\hbox{Id}$$
and $\varphi^0$ is homotopic to the identity.

In the other cases, we have to proceed differently. Let $S$ be the set of sequences $(c_1,c_2,c_3,c_4)$ of distinct reals such that the lemma is true for every movie
move of type MVM$_5(c_1,c_2,c_3,c_4)$. So we have:
$$c_1>c_2>c_3>c_4\ \ \Longrightarrow\ \ \ (c_1,c_2,c_3,c_4)\in S$$
On the other hand, it is easy to see that a movie move of type MVM$_5(c_1,c_2,c_3,c_4)$ is conjugate to a movie move of type MVM$_5(c_2,c_3,c_4,c_1)$. So we have:
$$(c_1,c_2,c_3,c_4)\in S\ \ \Longleftrightarrow\ \ (c_2,c_3,c_4,c_1)\in S$$

Let $c_1,c_2,c_3,c_4$ be distinct reals and $\varphi$ be a movie move of type MVM$_5(c_1,c_2,c_3,c_4)$ acting on a diagram $D$. It is possible to modify $D$ by a
Reidemeister move of type II in such a way to get the following configuration in $D$:
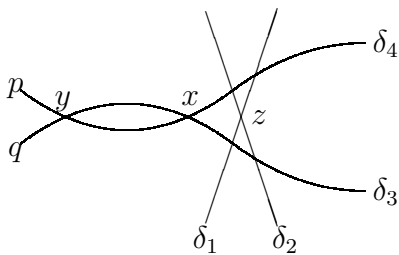
$$\begin{picture}(80,60)  \qbezier(0,0)(40,30)(80,0)\qbezier(80,0)(104,-18)(130,-18)\qbezier(0,20)(40,-10)(80,20)\qbezier(80,20)(104,38)(130,38)
\put(70,-30){\line(1,3){27}}\put(70,50){\line(1,-3){27}}\zput(90,10){$z$}\zput(64,17){$x$}\zput(16,17){$y$}\zput(70,-37){$\delta_1$}\zput(100,-37){$\delta_2$}
\zput(138,-19){$\delta_3$}\zput(138,37){$\delta_4$}\zput(-2,-2){$q$}\zput(-2,22){$p$}\end{picture}$$
\vskip 36pt
In this figure $\delta_i$ have constant height $c_i$ and, during the isotopy, $\delta_3$ and $\delta_4$ are fixed and $\delta_1$ and $\delta_2$
are moving by translation. So the isotopy of the figure is determined by the isotopy of $z$.

If $z$ makes one counterclockwise turn around $x$, we get a movie move of type MVM$_5(c_1,c_2,c_3,c_4)$. If $z$ makes one counterclockwise turn around $y$, we get
a movie move of type MVM$_5(c_1,c_2,c_4,c_3)$.

Let $\Delta$ be the loop containing $x$ and $y$. Suppose that $(c_1,c_2,c_4,c_3)$ is in $S$. Then the move $\varphi$ of type MVM$_5(c_1,c_2,c_3,c_4)$ is equivalent to
the move $\varphi'$ obtained by moving $z$ around $\Delta$. So we have:
$$\varphi'^0=D\circ C\circ B\circ A$$
where $A$ is obtained by moving $\delta_1$ through $\Delta$, $B$ by moving $\delta_2$ through $\Delta$, $C$ by moving back $\delta_1$ trough $\Delta$ and $D$ by
moving back $\delta_2$ through $\Delta$. Denote by $f$ the Reidemeister move of type II which creates crossings $x$ and $y$ and by $g$ its inverse move. So the morphism
$\varphi^0$ is homotopic to:
$$f^0\circ g^0\circ D\circ f^0\circ g^0\circ C\circ f^0\circ g^0\circ B\circ f^0\circ g^0\circ A\circ f^0\circ g^0$$
But each morphism $g^0\circ X\circ f^0$ (with $X=A,B,C,D$) is obtained by a type III Reidemeister move and a movie move of type MVM$_4$:
Let $\psi$ (resp. $\psi'$) be the Reidemeister move obtained by pushing up $z$ through $\delta_4$ (resp. pushing down $z$ through $\delta_3$). If $a,b,c$ are distinct
reals, we set
$\mu(a,b,c)=1$ is $a$ is between $b$ and $c$ and $\mu(a,b,c)=-1$ otherwise. So we have:
$$g^0\circ A\circ f^0\sim\mu(c_1,c_3,c_4)T(\omega^{-1}\otimes\overline\omega)\psi$$
$$g^0\circ B\circ f^0\sim\mu(c_2,c_3,c_4)\psi^{-1}T(\overline\omega^{-1}\otimes\omega)$$
$$g^0\circ C\circ f^0\sim\mu(c_1,c_3,c_4)T(\overline\omega\otimes\omega^{-1})\psi'$$
$$g^0\circ D\circ f^0\sim\mu(c_2,c_3,c_4)\psi'^{-1}T(\omega\otimes\overline\omega^{-1})$$
where $T$ is the map $u\otimes v\mapsto T_p(u)T_q(v)$. Therefore $\varphi^0$ is homotopic to the identity and $(c_1,c_2,c_3,c_4)$ belongs to $S$. So we have:
$$(c_1,c_2,c_3,c_4)\in S\ \ \Longleftrightarrow\ \ (c_1,c_2,c_4,c_3)\in S$$
and $S$ contains all sequences $(c_i)$. So the first part of the lemma is proven and the second one follows easily.\cqfd
\vskip 12pt
\noi{\bf Remark: } The general formula for $\varphi^{\rde K}$ is too complicated to be written down if $\varphi$ is a movie move of type MVM$_5$.
\vskip 12pt
\noi{\bf 4.8 Movie moves of type MVM$_6$.} Let $D$ be a link diagram. We can modify $D$ by adding a small circle near a branch of $D$ and then connect the circle to
$D$ near a point $p$ in the branch. So we get a movie sequence $\varphi=(f,g)$ where $f$ is a surgery move of index $0$ and $g$ a surgery move of index $1$. The
second sequence is $\psi=$Id. The move $(\varphi,\psi)$ is called a movie move of type MVM$_6$ near $p$.
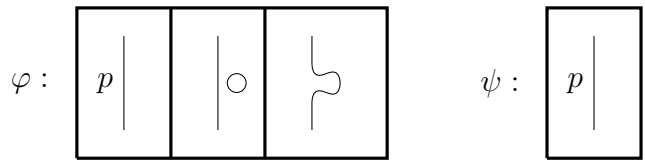
$$\begin{tikzpicture}[scale=1/8] \draw (0,0) -- (0,10) ; \draw (0,5) node[left] {$p$} ; \draw (10,0) -- (10,10) ; \draw (12,5) circle (1) ; 
\draw (20,0) -- (20,3) ; \draw (20,7) -- (20,10) ; \draw [domain=0:10] plot ({21.5-1.5*cos(36*\x)},{3+0.4*\x+0.6*sin(72*\x)}) ;
\draw [very thick] (-5,-3) -- (28,-3) -- (28,13) -- (-5,13) -- (-5,-3) ; \draw [very thick] (5,-3) -- (5,13) ; \draw [very thick] (15,-3) -- (15,13) ;
\draw (40,5) node {$\psi:$} ; \draw (-10,5) node {$\varphi:$} ;
\draw (50,0) -- (50,10) ; \draw (50,5) node[left] {$p$} ; \draw [very thick] (45,-3) -- (55,-3) -- (55,13) -- (45,13) -- (45,-3) ; 
\end{tikzpicture}$$
An easy computation gives the following:
\vskip 12pt
\noi{\bf 4.8.a Lemma: } {\sl Let $(\varphi,\psi)$ be a movie move of type MVM$_6$. Then we have:
$$\varphi^0=\psi^0\hskip 48pt (\overline\varphi)^0=(\overline\psi)^0$$

In the oriented case, suppose that the first move in $\varphi$ is a surgery move of type $(0,a,h)$. then we have:
$$\varphi^1=\psi^1\hskip 48pt (\overline\varphi)^1=(\overline\psi)^1$$
$$\varphi^{\rde K}=\widehat T_p(X(a,h)Y(-a,h))\psi^{\rde K}\hskip 48pt(\overline\varphi)^{\rde K}=(\overline\psi)^{\rde K}\widehat T_p(Z(a,h)Y(a,-h))$$}
\vskip 12pt
\noi{\bf 4.9 Movie moves of type MVM$_7$.} Let $D$ be a link diagram and $e$ be a sign. We have a movie move $(\varphi,\psi)$ described as follows:
$$
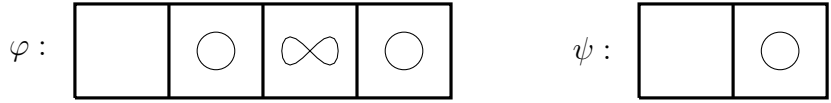
\begin{tikzpicture}[scale=1/8] \draw (5,5) node {$\varphi:$} ; \draw [very thick] (10,0)--(50,0)--(50,10)--(10,10)--(10,0) ; \draw [very thick] (20,0)--(20,10) ;
\draw [very thick] (30,0)--(30,10) ; \draw [very thick] (40,0)--(40,10) ; \draw (65,5) node {$\psi:$} ; \draw [very thick] (70,0)--(90,0)--(90,10)--(70,10)--(70,0) ;
\draw [very thick] (80,0)--(80,10) ;
\draw (25,5) circle (2) ; \draw (45,5) circle (2) ; \draw (85,5) circle (2) ; \draw [domain=0:10] plot ({35+3*cos(36*\x)},{5+1.5*sin(72*\x)}) ;
\end{tikzpicture}$$
We add a small circle to $D$ by a surgery of index $0$. Then we apply a Reidemeister move of type I$_e^+$ which creates a loop in the right part of the circle and we 
remove the other loop in the circle by a Reidemeister move of type I$_e^-$. So we get the movie sequence $\varphi$. The movie sequence $\psi$ is a surgery 
move of index $0$. The move $(\varphi,\psi)$ is called a movie move of type MVM$_7(e)$.

An easy computation gives the following:
\vskip 12pt
\noi{\bf 4.9.a Lemma: } {\sl Let $e$ be a sign and $(\varphi,\psi)$ be a movie move of type MVM$_7(e)$. Let $p$ be a point in the created circle. Then we have:
$$\varphi^0=-e T_p({\omega^{(e)}\over\overline\omega})\psi^0\hskip 48pt (\overline\varphi)^0=e(\overline\psi)^0 T_p({\omega\over\omega^{(e)}})$$

In the oriented case, suppose that the first move in $\varphi$ is a surgery move of type $(0,a,h)$. Then we have:
$$\varphi^1=-e\widehat T_p\Bigl({\omega^{(-eah)}\over\omega^{(ah)}}\Bigr)\psi^1\hskip 48pt (\overline\varphi)^1=e(\overline\psi)^1\widehat T_p\Bigl({\omega^{(eah)}\over
\omega^{(ah)}}\Bigr)$$
$$\varphi^{\rde K}=\widehat T_p\Bigl(X(a,h){A(e,-a,h)\over A(e,a,h)}\Bigr)\varphi^1\hskip 48pt \psi^{\rde K}=\widehat T_p\Bigl(X(-a,h)\Bigr)\psi^1$$
$$(\overline\varphi)^{\rde K}=(\overline\varphi)^1\widehat T_p\Bigl({A(e,a,h)\over A(e,-a,h)}Z(a,h)\Bigr)\hskip 48pt (\overline\psi)^{\rde K}=(\overline\psi)^1
\widehat T_p\Bigl(Z(-a,h)\Bigr)$$}
\vskip 12pt
\noi{\bf 4.10 Movie moves of type MVM$_8$.} Let $D$ be a link diagram and $e$ be a sign. We have a movie move $(\varphi,\psi)$ described as follows:
$$
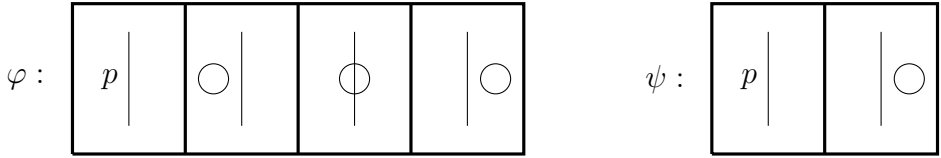
\begin{tikzpicture}[scale=1/8] \draw (-5,8) node {$\varphi:$} ; \draw [very thick] (0,0)--(48,0)--(48,16)--(0,16)--(0,0) ; \draw [very thick] (12,0)--(12,16) ;
\draw [very thick] (24,0)--(24,16);\draw [very thick] (36,0)--(36,16) ; \draw (63,8) node {$\psi: $} ; \draw [very thick] (68,0)--(92,0)--(92,16)--(68,16)--(68,0) ; 
\draw [very thick] (80,0)--(80,16) ;
\draw (6,3) -- (6,13) ; \draw (18,3) -- (18,13) ; \draw (30,3) -- (30,13) ; \draw (42,3) -- (42,13) ; \draw (74,3) -- (74,13) ; \draw (86,3) -- (86,13) ; 
\draw (4,8) node {$p$} ;\draw (72,8) node {$p$}; \draw (15,8) circle (1.6) ; \draw (30,8) circle (1.6) ; \draw (45,8) circle (1.6) ; \draw (89,8) circle (1.6) ; 
\end{tikzpicture}$$
The move $\varphi$ is obtained by creating a small circle of the left hand side of a branch of $D$ with a surgery move of index $0$  and then moving this circle
through the branch with two Reidemeister moves of type II. The move $\psi$ is a surgery move of index $0$ creating the circle on the right hand side of the branch.
This move depends on a sign $e$ where $e=+$ if and only if the circle goes over the branch. Let $p$ be a point in the branch. We say that $(\varphi,\psi)$ is a
movie move of type MVM$_8(e)$ near $p$.
\vskip 12pt
\noi{\bf 4.10.a Lemma:} {\sl Let $D$ be a link diagram and $e$ be a sign. Let $(\varphi,\psi)$ be a movie move of type MVM$_8(e)$ near a point $q\in D$. Let $p$ be a 
point in the created circle. Then we have:
$$\varphi^0=T_p({\omega\over\overline\omega})\psi^0\hskip 48pt (\overline\varphi)^0=-(\overline\psi)^0$$

In the oriented case, suppose that the first move of $\varphi$ is a surgery move of type $(0,a,h)$. Let $b$ be the sign such that $b=+$ if and only if the vertical branch
is oriented downward. Then we have:
$$\varphi^1=\widehat T\Bigl({\omega^{(ah)}\over\omega^{(-ah)}}\otimes 1\Bigr)\psi^1\hskip 48pt (\overline\varphi)^1=-(\overline\psi)^1$$
$$\varphi^{\rde K}=\widehat T\Bigl((X(a,h)\otimes 1){B^e(a,b,-abh)\over B^e(a,-b,-abh)}\Bigr)\varphi^1\hskip 48pt \psi^{\rde K}=\widehat T\Bigl(X(a,-h)\otimes1\Bigr)\psi^1$$
$$(\overline\varphi)^{\rde K}=(\overline\varphi)^1\widehat T\Bigr({B^e(a,-b,-abh)\over B^e(a,b,-abh)}(Z(a,h)\otimes1)\Bigr)\hskip 48pt
(\overline\psi)^{\rde K}=\widehat T\Bigl(Z(a,-h)\otimes1\Bigr)(\overline\psi)^1$$
where $\widehat T$ is the map $u\otimes v\mapsto \widehat T_p(u)\widehat T_q(v)$.}
\vskip 12pt
\noi{\bf Proof:} A straightforward computation gives the desired expressions of $\varphi^0$ and $(\overline\varphi)^0$ and therefore the expressions of $\varphi^1$ and
$(\overline\varphi)^1$.

The orientation of the diagram modified by the first move is described as follows:
$$
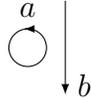
\begin{tikzpicture}[scale=1/8] \draw [->,>=latex] (18,13)--(18,3);\draw (14,8) circle (2);\draw [->,>=latex] (13.7,10)--(13.6,10);\draw (20,4) node {$b$};
\draw (14,12) node {$a$};\end{tikzpicture}$$

The types of the moves in $\varphi$ are:
$$(0,a,h),\ \ \hbox{II}^+(e,a,b,-abh),\ \ \hbox{II}^-(e,a,-b,-abh)$$
and the types of the moves in $\overline\varphi$ are:
$$\hbox{II}^+(e,a,-b,-abh),\ \ \hbox{II}^-(e,a,b,-abh),\ \ (2,a,h)$$
The result follows.\cqfd
\vskip 12pt
\noi{\bf 4.11 Movie moves of type MVM$_9$.} Let $D$ be a link diagram and $e$ be a sign. We have a movie move $(\varphi,\psi)$ described as follows:
$$
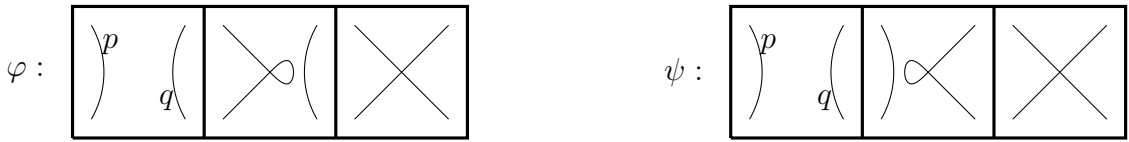
\begin{tikzpicture}[scale=1/8] \draw [very thick] (0,0)--(42,0)--(42,14)--(0,14)--(0,0); \draw [very thick] (14,0)--(14,14); \draw [very thick] (28,0)--(28,14);
\draw [very thick] (70,0)--(112,0)--(112,14)--(70,14)--(70,0); \draw [very thick] (84,0)--(84,14); \draw [very thick] (98,0)--(98,14);\draw (-5,7) node {$\varphi:$};
\draw (65,7) node {$\psi:$}; \draw (2,2) arc (-30:30:10); \draw (12,12) arc (150:210:10); \draw (72,2) arc (-30:30:10); \draw (82,12) arc (150:210:10);
\draw (30,2)--(40,12); \draw (40,2)--(30,12); \draw (100,2)--(110,12); \draw (110,2)--(100,12); \draw (26,12) arc (150:210:10); \draw (86,2) arc (-30:30:10);
\draw (16,2)--(21,7)--(16,12); \draw (96,12)--(91,7)--(96,2); \draw (4,10) node {$p$};\draw (74,10) node {$p$};\draw (10,4) node {$q$};\draw (80,4) node {$q$};
\draw [domain=0:10] plot ({21+2.5*sin(18*\x)},{7-1.25*sin(36*\x)}) ;\draw [domain=0:10] plot ({91-2.5*sin(18*\x)},{7-1.25*sin(36*\x)});\end{tikzpicture}$$
Both moves $\varphi$ and $\psi$ are obtained by a Reidemeister move of type I$_e^+$ followed by a surgery move of index $1$. We say that $(\varphi,\psi)$ is a movie move
of type MVM$_9(e)$ near $(p,q)$.
\vskip 12pt
\noi{\bf 4.11.a Lemma:} {\sl Let $D$ be a link diagram and $e$ be a sign. Let $(\varphi,\psi)$ be a movie move of type MVM$_9(e)$ near $(p,q)$. Then we have:
$$\varphi^0 T_p\Bigl({\omega^{(e)}\over\overline\omega}\Bigr)\sim -e\psi^0\hskip 48pt T_p\Bigl({\omega\over\omega^{(e)}}\Bigr)(\overline\varphi)^0\sim e(\overline\psi)^0$$

In the oriented case, suppose the first move of $\varphi$ is a Reidemeister move of type I$^+(e,a,h)$. Then we have:

$$\varphi^1\widehat T_p\Bigl({\omega^{(eah)}\over\omega^{(-ah)}}\Bigr)\sim -e\psi^1\hskip 48pt\widehat T_p\Bigl({\omega^{(ah)}\over\omega^{(eah)}}\Bigr)
(\overline\varphi)^1\sim e(\overline\psi)^1$$
$$\varphi^{\rde K}\sim\widehat T_p\Bigl(A(e,a,h)Y(-a,h)\Bigr)\varphi^1\hskip 48pt \psi^{\rde K}\sim\widehat T_p\Bigl(A(e,-a,h)Y(a,h)\Bigr)\psi^1$$
$$(\overline\varphi)^{\rde K}\sim(\overline\varphi)^1\widehat T_p\Bigl({Y(a,-h)\over A(e,a,h)}\Bigr)\hskip 48pt
(\overline\psi)^{\rde K}\sim (\overline\psi)^1\widehat T_p\Bigl({Y(-a,-h)\over A(e,-a,h)}\Bigr)$$}
\vskip 12pt
\noi{\bf Proof: } Set: $\e=(1+e)/2$. It is easy to check the following:
$$\varphi^0 T_p(\omega^\e)=-e\psi^0 T_q(\omega^\e)$$
$$T_p(\omega^{1-\e}(\overline\varphi)^0=e T_q(\omega^{1-\e}(\overline\psi)^0$$
So we get:
$$-e\psi^0=\varphi^0 T_p(\omega^\e)T_q(\omega^{-\e})=T_q(\omega^{-\e})\varphi^0 T_p(\omega^\e)\sim T_p(\overline\omega^{-\e})\varphi^0 T_p(\omega^\e)$$
$$=\varphi^0T_p(\overline\omega^{-\e}\omega^\e)=\varphi^0 T_p({\omega^{(e)}\over\overline\omega})$$
$$e(\overline\psi)^0=T_q(\omega^{\e-1})T_p(\omega^{1-\e})(\overline\varphi)^0=T_p(\omega^{1-\e})(\overline\varphi)^0T_q(\omega^{\e-1})\sim 
T_p(\omega^{1-\e})(\overline\varphi)^0T_p(\overline\omega^{\e-1})$$
$$=T_p(\overline\omega^{\e-1}\omega^{1-\e})(\overline\varphi)^0=T_p({\omega\over\overline\omega^{(e)}})(\overline\varphi)^0$$
and we get the first formulae and therefore the following:
$$\varphi^1 T_p\Bigl({\omega^{(e)}\over\overline\omega}\Bigr)\sim -e\psi^1\hskip 48pt T_p\Bigl({\omega\over\omega^{(e)}}\Bigr)(\overline\varphi)^1\sim e(\overline\psi)^1$$

Since the first move of $\varphi$ is a Reidemeister move of type I$^+(e,a,h)$, the orientation of $D$ is as follows:
$$
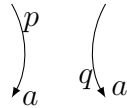
\begin{tikzpicture}[scale=1/8] \draw [->,>=latex] (2,12) arc (30:-30:10);\draw [->,>=latex] (12,12) arc (150:210:10);\draw (4.1,10) node {$p$};\draw (9.9,4) node {$q$};
\draw (4,2) node {$a$};\draw (13.5,3) node {$a$};\end{tikzpicture}$$
and the types of the moves in $\varphi$ (resp. $\psi$) are I$^+(e,a,h)$ and $(1,-a,h)$ (resp. I$^+(e,-a,h)$ and $(1,a,h)$). For $\overline\varphi$ (resp. 
$\overline\psi$), the types are $(1,a,-h)$ and I$^-(e,a,h)$ (resp. $(1,-a,-h)$ and I$^-(e,-a,h)$). The result follows.\cqfd
\vskip 12pt
\noi{\bf 4.12 Movie moves of type MVM$_{10}$.}  Consider a movie move $(\varphi,\psi)$ acting on a link diagram $D$ as follows:
$$
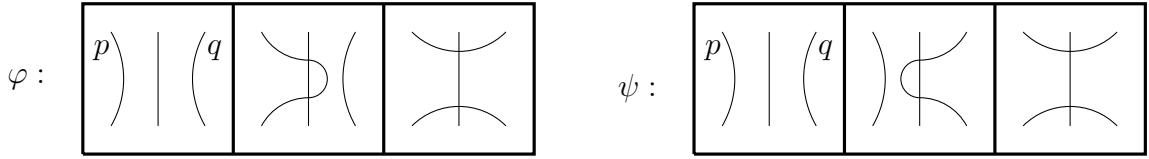
\begin{tikzpicture}[scale=1/8] \draw (-9,5) node {$\varphi:$};\draw (56,4) node {$\psi:$};
\draw (0,0) arc (-30:30:10); \draw (10,10) arc (150:210:10);\draw (5,0)--(5,10);\draw (-1,8) node {$p$};\draw (11,8) node {$q$};
\draw (65,0) arc (-30:30:10); \draw (75,10) arc (150:210:10);\draw (70,0)--(70,10);\draw (64,8) node {$p$};\draw (76,8) node {$q$};
\draw (26,10) arc (150:210:10);\draw (21,0)--(21,10);
\draw (81,0) arc (-30:30:10);\draw (86,0)--(86,10);\draw (21,3) arc (-90:90:2);\draw (21,3) arc (90:151.93:17/3);\draw (21,7) arc (-90:-151.93:17/3);
\draw (32,0) arc (135:45:7.071);\draw (42,10) arc (-45:-135:7.071);\draw (37,0)--(37,10);\draw (86,3) arc (270:90:2);\draw (86,3) arc (90:28.07:17/3);
\draw (86,7) arc (-90:-28.07:17/3);\draw (97,0) arc (135:45:7.071);\draw (107,10) arc (-45:-135:7.071);\draw (102,0)--(102,10);
\draw [very thick] (-3,-3)--(45,-3)--(45,13)--(-3,13)--(-3,-3);\draw [very thick] (13,-3)--(13,13);\draw [very thick] (29,-3)--(29,13);
\draw [very thick] (62,-3)--(110,-3)--(110,13)--(62,13)--(62,-3);\draw [very thick] (78,-3)--(78,13);\draw [very thick] (94,-3)--(94,13);\end{tikzpicture}$$
Both moves $\varphi$ and $\psi$ are obtained by a Reidemeister move of type II$^+$ followed by a surgery move of index $1$. Let $e$ be the sign such that 
$e=+$ (resp. $e=-$) if the middle branch is under (resp. over) the other ones in the modified diagram.
We say that $(\varphi,\psi)$ is a movie move of type MVM$_{10}(e)$ near $(p,q)$. 
\vskip 12pt
\noi{\bf 4.12.a Lemma:} {\sl Let $(\varphi,\psi)$ be a movie move of type MVM$_{10}(e)$ near $(p,q)$. Then we
have:
$$\varphi^0 T_p({\omega\over\overline\omega})\sim\psi^0\hskip 48pt (\overline\varphi)^0=-(\overline\psi)^0$$

In the oriented case, suppose that the first move of $\varphi$ is a Reidemeister move of type II$^+(e,a,b,h)$. Then we have:
$$\varphi^1\widehat T({\omega^{(bh)}\over\omega^{(-bh)}}\otimes1)\sim\psi^1\hskip 48pt (\overline\varphi)^1=-(\overline\psi)^1$$
$$\varphi^{\rde K}\sim\varphi^1\widehat T(B^e(a,b,h)(Y(-a,abh)\otimes 1))\hskip 24pt \psi^{\rde K}\sim\psi^1\widehat T(B^e(a,-b,h)(Y(-a,-abh)\otimes 1))$$
$$(\overline\varphi)^{\rde K}\sim\widehat T\Bigl({Y(a,-abh)\otimes 1\over B^e(a,b,h)}\Bigr)(\overline\varphi)^1\hskip 48pt
(\overline\psi)^{\rde K}\sim\widehat T\Bigl({Y(a,abh)\otimes1\over B^e(a,-b,h)}\Bigr)(\overline\psi)^1$$
where $r$ is a point in the middle branch and $\widehat T$ the map $u\otimes v\mapsto \widehat T_p(u)\widehat T_r(v)$.}
\vskip 12pt
\noi{\bf Proof:} It is easy to check the following:
$$\psi^0=\varphi^0 T_p(\omega) T_q(\omega^{-1})\hskip 48pt (\overline\varphi)^0=-(\overline\psi)^0$$
and that implies:
$$\psi^0=T_p(\omega) T_q(\omega^{-1})\varphi^0\sim T_p(\omega) T_p(\overline\omega^{-1})\varphi^0=\varphi^0 T_p(\omega\overline\omega^{-1})$$

By assumption the diagram $D$ is oriented as follows:
$$
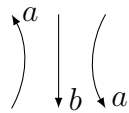
\begin{tikzpicture}[scale=1/8] \draw [->,>=latex] (0,0) arc (-30:30:10); \draw [->,>=latex] (10,10) arc (150:210:10);\draw [->,>=latex] (5,10)--(5,0);
\draw (2,10) node {$a$};\draw (11.5,1) node {$a$}; \draw (6.8,1) node {$b$};
\end{tikzpicture}$$
Therefore the types of the moves in $\varphi$ (resp. $\psi$) are: 

II$^+(e,a,b,h)$ and $(1,b,abh)$ (resp. II$^+(e,a,-b,abh)$ and $(1,-a,-abh)$). 

\noi For the move $\overline\varphi$ (resp. $\overline\psi$), the types are:

$(1,-b,-abh)$ and II$^-(e,a,b,h)$ (resp. $(1,a,abh)$ and II$^-(e,a,-b,abh)$). 

\noi The results follows.\cqfd
\vskip 12pt
\noi{\bf 4.13 Remark:} The conditions induced by all these movie moves are essentially the same as the conditions induced by all the movie moves of Carter and Saito [CS].

Denote by $C(i)$ the condition induced on the correspondance $f\mapsto f^{\rde K}$ by all the movie moves of type MVM$_i$ and by $C'(j)$ the condition induced on it  by
all the movie moves of type MM$_j$ (by using the Bar Natan classification of these moves [BN2]).

The conditions $C'_j$ for $j=1,2,3,4,5$ are automatically satisfied  here because, for every Reidemeister move $f$ with inverse move $g$, $g^{\rde K}$ is a homotopy
inverse of $f^{\rde K}$.

The movie moves of type MVM$_1$, MVM$_2$, MVM$_4$, MVM$_5$, MVM$_6$, MVM$_9$, MVM$_{10}$ are essentially the movie moves of type MM$_7$, MM$_9$, MM$_6$, MM$_{10}$, 
MM$_{11}$, MM$_{13}$ and MM$_{15}$. The movie moves of type MVM$_7$ and MVM$_8$ are more or less equivalent to the movie moves of type MM$_{12}$ and MM$_{14}$, that is
the conditions $C(7)$ and $C(8)$ are exactly the conditions $C'(12)$ and $C'(14)$. Finally the condition $C(3)$ is equivalent, modulo the conditions $C(0)$ and 
$C(4)$, to the condition $C'(8)$.
\vskip 24pt
\noi{\bf 5. Functoriality.}
\vskip 24pt
\noi{\bf 5.1 The category of cobordisms of oriented links.}
\vskip 12pt
Denote by $\pi$ the projection map $(x,y,z)\mapsto (x,y)$ from $\R^3$ onto $\R^2$ and by ${\rde E}$ the space of oriented links in $\R^3$. We say that a link $L\in 
{\rde E}$ is generic if $\pi$ sends $L$ by an immersion onto a link diagram.

The space ${\rde E}$ is a Fr\'echet manifold and the space $Z\subset{\rde E}$ of non generic links is a closed stratified subspace of ${\rde E}$.

Let $L_0$ and $L_1$ be two oriented links. A cobordism $C$ from $L_0$ to $L_1$ is an oriented compact surface $C$ contained in $\R^3\times[0,1]$ and meeting transversally
$\R^3\times\{0,1\}$ in $\partial C=L_1\times\{1\}\cup (-L_0)\times\{0\}$. Two cobordisms $C$ and $C'$ from $L_0$ to $L_1$ are called isotopic if there exist an isotopy 
$f_t:\R^3\times[0,1]\rightarrow\R^3\times[0,1]$, with  $0\leq t\leq 1$, such that:
$$f_0=\hbox{Id}\hskip 48pt f_1(C)=C'$$
$$\forall t\in[0,1],\ \  f_t\hbox{\ is the identity on } \R^3\times\{0,1\}$$

Let $L_0$ and $L_1$ be two generic oriented links and $D_0$ and $D_1$ be the corresponding link diagrams. Let $f: D_0\rightarrow D_1$ be an elementary move. This move
induces a cobordism $C$ from $L_0$ to $L_1$ which is well defined up to isotopy. Such cobordisms are called elementary cobordisms.

Let $L_0$ and $L_1$ be two generic oriented links and $C$ be a cobordism from $L_0$ to $L_1$. We say that $C$ is generic if there is finitely many elements $t_i\in(0,1)$
such that: $C\cap \R^3\times\{t\}$ is generic for all $t\in[0,1]$ except the $t_i$'s and the cobordism $C\cap\R^3\times[t_i-\e,t_i+\e]$ is elementary for each $i$ and
$\e$ small enough.

By transversality we see that every link $L\in{\rde E}$ is isotopic to a generic link by a small isotopy. Moreover every cobordism between two generic links in ${\rde E}$
is isotopic to a generic cobordism by a small isotopy.

Denote by ${\rde L}$ the category of cobordisms of oriented links in $\R^3$. The objects of ${\rde L}$ is the links in ${\rde E}$ and the morphisms are the isotopy 
classes of cobordisms. Every morphism in ${\rde L}$ from a generic link $L_0$ to a generic link $L_1$ is a composite of elementary morphisms and can be described by
a movie sequence.
\vskip 12pt
\noi{\bf 5.2 Definition: } Let ${\rde K}$ be a Khovanov data. This data is said to be functorializable if there exists a monoidal functor $\Psi$ from the category of 
cobordisms of oriented links ${\rde L}$ to the homotopy category of $K$-complexes satisfying the following:

--- if $D$ is the diagram of a generic oriented link $L$, we have: $\Psi(L)=KH(D)$

--- if $C$ is an elementary cobordism from $L_0$ to $L_1$ associated to an elementary move $f:D_0\rightarrow D_1$, $\Psi(C)$ is homotopic to the morphism $f^{\rde K}$.

If these conditions are satisfied the functor $\Psi$ will be called a Khovanov functor associated to ${\rde K}$.
\vskip 12pt
\noi{\bf 5.3 Lemma: } {\sl Let ${\rde K}$ be a Khovanov data. Suppose ${\rde K}$ is functorializable. Then all the Khovanov functors associated to ${\rde K}$ are 
isomorphic.}
\vskip 12pt
\noi{\bf Proof:} Let ${\rde L}_0$ be the full subcategory of ${\rde L}$ generated by generic links. It is clear that the inclusion ${\rde L}_0\subset{\rde L}$ is a
equivalence of categories. The result follows.\cqfd
\vskip 12pt
\noi{\bf 5.4 Lemma:} {\sl Let ${\rde K}$ be a Khovanov data. Then ${\rde K}$ is functorializable if and only if all conditions $C(i)$, for $i=0,\dots,10$ are satisfied
for the correspondance $f\mapsto f^{\rde K}$.}
\vskip 12pt
\noi{\bf Proof:} For every movie move $(\varphi,\psi)$, the condition $\varphi^{\rde K}\sim\psi^{\rde K}$ will be denoted by ${\rde K}(\varphi,\psi)$.
It is clear that, if the correspondance $f\mapsto f^{\rde K}$ comes from a functor, every movie move induces a trivial condition. So all conditions $C(i)$ are
satisfied.

Conversely, suppose all conditions $C(i)$, for $i=0,\dots,10$ are satisfied. The only thing to do is to define the functor
$\Psi$ on the category of cobordisms of generic oriented links ${\rde L}_0$. This functor is well defined on the objects. 

Let $L_0$ and $L_1$ be two generic oriented links and $C_0$ be a cobordism from $L_0$ to $L_1$. This cobordism represents a morphism $f$ from $L_0$ to $L_1$. We have
to define $\Psi(f)$. Since $C_0$ is isotopic to a generic cobordism, $\Psi(f)$ is homotopic to a morphism $\varphi^{\rde K}$, where $\varphi$ is some movie sequence. So we
don't have any choice for $\Psi(f)$. The only thing to do is to prove that $\varphi^{\rde K}$ depends only on the isotopy class of $C_0$.

Let $M={\rde E}(C_0)$ be the space of cobordisms isotopic to $C_0$ and $Z_0=Z(C_0)$ be the space of cobordisms in $M$ which are not generic. The space $M$ is
a connected Fr\'echet manifold and $Z_0$ is a closed stratified subspace of $M$ of codimension $1$. Every $C\in M$ which is not in $Z_0$ is determined by a
movie sequence $\varphi=\widehat C$, and $\Psi(C)=\varphi^{\rde K}$ is well defined. The last thing to do is to prove the following: if $C$ and $C'$ are in 
$M\setminus Z_0$, $\Psi(C)$ and $\Psi(C')$ are homotopic.

Let $C$ be a cobordism in $M$. For every $t\in[0,1]$ denote by $C_t$ the intersection $C\cap\R^3\times\{t\}$. We say that $t$ is $C$-regular if $C$ is transverse
to $\R^3\times\{t\}$ and $C_t$ is a generic link in $\R^3\times\{t\}\simeq\R^3$. We say that $t$ is $C$-singular if $t$ is not $C$-regular.

If $C\in M$ is generic, there is only finitely many $C$-singular elements in $[0,1]$, and every $C$-singular parameter $t$ correspond to a Reidemeister move if $C$ is
transverse to $\R^3\times\{t\}$ and a surgery move if it is not the case.

Since $M$ is connected, there is a path $\gamma$ in $M$ joining $C$ and $C'$. Up to modify $\gamma$ we may as well suppose that $\gamma$ is transverse to $Z_0$ and we
have to prove that $\Psi(\gamma(s))$ doesn't change when $\gamma(s)$ goes through $Z_0$. Let $s_0$ be a parameter such that $\gamma(s_0)$ is in $Z_0$. Since $\gamma$
is transverse to $Z_0$, $\gamma(s_0)$ is in a stratum $S\subset Z_0$ of codimension $1$. For $\e$ small enough, $\gamma(s_0-\e)$ and $\gamma(s_0+\e)$ correspond to
movie sequences $\varphi$ and $\psi$ and we get a movie move $(\varphi,\psi)$. So every codimension $1$ stratum of $Z_0$ induces a movie move $(\varphi,\psi)$. The only
condition to check is the conditions ${\rde K}(\varphi,\psi)$ for all these movie moves.

Let $S$ be a codimension $1$ stratum in $Z_0$. An element in $S$ is a cobordism $C$ which is not generic but for only one reason. So there is a unique $t_0\in(0,1)$ such
that $C\cap (\R^3\times[t_0-\e,t_0+\e])$ is not generic. Let $(\varphi,\psi)$ be the movie move associated to $C$. The cobordism $C\cap (\R^3\times[t_0-\e,t_0+\e])$
induces a movie move $(\varphi_0,\psi_0)$ and there are two movies sequences $\alpha$ and $\beta$ such that: $\varphi=\alpha\circ\varphi_0\circ\beta$ and 
$\psi=\alpha\circ\psi_0\circ\beta$. So the condition ${\rde K}(\varphi_0,\psi_0)$ implies the condition ${\rde K}(\varphi,\psi)$ and it is enough to consider the case
where $t_0$ is the only $C$-singular element in $[0,1]$ and $C_t$ is generic for every $t\not=t_0$.

Suppose $C$ is transverse to $\R^3\times\{t_0\}$. Then the problem reduces to an isotopy problem and then to an isotopy corresponding to a loop around a codimension $2$
stratum of $Z\subset{\rde E}$. These loops correspond to all movie moves of type MVM$_i$ for $0\leq i\leq 5$ involving only Reidemeister moves.

Suppose $C$ is not transverse to $\R^3\times\{t_0\}$. Then the function $C\subset \R^3\times[0,1]\rightarrow[0,1]$ has critical points $u_i$ is $\R^3\times\{t_0\}$.
If one of these points is degenerated, the problem reduces to a movie move of type MVM$_6$. If there is at least $2$ critical points, the problem reduces to a movie
move of type MVM$_0$ involving two surgery moves. In the other cases, we have only one critical point $u$ and this point is non degenerated with index $d$. But we have
an extra condition because $C$ belongs to $Z_0$. There is two possibilities for this condition: $u$ is a critical point for $\pi: C\rightarrow \R^2$ or $\pi(u)$ is a
multiple point in $\pi(C_{t_0})$. In the first case, we get a movie move of type MVM$_7$ if $d=0$ or $d=2$ and a movie move of type MVM$_9$ if $d=1$.
In the second case, we get a movie move of type MVM$_0$ involving one Reidemeister move and one surgery move or a movie move of type MVM$_8$ if $d=0$ or $d=2$
and a movie move of type MVM$_{10}$ if $d=1$.

Thus, if all conditions $C(i)$ are satisfied, the functor $\Psi$ is well defined.\cqfd
\vskip 12pt
For simplicity we'll identify the sign $+$ with $1$ and the sign $-$ with $-1$ and we define a map $<?|?>$ from $\{\pm\}^2$ to $\{\pm\}$, a map $\mu$ from $\{\pm\}$ to
$R^*$ and a map $\delta$ from $\{\pm\}^2$ to $R^*$ by:
$$<(-1)^p|(-1)^q>=(-1)^{pq}$$
$$\mu(e)=\omega^{(1+e)/2}=\left\{\matrix{\omega&\hbox{if}\ e=+\cr 1&\hbox{if}\ e=-\cr}\right.$$
$$\delta(a,b)=\left\{\matrix{\omega\overline\omega&\hbox{if}\ a=b=+\cr 1&\hbox{otherwise}\cr}\right.$$
The map $<?|?>$ is symmetric and satisfy the following property:
$$\forall a,b,c\in\{\pm\},\ \ <a|bc>=<a|b><a|c>$$
If $H(h)$ is an element of $R^*$ depending on a sign $h$, we set:
$$\widetilde H={H(+)\over H(-)}$$
\vskip 12pt
\noi{\bf 5.5 Theorem:} {\sl Let ${\rde K}=(A,B,C,X,Y,Z)$ be a Khovanov data. Then ${\rde K}$ is functorializable if and only if there exist elements $\sigma(h)$ in $K^*$, $E_e(h)$,
$F(e)$ in $R^*$, $W(h)$ in $(R\otimes R)^*$ depending on signs $e$ and $h$ such that the following holds, for every signs $e,h,a,b,c$:
$$\sigma(+)\sigma(-)=\omega\overline\omega$$
$$E_+(h)E_-(h)=\omega F(+)F(-)U(h)V(h)$$
$$A(e,a,h)=<e|a><e|h><a|h>\Bigl({\sigma(h)\over\omega^{(a)}}\Bigr)^{(1-eh)/2}{E_e(h)\over F(a)}$$
$$B(a,b,h)={<a|b>\over\delta(ah,bh)}\sigma(h)^{(1+ab)/2}U(h)\mu(bh)\otimes V(h)\mu(ah)$$
$$X(a,h)=<a|-h>\Bigl({\sigma(h)\over\omega^{(a)}}\Bigr)^{(1+h)/2}F(-a)$$
$$Y(a,h)=X(-a,h)^{-1}\hskip 48pt Z(a,h)=X(-a,-h)$$
$$C(e,a,b,c)=-ac(1\otimes \bigl(\widetilde E_{-e}\bigr)^{-1}H^{(a+b)(b+c)/4}\otimes 1)\widehat U^{a(b+c)/2}\widehat V^{(a+b)c/2}$$
with: $W(h)=U(h)\otimes V(h)$, $H=\widetilde\sigma\bigl(\widetilde U\widetilde V\bigr)^2$, $\widehat U=\widetilde U\otimes\widetilde U^{-1}\otimes 1$,
$\widehat V=1\otimes\widetilde V^{-1}\otimes\widetilde V$.}
\vskip 12pt
\noi{\bf Remark:} Suppose  ${\rde K}$ is functorializable. Denote by $\Psi$ the associated functor. The system $(\sigma,E,F,W)$ will be called a parametrization of
the Khovanov functor $\Psi$. It is easy to see that such a parametrization is unique.
\vskip 12pt
\noi{\bf 5.6 Remark:} In the classical case, $\omega$ is equal to $1$ and the Khovanov data given by:
$$A(e,a,h)=<e|a><e|h><a|h>\hskip 24pt B(a,b,h)=<a|b>\hskip 24pt C(e,a,b,c)=-ac$$
$$X(a,h)=<a|-h>\hskip 24pt Y(a,h)=<-a|-h>\hskip 24pt Z(a,h)=<-a|h>$$
is functorializable.
\vskip 12pt
\noi{\bf Proof: } Because of lemma 4.3.a, the condition $C(2)$ is equivalent to:
$${B^e(-a,-b,h)\over B^e(a,b,h)}=-ab\omega^{(-bh)}\otimes{1\over\omega^{(ah)}}$$
for every signs $e,a,b,h$. It is easy to see that these conditions are equivalent to:
$${B(-a,-b,h)\over B(a,b,h)}=-ab\omega^{(-bh)}\otimes{1\over\omega^{(ah)}}$$
for every $a,b,h$. Define the elements $B'(a,b,h)$ in $(R\otimes R)^*$ by:
$$B(a,b,h)={<a|b>\over\delta(ah,bh)}\bigl(\mu(bh)\otimes\mu(ah)\bigr)B'(a,b,h)$$
With these new elements, the condition $C(2)$ is equivalent to:
$$B'(-a,-b,h)=B'(a,b,h)$$
and $C(2)$ is equivalent to the fact that $B'(a,b,h)$ depends only on $ab$ and $h$: 
$$B'(a,b,h)=B''(ab,h)$$

Because of lemma 4.10.a, the condition $C(8)$ is equivalent to:
$${B^e(a,b,-abh)\over B^e(a,-b,-abh)}=\Bigl({X(a,-h)\over X(a,h)}{\omega^{(-ah)}\over\omega^{(ah)}}\Bigr)\otimes1$$
$${B^e(a,b,-abh)\over B^e(a,-b,-abh)}=-\Bigl({Z(a,h)\over Z(a,-h)}\Bigr)\otimes1$$
These conditions imply the following:
$${B(a,b,h)\over B(a,-b,h)}\in R^*\otimes1\hskip 48pt{B(a,b,h)\over B(-a,b,h)}\in 1\otimes R^*$$
and that's equivalent to the fact that
$${B''(e,h)\over B''(-e,h)}$$
lies in $R^*\otimes 1$ and in $1\otimes R^*$ and therefore in $K^*(1\otimes 1)$.

Set: 
$$W(h)=B''(-,h)\hskip 24pt \sigma(h)=B''(+,h)W(h)^{-1}$$
The elements $\sigma(h)$ are in $K^*$ and we have:
$$B(a,b,h)={<a|b>\over\delta(ah,bh)}\sigma(h)^{(1+ab)/2}(\mu(bh)\otimes\mu(ah))W(h)$$
With this expression, we have:
$${B^e(a,b,-abh)\over B^e(a,-b,-abh)}=a\sigma(-abh)^{ab}\bigl(\omega^{(bh)}\bigr)^{-ab}\otimes 1$$
and the condition $C(8)$ is equivalent to:
$$a\sigma(-abh)^{ab}\bigl(\omega^{(bh)}\bigr)^{-ab}={X(a,-h)\over X(a,h)}{\omega^{(-ah)}\over\omega^{(ah)}}=-{Z(a,h)\over Z(a,-h)}$$
So the left hand side term is independant of $b$ and we get: $\sigma(+)\sigma(-)=\omega\overline\omega$. Using that, the condition $C(8)$ is equivalent to:
$${X(a,h)\over X(a,-h)}=a{\sigma(h)\over\omega^{(ah)}}\hskip 48pt {Z(a,h)\over Z(a,-h)}=-a{\sigma(-h)\over\omega^{(ah)}}$$

Set: $F(a)=X(-a,-)$. Then we have:
$$X(a,h)=<a|-h>\Bigl({\sigma(h)\over\omega^{(a)}}\Bigr)^{(1+h)/2}F(-a)$$

The condition $C(6)$ is equivalent to:
$$X(a,h)Y(-a,h)=1\hskip 48pt Z(a,h)Y(a,-h)=1$$
and $B,X,Y,Z$ can be describe in term of $\sigma,W,F$. Using these descriptions one can check that all conditions $C(i)$, for $i=0,2,6,8,10$, are satisfied.

Define the elements $A'(e,a,h)$ by:
$$A(e,a,h)=<e|a><ea|h>\Bigl({\sigma(h)\over\omega^{(a)}}\Bigr)^{(1-eh)/2}{A'(e,a,h)\over F(a)}$$
Using these new elements, the condition $C(7)$ becomes: $A'(e,a,h)=A'(e,-a,h)$ and $A'(e,a,h)$ depends only on $(e,h)$. So by setting: $A'(e,a,h)=E_e(h)$, we have:
$$A(e,a,h)=<e|a><ea|h>\Bigl({\sigma(h)\over\omega^{(a)}}\Bigr)^{(1-eh)/2}{E_e(h)\over F(a)}$$

Using that, all conditions $C(i)$, for $i=0,2,6,7,8,9,10$, are now satisfied. The condition $C(1)$ becomes the relation:
$$E_+(h)E_-(h)=\omega F(+)F(-)U(h)V(h)$$
with: $U(h)\otimes V(h)=W(h)$ and the last thing to do is to compute the elements $C(e,a,b,c)$ and verify the conditions $C(i)$, for $i=3,4,5$.

Because of lemma 4.4.a, the condition $C(3)$ is equivalent to:
$$C(e,-ah,-ah,bh)=u\otimes v\otimes w\ \Longrightarrow\ uv\otimes w=ab\widetilde\sigma^{(1-ab)/2}\widetilde E_e\otimes 1(\widetilde U\otimes\widetilde V)^{-ab}$$
$$C(e,bh,-ah,-ah)=u\otimes v\otimes w\ \Longrightarrow\ wv\otimes u=ab\widetilde\sigma^{(1-ab)/2}\widetilde E_e\otimes 1(\widetilde V\otimes\widetilde U)^{-ab}$$
In the case: $b=-a$, we get (with: $C(e,-ah,-ah,-ah)=u\otimes v\otimes w$):
$$uv\otimes w=-\widetilde\sigma(\widetilde E_e\otimes 1)(\widetilde U\otimes\widetilde V)$$
$$wv\otimes u=-\widetilde\sigma(\widetilde E_e\otimes 1)(\widetilde V\otimes\widetilde U)$$
and that implies:
$$C(e,-ah,-ah,-ah)=-\widetilde U\otimes\widetilde\sigma\widetilde E_e\otimes\widetilde V$$
and then:
$$C(e,a,a,a)=-\widetilde U\otimes\widetilde\sigma\widetilde E_e\otimes\widetilde V$$

Consider the condition $C(4)$ given by lemma 4.5.a. This condition depends on signs $a,b,c,h$ and on an element $e$ in $\{+,-,0\}$. We consider this condition
in the case: $a=b=1$.

If $e=-$, this condition is the following:
$$C(+,h,h,ch)C(-,-h,-h,ch)=-\widetilde\sigma(\widetilde U\otimes\widetilde V\otimes 1)$$
or:
$$C(c,ch,ch,ch)C(-c,-ch,-ch,ch)=-\widetilde\sigma(\widetilde U\otimes\widetilde V\otimes 1)$$
$$\Longrightarrow C(-c,-ch,-ch,ch)=1\otimes {\widetilde V\over\widetilde E_c}\otimes{1\over\widetilde V}$$
So we get:
$$C(e,a,a,-a)=1\otimes{\widetilde V\over\widetilde E_{-e}}\otimes{1\over\widetilde V}$$

If $e=0$ the condition is:
$$C(+,h,ch,h)C(-,-h,ch,-h)=\widetilde\sigma(\widetilde U\otimes1\otimes\widetilde V)$$
or:
$$C(c,ch,ch,ch)C(-c,-ch,ch,-ch)=\widetilde\sigma(\widetilde U\otimes 1\otimes\widetilde V)$$
$$\Longrightarrow C(-c,-ch,ch,-ch)=-1\otimes{1\over \widetilde E_c}\otimes 1$$
and we have:
$$C(e,a,-a,a)=-1\otimes{1\over \widetilde E_{-e}}\otimes 1$$

If $e=+$ the condition is:
$$C(+,ch,h,h)C(-,ch,-h,-h)=-\widetilde\sigma(1\otimes\widetilde U\otimes\widetilde V)$$
or:
$$C(c,ch,ch,ch)C(-c,ch,-ch,-ch)=-\widetilde\sigma(1\otimes\widetilde U\otimes\widetilde V)$$
$$\Longrightarrow C(-c,ch,-ch,-ch)={1\over\widetilde U}\otimes{\widetilde U\over\widetilde E_c}\otimes 1$$
and we have:
$$C(e,-a,a,a)={1\over\widetilde U}\otimes{\widetilde U\over\widetilde E_{_e}}\otimes 1$$
So we get the general formula:
$$C(e,a,b,c)=-ac(1\otimes \bigl(\widetilde E_{-e}\bigr)^{-1}H^{(a+b)(b+c)/4}\otimes 1)\widehat U^{a(b+c)/2}\widehat V^{(a+b)c/2}$$
with: $H=\widetilde\sigma\bigl(\widetilde U\widetilde V\bigr)^2$, $\widehat U=\widetilde U\otimes\widetilde U^{-1}\otimes 1$,
$\widehat V=1\otimes\widetilde V^{-1}\otimes\widetilde V$. 

Using this expression, it is easy to check the conditions $C(3)$ and $C(4)$ and the last condition to check is $C(5)$.

Consider a movie move of type MVM$_5(c_1,c_2,c_3,c_4)$, where the $c_i$'s are distinct reals. This move involves four lines $\delta_1,\delta_2,\delta_3,\delta_4$ and
each $c_i$ is the height of $\delta_i$. Since we are consider oriented links, each line has to be oriented. Let $a,b,c,d$ be four signs and suppose the singular link
diagram is oriented as follows:
$$
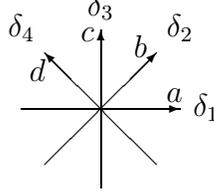
\begin{picture}(60,60) \put(30,0){\vector(0,1){60}}\put(0,30){\vector(1,0){60}}\put(8.79,8.79){\vector(1,1){42.42}}\put(51.21,8.79){\vector(-1,1){42.42}}
\zput(70,30){$\delta_1$}\zput(60,60){$\delta_2$}\zput(30,68){$\delta_3$}\zput(0,60){$\delta_4$}\zput(58,35){$a$}\zput(45,52){$b$}\zput(25,58){$c$}\zput(6,43){$d$}
\end{picture}$$
This figure induces a movie move of type MVM$_5(c_1,c_2,c_3,c_4)$ and then an element $M(a,b,c,d,h)\in(R^4)^*$ depending on the signs $a,b,c,d$ and the $D$-sign of
the center of the triangle corresponding to the first Reidemeister move. So we have to prove that $M(a,b,c,d,h)$ is allways equal to $1$.

First of all, consider the case where the sequence $(c_1,c_2,c_3,c_4)$ is decreasing.

For each $k$ in $\{1,2,3,4\}$ denote by $F_k:R^{\otimes3}\rightarrow R^{\otimes4}$ the tensorization by $1$ at the $k$-th position. For example, $F_2$ is the map:
$$u\otimes v\otimes w\mapsto u\otimes 1\otimes v\otimes w$$

We have the following:
$$M(a,b,c,d,h)=F_4(C(-,ah,-bh,ch)^h)F_3(C(-,-ah,bh,-dh)^{-h})\times$$
$$F_2(C(-,ah,-ch,dh)^h)F_1(C(-,-bh,ch,-dh)^{-h})\times$$
$$F_4(C(-,-ah,bh,-ch)^h)F_3(C(-,ah,-bh,dh)^{-h})\times$$
$$F_2(C(-,-ah,ch,-dh)^h)F_1(C(-,bh,-ch,dh)^{-h})$$
But $C(e,ah,bh,ch)$ doesn't depend on $h$. So we have:
$$M(a,b,c,d,h)^h=M^2$$
with:
$$M={F_4(C(-,a,-b,c))F_2(C(-,a,-c,d))\over F_3(C(-,a,-b,d))F_1(C(-,b,-c,d))}$$
A straightforward computation gives the following:
$$M=\Bigl(1\otimes H^{(a-b)(c-d)/4}\otimes H^{-(a-b)(c-d)/4}\otimes 1\Bigr)\Bigl(1\otimes \widetilde U^{-(a-b)(c-d)/2}\otimes\widetilde U^{(a-b)(c-d)/2}\otimes1\Bigr)\times$$
$$\Bigl(1\otimes\widetilde V^{-(a-b)(c-d)/2}\otimes\widetilde V^{(a-b)(c-d)/2}\otimes1\Bigr)=1\otimes X\otimes X^{-1}\otimes1$$
with:
$$X=H^{(a-b)(c-d)/4}\widetilde U^{-(a-b)(c-d)/2}\widetilde V^{-(a-b)(c-d)/2}=\Bigl(H\widetilde U^{-2}\widetilde V^{-2}\Bigr)^{(a-b)(c-d)/4}=\widetilde\sigma^{(a-b)(c-d)/4}$$
Therefore $M$ is equal to $1$ and the condition $C(5)$ is allways satisfied when the sequence $(c_i)$ is decreasing. For a general sequence, we use the fact that the
condition $C(4)$ is allways satisfied and we check the condition with exactly the same proof as the proof of lemma 4.6.a. So all the conditions $C(i)$ are satisfied and
the theorem is proven. Therefore theorem A is also proven.\cqfd
\vskip 12pt
Let $(\sigma,E,F,W)$ be a parametrization of a Khovanov functor $\Psi$. Consider the following element in $R^*$:
$$\pi=-{\omega F(+)F(-)\over\sigma(-)}$$
This element is called the weight of $\Psi$. By construction, every invertible element in $R$ is the weight of a Khovanov functor. For example the weight of the functor
described in remark 5.6 is $-1$.
\vskip 12pt
\noi{\bf 5.7 Proposition:} {\sl Let $\Psi$ be a Khovanov functor and $\pi$ be its weight. For every integer $p\geq0$, denote by $\Sigma_p$ an 
unknotted oriented surface of genus $p$ in $\R^4$. Then we have:
$$\Psi(\Sigma_p)=\e(\delta^p\pi^{1-p})$$
$$\sum_{p\geq0} x^p \Psi(\Sigma_p)=\e\Bigl({\pi\over1-x\delta/\pi}\Bigr)={\e(\pi)+x(2-\e(\pi)^2u)\over 1-x\e(\pi)u+x^2u}\in K[[x]]$$
with: $u=-\omega\overline\omega(\alpha-\overline\alpha)^2/(\pi\overline\pi)=\delta\overline\delta/(\pi\overline\pi)$.}
\vskip 12pt
\noi{\bf 5.8 Remark:} For $\pi=1$, this formula is exactly the same as the formula in the lemma 1.6. Actually, the right hand side part of this formula is
the image of the corresponding part in lemma 1.6 by the morphism sending $\e(1)$ to $\e(\pi)$ and $\delta$ to $\delta/\pi$.
\vskip 12pt
\noi{\bf Proof:} Let $(A,B,C,X,Y,Z)$ be a Khovanov data such that $\Psi$ is the corresponding functor. Let $p\geq0$ be an integer. An unknotted surface of
genus $p$ in $\R^3\times[0,1]$ can be describe by the following movie sequence:
$$\varphi=(f,g,g',\dots,f')$$
where $(g,g')$ is repeated $p$ times. In this sequence, $f$, $g$, $g'$ and $f'$ are surgery moves of type $(0,-,+)$, $(1,-,-)$, $(1,+,+)$ and $(2,-,+)$. So we get:
$$\Psi(\Sigma_p)=\e\Bigl(X(-,+)Z(-,+)\Bigl(Y(-,-)Y(+,+)\omega(\alpha-\overline\alpha)\Bigr)^p\Bigr)$$
But it is easy to check the following:
$$X(-,+)Z(-,+)=X(-,+)X(+,-)=\pi$$
$$Y(-,-)Y(+,+)=\Bigl(X(+,-)X(-,+)\Bigr)^{-1}={1\over\pi}$$
So we have:
$$\Psi(\Sigma_p)=\e(\delta^p\pi^{1-p})$$
$$\sum_{p\geq0}x^p\Psi(\Sigma_p)=\e\Bigl({\pi\over 1-x\delta/\pi}\Bigr)={\e(\pi-x\overline\delta\pi/\overline\pi)\over(1-x\delta/\pi)(1-x\overline\delta/\overline\pi)}
={\e(\pi-x\overline\delta\pi/\overline\pi)\over1-x\e(\pi)\sigma\overline\sigma+x^2\sigma\overline\sigma}$$
with: $\sigma=\delta/\pi$. So we get the desired formula.\cqfd
\vskip 12pt
\noi{\bf 5.9 Theorem:} {\sl Two Khovanov functors with the same weight are isomorphic.}
\vskip 12pt
\noi{\bf Proof:} Consider two Khovanov functors $\Psi$ and $\Psi'$ with the same weight $\pi$. They are described by two Khovanov data $(A,B,C,X,Y,Z)$ and
$(A',B',C',X',Y',Z')$. Since these two Khovanov data are functorializable there exist elements $\sigma(h)$ in $K^*$, $E_e(h),F(e)$ in $R^*$ and 
$W(h)$ in $(R\otimes R)^*$ such that:
$$\sigma(+)\sigma(-)=1$$
$$E_+(h)E_-(h)=F(+)F(-)U(h)V(h)$$
$$A'(e,a,h)=A(e,a,h)\sigma(h)^{(1-eh)/2}{E_e(h)\over F(a)}$$
$$B'(a,b,h)=B(a,b,h)\sigma(h)^{(1+ab)/2}W(h)$$
$$X'(a,h)=X(a,h)\sigma(h)^{(1+h)/2}F(-a)$$
$$C'(e,a,b,c)=C(e,a,b,c)(1\otimes \widetilde E_{-e}^{-1}H^{(a+b)(b+c)/4}\otimes 1)\widehat U^{a(b+c)/2}\widehat V^{(a+b)c/2}$$
with: $W(h)=U(h)\otimes V(h)$, $H=\widetilde\sigma\bigl(\widetilde U\widetilde V\bigr)^2$, $\widehat U=\widetilde U\otimes\widetilde U^{-1}\otimes 1$,
$\widehat V=1\otimes\widetilde V^{-1}\otimes\widetilde V$.

Because of the first relation, there exists an element $\sigma\in K^*$ with: $\sigma(h)=\sigma^h$. So we have:
$$A'(e,a,h)=A(e,a,h)\sigma^{(h-e)/2}{E_e(h)\over F(a)}$$
$$B'(a,b,h)=B(a,b,h)\sigma^{h(1+ab)/2}W(h)$$
$$X'(a,h)=X(a,h)\sigma^{(1+h)/2}F(-a)$$
$$H=(\sigma\widetilde U\widetilde V)^2$$
Moreover, since $\Psi$ and $\Psi'$ have the same weight, we have also:
$$F(+)F(-)\sigma=1$$
So, by setting: $F=F(+)$, we have:
$$F(e)=F^e\sigma^{(e-1)/2}$$

Consider elements $x_\e(e,h)$ in $R^*$ depending on signs $\e$, $e$ and $h$.

Let $D$ be an oriented link diagram and $\widehat D$ be the oriented resolution of $D$. Denote by $d(D)$ the winding number of $D$. For each component 
$\widehat C$ of $\widehat D$ denote by $g'(\widehat C)$ the $(\widehat D\minus\widehat C)$-sign of a point in $\widehat C$. Denote also by $g(D)$ the sum of all 
$g'(\widehat C)$. It is clear that $d(D)$ and $g(D)$ are both congruent to the number of components of $\widehat D$ mod $2$ and $d(D)-g(D)$ is even.

For each signs $e$ and $h$, denote by $X(e,h)$ the set of crossings of $D$ with sign $e$ and $D$-sign $h$. If $C$ is a component of $D$, denote by $N_+(e,h,C)$
(resp. $N_-(e,h,C)$) the number of crossing $x$ in $X(e,h)$ such that the over branch (resp. the under branch) containing $x$ is in $C$.

Set:
$$G(D)=\sigma^{(d(D)-g(D))/2}\build\otimes_C^{}\Bigl(F^{d(C)}\prod_{\e,e,h}x_\e(e,h)^{N_\e(e,h,C)}\Bigr)$$
Let $D_0$ be the set of components of $D$. The element $G(D)$ belongs to $R^{\otimes D_0}$ and induces, via the maps $\widehat T$, an automorphism 
$A(D):KH(D)\rightarrow KH(D)$ well defined up to homotopy.

By conjugation with these automorphisms, the functor $\Psi'$ is transformed into a new fonctor $\Psi''$. So, for each morphism $f:D\rightarrow D'$, we have a
diagram which is commutative up to homotopy:
$$\diagram{KH(D)&\hfl{\Psi'(f)}{}&KH(D')\cr \vfl{A(D)}{}&&\vfl{A(D')}{}\cr KH(D)&\hfl{\Psi''(f)}{}&KH(D')\cr}$$
This new functor is clearly isomorphic to $\Psi'$. Let's choose the elements $x_*(*,*)$ such that:
$$x_+(+,h)x_+(-,h)={1\over U(h)}\hskip 24pt x_-(+,h)x_-(-,h)={1\over V(h)}$$
$$x_+(e,h)x_-(e,h)={\sigma^{(e-1)/2}\over E_e(h)}$$
It is easy to see that this choice is possible. In this case, a straightforward computation shows that $\Psi$ and $\Psi''$ agree on every elementary move and therefore
on the category ${\rde L}_0$. Thus $\Psi'$ is isomorphic to $\Psi$ and theorem 5.9 (and theorem B) is proven.\cqfd
\vskip 12pt
\noi{\bf 5.10 Remark:} Consider two isomorphic Khovanov functors of weight $\pi$ and $\pi'$. Because of Proposition
5.7, we have: $\e(\pi')=\e(\pi)$. In the case $R=R_0$, that implies: $\pi'=\pi$ or $\pi'=-\theta\overline\pi$. So we may ask the question:

Two Khovanov functors of weight $\pi$ and $-\theta\overline\pi$ are they isomorphic?
 
\vskip 12pt
\noi{\bf 5.11 Khovanov functors and Frobenius endomorphisms.} If $f$ is an endomorphism of the Frobenius algebra $R$, it acts on the complexes $KH(D)$ and transforms a
Khovanov functor $\Psi$ to a Khovanov functor $\Psi'$. Actually, if $(\sigma,E,F,W)$ and $(\sigma',E',F',W')$ are parametrizations of $\Psi$ and $\Psi'$, we have:
$$f(\sigma(h))=\lambda\sigma'(h)\hskip 48pt f(E_e(h))=E'_e(h)$$
$$f(F(h))=F'(h)\hskip 48pt f(W(h))={1\over\lambda}W'(h)$$
where $\lambda$ is the element in $K^*$ such that: $f(\omega)=\lambda\omega$.

Therefore, if $\pi$ is the weight of $\Psi$, the weight of $\Psi'$ is $f(\pi)$. Another consequence is the fact that there is no Khovanov functor invariant under
the endomorphisms of $R$ (at least if $R$ is the universal Frobenius algebra $R_0)$.
\vskip 12pt
\noi{\bf 5.12 Extensions of Khovanov functors.} In section 1.7, a category of mixed cobordisms ${\rde C}'$ was introduced. This category is a monoidal category containing
the category of cobordisms of closed oriented curves ${\rde C}$. Moreover the functor associated to $R$ extends to this category. Actually it is possible to define
in the same way a category of mixed cobordisms of oriented links ${\rde L}'$ containing the category ${\rde L}$. So we may ask the following:
\vskip12pt
\noi{\bf 5.13 Question:} Is it possible to extend a Khovanov functor to a monoidal functor from the category  ${\rde L}'$ of mixed cobordisms of oriented links to the
homotopy category of $K$-complexes? 

If such an extension $\Psi$ exists, the morphism $\Psi(f)$ associated to a mixed cobordism $f$ would be a chain map with a non necessarily
zero degree. Notice that, if $L$ is an oriented link and $L'$ is the same link but where the orientation of a sublink $L_1$ of $L$ is changed, there is a homotopy
equivalence from $KH(L)$ to $KH(L')$ of degree $2\lambda$, where $\lambda$ is the linking number between $L_1$ and $L\setminus L_1$.

Actually there is another category ${\rde C}''$ between ${\rde C}$ and ${\rde C}'$: the category of decorated cobordisms, where a decorated cobordism 
is a mixed cobordism on the form $(C,\emptyset,u)$ or equivalently a pair $(C,u)$ where $C$ is a cobordism decorated by a map $u:\pi_0(C)\rightarrow R$. Similarly there
is a category ${\rde L}''$ with: ${\rde L}\subset{\rde L}''\subset{\rde L}'$. Using the operators $\widehat T$ it is easy to extend every Khovanov
functor to the category  ${\rde L}''$, but the extension to ${\rde L}'$ is much more problematic.
\vskip 24pt
\noi{\bf 6. Invariant of knotted surfaces.}
\vskip 24pt
This section is devoted to the proof of theorem C.

Consider a closed oriented surface $S$ contained in $\R^4$. Up to isotopy, we may as well suppose that $S$ is contained in $\R^3\times(0,1)$. So this surface is a
cobordism in ${\rde L}$ from the empty link to itself. Therefore, any Khovanov functor $\Psi$ sends this surface to a morphism from $K$ to $K$ which is the
multiplication by an element $\Psi(S)\in K$. Because of proposition 5.7, we have: $\Psi(S)=\e(\delta^p\pi^{1-p})$ if $S$ is an unknotted connected surface of genus $p$
and $\pi$ is the weight of $\Psi$. That proves theorem C if $S$ is unknotted. 

In the classical case, the functor $\Psi$ was well defined up to sign and Tanaka [Ta] and Rasmussen [Ra2] proved that $\Psi(S)$ is $\pm2$ if $S$ is the torus and $0$ if
$S$ is any other connected surface. Notice that $\e(\delta^p)$ in the classical case is equal to $2$ if $p=1$ and to $0$ otherwise.

Let's say that $(R,\pi)$ is special if the following conditions hold:

$\bullet$ $\delta$ is invertible in $R$

$\bullet$ $\pi$ is a square in $R$

$\bullet$ $R$ has a twisting element which is the square of a element $x\in R$ with: $\overline x=x^{-1}$.
\vskip 12pt
\noi{\bf 6.1 Lemma:} {\sl Suppose theorem C is true if $(R,\pi)$ is special. Then the theorem is true in any case.}
\vskip 12pt
\noi{\bf Proof:} Consider a Frobenius algebra $R$ and a Khovanov functor $\Psi$ with weight $\pi\in R^*$. Let $R_1$ be the following ring:
$$R_1=\Z[\alpha,\overline\alpha,a,b,c,d,(a+b\alpha)^{-1},(a+b\overline\alpha)^{-1},(c+d\alpha)^{-1},(c+d\overline\alpha)^{-1}]$$
This ring is a localization of a polynomial ring with $6$ variables. It is equipped with an involution keeping $a,b,c,d$ fixed and exchanging $\alpha$ and 
$\overline\alpha$. Set:
$$\omega_1=a+b\alpha\hskip 24pt \pi_1=c+d\alpha$$
Then $R_1$ is a Frobenius algebra with generator $\alpha$ and twisting element $\omega_1$. It is easy to see that there is a unique morphism of Frobenius algebras $f$
from $R_1$ to $R$ sending $\alpha,\omega_1,\pi_1$ to $\alpha,\omega,\pi$ respectively.

Consider the following ring:
$$R_2=\Z[i][p,q,p_1,q_1,p_2,q_2,p^{-1},(1+q^2)^{-1},1/2,q_1^{-1},(p_2^2+q_2^2)^{-1}]\subset\C(p,q,p_1,q_1,p_2,q_2)$$
The complex conjugation induces an involution on $R_2$ and $R_2$ is a Frobenius algebra with generator $(p_1+iq_1)/p$ and twisting
element $p(1+iq)^2(1-iq)^{-2}$. It is easy to see that there is a unique morphism $g$ of Frobenius algebras $ĝ$ from $R_1$ to $R_2$ such that:
$$g(\alpha)={1\over p}(p_1+iq_1)\hskip 48ptg(\omega_1)=p\left({1+iq\over1-iq}\right)^2$$
$$g(\pi_1)=\pi_2=(p_2+iq_2)^2$$

It is clear that the fraction field of $R_2$ is, via $g$, an algebraic extension of the fraction field of $R_1$. Then $g$ is injective.

On the other hand, $R_2$ has $iq_1$ as generator and the corresponding twisting element is: $\omega_2=x^2$ with:
$$x={1+iq\over1-iq}$$
and $(R_2,\pi_2)$ is special.

Consider a Frobenius functor $\Psi_1$ associated with the Frobenius algeba $R_1$ with weight $\pi_1$. The morphisms $f$ and $g$ send $\Psi_1$ to functors $\Psi'$ and
$\Psi_2$ with weight $\pi$ and $\pi_2$.

If the theorem is true for special pairs, it is true for $\Psi_2$. Since $g:R_1\rightarrow R_2$ is injective, the theorem is also true for $\Psi_1$ and then for
$\Psi'$. But $\Psi'$ is isomorphic to any Khovanov functor with weight $\pi$. Therefore the theorem is true for $\Psi$.\cqfd
\vskip 12pt
Using this lemma, we may as well suppose that $\Psi$ is a Khovanov functor of weight $\pi$ and that $(R,\pi)$ is special. So $\beta=\alpha-\overline\alpha$ is invertible
and there exist $x$ and $y$ in $R^*$ such that:
$$\omega=x^2\hskip 48pt x\overline x=1\hskip 48pt \pi=y^2$$
and we may also suppose that $\Psi$ is parametrized by $(\sigma,E,F,W)$ with:
$$\sigma(h)=1\hskip 40pt E_e(h)=ehyx^{-1}\hskip 40pt F(a)=-ayx^{-1}\hskip 40pt W(h)=x^{-1}\otimes x^{-1}$$

From now on, we will suppose that all these properties are satisfied. The Khovanov data $(A,B,C,X,Y,Z)$ associated to $\Psi$ will be denoted by ${\rde K}$.

Consider an oriented link diagram $D$. Let $D^\circ$ be the union of $D$ and a trivial circle oriented clockwise and contained in half a plan disjoint from $D$. Let
$p$ be a point in this circle. The complex $KH(D^\circ,p,R)$ is naturally isomorphic to $R\otimes_K KH(D,R)$ and the graded module $E(D^\circ,p,R)$ is naturally isomorphic
to $E(D,R)$. Denote by $\varphi_D$ the composite map:
$$R\otimes_K KH(D,R)\ \ \build\longrightarrow_{}^\sim\ \ KH(D^\circ,p,R)\ \ \build\longrightarrow_{}^{\varphi(R)}\ \ E(D^\circ,p,R)\ \ \build\longrightarrow_{}^\sim\ \ 
E(D,R)$$
Because of theorem 2.15, this map is a homotopy equivalence and we have an explicit homotopy inverse of it (see remarks 2.14 and 2.16).

Then for every elementary move $f:D\rightarrow D'$ (compatible with the orientations) the morphism $f^{\rde K}: KH(D)\rightarrow KH(D')$
induces a well defined $R$-linear map $\widehat f: E(D,R)\rightarrow E(D',R)$. Since ${\rde K}$ is functorializable the correspondance $f\mapsto\widehat f$ extends
to composite of elementary moves (i.e. to movie sequences).

Consider an oriented curve $\Gamma$ in the plane and a component $C$ of $\Gamma$. We set:
$$\lambda(\Gamma,C)=\omega^{a(1+h)/2}<a|h>$$
where $a$ is the winding number of $C$ and $h$ is the $(\Gamma\minus C)$-sign of any point in $C$. 

If $D$ is an oriented link diagram, we set:
$$g(D)=\beta^{(m+q-n)/2}\prod_C\lambda(\widetilde D,C)$$
where the product holds for every component $C$ of the oriented resolution $\widetilde D$ of $D$, $q$ being the algebraic number of crossings of $D$ and $m$ (resp. $n$)
the number of components of $\widetilde D$ (resp. $D$).

Consider an elementary move $f:D\rightarrow D'$ between oriented diagrams. Denote by $\pi_0(D)$ the set of components of $D$ and by $\widehat D$ the set of maps from 
$\pi_0(D)$ to $\{\pm\}$. Such a map is called a $D$-state. For each state $\sigma\in\widehat D$, denote by $Y(\sigma)$ the set of crossings of $D$ between two 
components $c_1$ and $c_2$ with $\sigma(c_1)\not=\sigma(c_2)$, by $\widetilde D(\sigma)$ the oriented resolution of $D(\sigma)$ and by $d(\sigma)$ the number of 
components of $\widetilde D(\sigma)$. By using the same notations with the diagram $D'$, we have sets $\pi_0(D')$ and $\widehat D'$ and, for each $\sigma\in\widehat D'$,
a set $Y'(\sigma)$, a diagram $\widetilde D'(\sigma)$ and an integer $d'(\sigma)$.

This elementary move corresponds to a cobordism $C$ between links associated to $D$ and $D'$. So we have two maps from $\widehat C$ to $\widehat D$ and $\widehat D'$,
where $\widehat C$ is the set of $C$-states that is the set of maps from $\pi_0(C)$ to $\{\pm\}$. Let $\sigma$ (resp. $\sigma'$) be a $D$-state (resp. a $D'$-state).
We say that $\sigma$ and $\sigma'$ are compatible (or: $\sigma\sim\sigma'$) if they are coming from a state of $C$. A straightforward computation shows the following:
\vskip 12pt
\noi{\bf 6.2 Lemma:} {\sl Let $f:D\rightarrow D'$ and $\sigma$ (resp. $\sigma'$) be a $D$-state (resp. a $D'$-state). Suppose that $\sigma$ and $\sigma'$ aren't
compatible. Then the morphism:
$$\Lambda^{-e}(Y(\sigma))\otimes R v(\sigma)\rightarrow E(D,R)\build\longrightarrow_{}^{\widehat f} E(D',R)\rightarrow  \Lambda^{-e}(Y'(\sigma'))\otimes R v(\sigma')$$
is trivial.}
\vskip 12pt
An immediate consequence of this lemma is the following:
\vskip 12pt
\noi{\bf 6.3 Lemma:} {\sl Let $f:D\rightarrow D'$ be a movie sequence represented by a cobordism $C$. Let $\sigma$ and $\sigma'$ be a $D$-state and a $D'$-state. Suppose
the composite map:
$$\Lambda^{-e}(Y(\sigma))\otimes R v(\sigma)\rightarrow E(D,R)\build\longrightarrow_{}^{\widehat f} E(D',R)\rightarrow  \Lambda^{-e}(Y'(\sigma'))\otimes R v(\sigma')$$
is not trivial. Then the two states $\sigma$ and $\sigma'$ are coming from a $C$-state.}
\vskip 12pt
Let $f:D\rightarrow D'$ be a movie sequence represented by a cobordism $C$ and $\tau$ be a $C$-state. This state restricts to a $D$-state $\sigma$ and
a $D$'-state $\sigma'$. Define the map $\widehat f(\tau)$ by:
$$\widehat f(\tau)=pr'\circ\widehat f\circ pr$$
where $pr$ (resp. $pr'$) is the projection $E(D,R)\rightarrow \Lambda^{-e}(Y(\sigma))\otimes R v(\sigma)\subset E(D,R)$ (resp. $E(D',R)\rightarrow 
\Lambda^{-e}(Y'(\sigma'))\otimes R v(\sigma')\subset E(D',R)$). The map $\widehat f(\tau)$ vanishes on $\Lambda^{-e}(Y(\sigma_1))\otimes R v(\sigma_1)$ for 
$\sigma_1\not=\sigma$ and is the composite:
$$\Lambda^{-e}(Y(\sigma))\otimes R v(\sigma)\rightarrow E(D,R)\build\longrightarrow_{}^{\widehat f} E(D',R)\rightarrow  \Lambda^{-e}(Y'(\sigma'))\otimes R v(\sigma')
\subset E(D',R)$$
on $\Lambda^{-e}(Y(\sigma))\otimes R v(\sigma)$. 

We have clearly the following:
$$\widehat f=\sum_\tau\widehat f(\tau)$$

If $f$ corresponds to a closed surface $S$, $\widehat f(\tau)$ is the multiplication by an element $\widehat S(\tau)\in R$. Clearly $\Psi(S)$ is the sum of all 
$\widehat S(\tau)$.

Consider now an elementary move $f:D\rightarrow D'$ corresponding to a cobordism $C$. Consider a $C$-state $\tau$. This state restricts to a $D$-state $\sigma$ and
a $D$'-state $\sigma'$. A straightforward computation, case by case, gives the following:
\vskip 12pt
\noi{\bf 6.4 Lemma:} {\sl Suppose $f$ is a Reidemeister move of type I$^+(e,a,h)$. Denote by $c$ the component of $D$ which is modified by $f$.
Then the morphism $\widehat f(\tau)$ is the map:
$$u\otimes v(\sigma)\mapsto u\otimes wv(\sigma')$$
with:
$$w=\Bigl({g(D')\over g(D)}\Bigr)^{(\sigma(c))}={g(D'(\sigma'))\over g(D(\sigma))}<\sigma(c)|-eh>$$}
\vskip 12pt
\noi{\bf 6.5 Lemma:} {\sl Suppose $f$ is a Reidemeister move of type II$^+(a,b,h)$. Let $x_+$ and $x_-$ be the created crossings with sign $+$ and $-$. Denote by $c_+$ 
(resp. $c_-$) the component of $D$ containing the over (resp. under) branch of the move. Then the morphism $\widehat f(\tau)$ is the map:
$$u\otimes v(\sigma)\mapsto u\otimes wv(\sigma')$$
with:
$$w=\Bigl({g(D')\over g(D)}\Bigr)^{(\sigma(c_+))}={g(D'(\sigma'))\over g(D(\sigma))}<\sigma(c_+)|-ab>$$
if $\sigma(c_+)=\sigma(c_-)$ and the map:
$$u\otimes v(\sigma)\mapsto x_+\vc x_-\vc u\otimes wv(\sigma')$$
with:
$$w={g(D'(\sigma'))\over g(D(\sigma))}<\sigma(c_+)|b><\sigma(c_-)|a><\sigma(c_+)|\sigma(c_-)>x^{-h(a\sigma(c_+)+b\sigma(c_-))(1-<\sigma(c_+)\sigma(c_-)|ab>)}$$
otherwise.}
\vskip 12pt
\noi{\bf 6.6 Lemma:} {\sl Suppose $f$ is a Reidemeister move of type III$(e,a,b,c,h)$. Let $c_1$ (resp. $c_2$, $c_3$) be the over branch (resp. the middle branch, the
under branch) of the move. Then the morphism $\widehat f(\tau)$ is the map:
$$u\otimes v(\sigma)\mapsto u\otimes wv(\sigma')$$
with:
$$w=\Bigl({g(D')\over g(D)}\Bigr)^{(\sigma(c_1))}={g(D'(\sigma'))\over g(D(\sigma))}$$
if $\sigma(c_1)=\sigma(c_2)=\sigma(c_3)$, and:
$$w={g(D'(\sigma'))\over g(D(\sigma))}\sigma(c_1)\sigma(c_3)$$
in the general case.}
\vskip 12pt
\noi{\bf 6.7 Lemma:} {\sl Suppose $f$ is a surgery move of type $(0,a,h)$. Let $c$ be the created circle. Then the morphism $\widehat f(\tau)$ is the map:
$$u\otimes v(\sigma)\mapsto u\otimes wv(\sigma')$$
with:
$$w=\Bigl({g(D')\over g(D)}yx^{-1}\Bigr)^{(\sigma'(c))}={g(D'(\sigma'))\over g(D(\sigma))}\Bigl(yx^{-1}\Bigr)^{(\sigma'(c))}$$}
\vskip 12pt
\noi{\bf 6.8 Lemma:} {\sl Suppose $f$ is a surgery move of type $(1,a,h)$. Let $c$ be the component of the cobordism containing the modified branches. Let: $e=1$ (resp.
$e=0$) if the surgery increases (resp. decreases) the number of components of the link. Then the morphism $\widehat f(\tau)$ is the map:
$$u\otimes v(\sigma)\mapsto u\otimes wv(\sigma')$$
with:
$$w=\Bigl({g(D')\over g(D)}\beta^e xy^{-1}\Bigr)^{(\tau(c))}={g(D'(\sigma'))\over g(D(\sigma))}\beta^e\Bigl(xy^{-1}\Bigr)^{(\tau(c))}<\tau(c)|h>$$}
\vskip 12pt
\noi{\bf 6.9 Lemma:} {\sl Suppose $f$ is a surgery move of type $(2,a,h)$. Let $c$ be the component removed by the surgery. Then the morphism $\widehat f(\tau)$ is the
map:
$$u\otimes v(\sigma)\mapsto u\otimes wv(\sigma')$$
with:
$$w=\Bigl({g(D')\over g(D)}\beta^{-1}yx^{-1}\Bigr)^{(\sigma(c))}={g(D'(\sigma'))\over g(D(\sigma))}\beta^{-1}\Bigl(xy^{-1}\Bigr)^{(\sigma(c))}<\sigma(c)|-h>$$}
\vskip 12pt
An easy consequence of these lemmas is the following:
\vskip 12pt
\noi{\bf 6.10 Lemma:} {\sl Let $f:D\rightarrow D'$ be a movie sequence corresponding to a cobordism $C$. Let $\tau$ be a constant $C$-state sending each component of $C$
to a sign $e$. Let $n_0$ (resp. $n_2$) be the number of surgery moves in $f$ of index $0$ (resp. $2$). Let $n_1^+$ (resp. $n_1^-$) be the number of surgery moves in $f$
of index $1$ which increases (resp. decreases) the number of components of the link. Then the morphism $\widehat f(\tau)$ is the map:
$$u\otimes v(\sigma)\mapsto u\otimes w^{(e)}v(\sigma')$$
with:
$$w={g(D')\over g(D)}\beta^{n_1^+-n_2}(yx^{-1})^{n_0+n_2-n_1^+-n_1^-}$$}
\vskip 12pt
\noi{\bf 6.11 Corollary:} {\sl Let $S$ be a surface in $\R^4$. Let $p_i$ be the genus of the $i$-th component of $S$. Let $\tau$ be a constant $S$-state sending each
component of $S$ to a sign $e$. Then we have:
$$\widehat S(\tau)=\Bigl(\prod_i\Bigl({\delta\over\pi}\Bigr)^{p_i-1}\Bigr)^{(e)}$$}
\vskip 12pt
\noi{\bf Proof:} Let $k$ be the number of components of $S$ and $p$ be the sum of the $p_i$'s. Numbers $n_0,n_1^+,n_1^-,n_2$ are related with a handle decomposition of
$S$ and there exist two integers $a,b\geq0$ such that:
$$n_0=k+a\hskip 24pt n_1^+=p+b\hskip 24pt n_1^-=p+a\hskip 24pt n_2=k+b$$
So we have:
$$\beta^{n_1^+-n_2}(yx^{-1})^{n_0+n_2-n_1^+-n_1^-}=\beta^q(xy^{-1})^{2q}=(\beta x^2y^{-2})^q=(\delta/\pi)^q$$
with $q=p-k=\sum_i(p_i-1)$. The result follows.\cqfd
\vskip 12pt
Let $f:D\rightarrow D'$ be a movie sequence corresponding to a cobordism $C$ and $\tau$ be a $C$-state. Define the sign $s(f,\tau)$ by the following:

Suppose $f$ is a Reidemeister mode of type III. Denote by $c_1$ (resp. $c_2$, $c_3$) the component of $C$ containing the top branch (resp. the middle
branch, the bottom banch) of the move. In this case we set: $s(f,\tau)=-1$ if $\tau(c_1)=\tau(c_3)=-\tau(c_2)$ and  $s(f,\tau)=1$ otherwise.

If $f$ is another elementary move we set: $s(f,\tau)=1$.

If $f$ is a movie sequence: $f=(f_1,f_2,\dots,f_p)$, we set: $s(f,\tau)=\prod_i s(f_i,\tau)$.

If $S$ is a closed oriented surface in $\R^4$, $S$ corresponds to a movie sequence $f$ and we set: $s(S,\tau)=s(f,\tau)$.
\vskip 12pt
\noi{\bf 6.12 Lemma:} {\sl Let $S$ be a surface in $\R^4$ and $\tau$ be a $S$-state. Let $S_i$ be the $i$-th component of $S$, $p_i$ be the genus of $S_i$ and $e_i$ be
the sign $\tau(S_i)$. Then we have:
$$\widehat S(\tau)=s(S,\tau)\prod_i\Bigl(\Bigl({\delta\over\pi}\Bigr)^{p_i-1}\Bigr)^{(e_i)}$$}
\vskip 12pt
\noi{\bf Proof:} Denote by $S_+$ (resp. $S_-$) the submanifold of $S$ where $\tau$ is equal to $+$ (resp. $-$). By moving down $S_-$ along the vertical axis in $\R^3$,
we get a new surface $S'$ which is isotopic to $S_+\coprod S_-$. The $S$-state $\tau$ induces a $S'$-state $\tau'$. Because of the corollary we have:
$$\widehat S'(\tau')=\widehat S_+(+)\widehat S_-(-)=\prod_i\Bigl(\Bigl({\delta\over\pi}\Bigr)^{p_i-1}\Bigr)^{(e_i)}$$

Suppose the morphism $f$ corresponding to $S$ is a movie sequence: $f=(f_1,\dots,f_k)$ where $f_i$ is an elementary move from a diagram $D_{i-1}$ to a diagram $D_i$.
Then the morphism corresponding to $S'$ is a movie sequence: $f'=(f'_1,f'_2,\dots,f'_k)$ where $f'_i$ is an elementary move from a diagram $D'_{i-1}$ to a diagram $D'_i$.
For every $i=1,2,\dots,k$ we have:
$$\widehat f_i(\tau)={g(D_i(\sigma_i))\over g(D_{i-1}(\sigma_{i-1}))}\varphi_i(\tau)$$
$$\widehat f'_i(\tau')={g(D'_i(\sigma'_i))\over g(D'_{i-1}(\sigma'_{i-1}))}\varphi'_i(\tau')$$
Using lemmas 6.4 to 6.9, we check that $\varphi_i(\tau)$ and $\varphi'_i(\tau')$ are allways the same except for type III Reidemeister moves.  In these cases, we have:
$$\varphi'_i(\tau')=s(f_i,\tau)\varphi_i(\tau)$$
Thus we have:
$$\widehat S(\tau)=\prod_i\varphi_i(\tau)=\prod_i s(f_i,\tau)\varphi'_i(\tau')=s(f,\tau)\prod_i\varphi'_i(\tau')=s(S,\tau)\widehat S'(\tau')$$
and the result follows.\cqfd
\vskip 12pt
\noi{\bf 6.13 Lemma:} {\sl Let $S$ be a closed oriented surface in $\R^4$ and $\tau$ be a $S$-state. Then we have:
$$s(S,\tau)=1$$}
\vskip 12pt
\noi{\bf Proof:} Because of lemma 6.12, $s(S,\tau)$ is invariant under isotopy. But $s(S,\tau)$ is also invariant under surgery moves. Therefore $s(S,\tau)$ depends
only on the cobordism class of $(S,\tau)$ in the group $\Omega$ of cobordisms of bicolored surfaces in $\R^4$. The Pontryagin-Thom construction implies:
$$\Omega=\pi_4(MSO_2\vee MSO_2)=\pi_4(MU_1\vee MU_1)=\pi_4(BU_1\vee BU_1)$$
Denote by $E$ the space of paths in $BU_1$ ending at the base point. The homotopy fiber $F$ of the inclusion $BU_1\vee BU_1\subset BU_1\times BU_1$ is the following:
$$F=E\times\Omega BU_1\build\cup_{\Omega BU_1\times\Omega BU_1}^{}\Omega BU_1\times E$$
But we have a homotopy equivalence: $(E,\Omega BU_1)\sim (B^2,S^1)$. So we get a homotopy equivalence:
$$F\sim B^2\times S^1\build\cup_{S^1\times S^1}^{} S^1\times B^2=S^3$$
and we have:
$$\Omega\simeq \pi_4(F)\simeq\pi_4(S^3)\simeq\Z/2$$
Since $s(S,\tau)$ depends only on the class of $(S,\tau)$ in $\Omega$, we have a map $\varphi: \Omega\rightarrow \{\pm\}$ such that:
$$s(S,\tau)=\varphi([S,\tau])$$
where $[S,\tau]$ is the cobordism class of $(S,\tau)$. By testing this formula for the empty surface we get: $\varphi(0)=1$. Then the last thing to do is to determine
$s(S,\tau)$ for some bicolored surface which is not trivial in $\Omega$. Following [CKSS], the non zero element in $\Omega$ is represented by two tori $T_+$ and $T_-$
where $T_+$ intersects $\R^3\times\{0\}$ in a Hopf link $H$ and $T_-$ is the boundary of a regular neighborough of one component of $H$ in $\R^3\times\{0\}$.

Consider a movie sequence $f$ from the empty diagram to the diagram $D$ of a Hopf link, given by a $0$-surgery, two Reidemeister moves of type I$^+_+$ and a $1$-surgery.
Denote by $\overline f$ the inverse move. So $T_+$ is represented by the movie move $(f,\overline f)$. Consider a movie sequence $f_1$ from $D$ to a diagram $D_1$
given by a $0$ surgery, two Reidemeister moves of type II$^+$ and a $1$-surgery. This movie creates two circles $C_1$ and $C_2$ in a neighborough of a component $C$ of
the diagram $D$. By moving $C_1$ around $C$, this circle goes through the other component of $D$ and we have a movie sequence $f_2$ from $D_1$ to a diagram $D_2$ given
by a Reidemeister move of type II$^+$, two Reidemeister moves of type III and a Reidemeister move of type II$^-$. By moving $C_2$ around the other part of $C$, we get
a movie sequence $f_3$ from $D_2$ to $D_3$ given also by a Reidemeister move of type II$^+$, two Reidemeister moves of type III and a Reidemeister move of type II$^-$.
Finally we have a movie sequence $f_4$ from $D_3$ to $D$ given by a $1$-surgery, two Reidemeister moves of type II$^-$ and a $2$-surgery. The diagrams $D$, $D_1$, $D_2$
and $D_3$ are the following:
$$
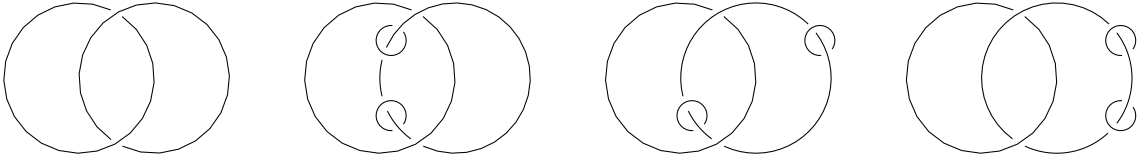
\begin{tikzpicture}[scale=1/10] \draw [domain=6.5:41.5] plot({10+10*cos(10*\x)},{10+10*sin(10*\x)});
\draw [domain=-11.5:23.5] plot({20+10*cos(10*\x)},{10+10*sin(10*\x)});
\draw [domain=6.5:41.5] plot({50+10*cos(10*\x)},{10+10*sin(10*\x)});\draw [domain=-11.5:15.6] plot({60+10*cos(10*\x)},{10+10*sin(10*\x)});
\draw [domain=16.6:19.4] plot({60+10*cos(10*\x)},{10+10*sin(10*\x)});\draw [domain=20.6:23.4] plot({60+10*cos(10*\x)},{10+10*sin(10*\x)});
\draw [domain=8.5:39.5] plot ({51.5+2*cos(10*\x)},{15+2*sin(10*\x)});\draw [domain=8.5:39.5] plot ({51.5+2*cos(10*\x)},{5-2*sin(10*\x)});
\draw [domain=6.5:41.5] plot({90+10*cos(10*\x)},{10+10*sin(10*\x)});\draw [domain=-11.5:3.7] plot({100+10*cos(10*\x)},{10+10*sin(10*\x)});
\draw [domain=4.6:19.4] plot({100+10*cos(10*\x)},{10+10*sin(10*\x)});\draw [domain=20.6:23.4] plot({100+10*cos(10*\x)},{10+10*sin(10*\x)});
\draw [domain=8.5:39.5] plot ({108.5+2*cos(10*\x)},{15-2*sin(10*\x)});\draw [domain=8.5:39.5] plot ({91.5+2*cos(10*\x)},{5-2*sin(10*\x)});
\draw [domain=6.5:41.5] plot({130+10*cos(10*\x)},{10+10*sin(10*\x)});\draw [domain=-11.5:-4.7] plot({140+10*cos(10*\x)},{10+10*sin(10*\x)});
\draw [domain=-3.7:3.7] plot({140+10*cos(10*\x)},{10+10*sin(10*\x)});\draw [domain=4.6:23.4] plot({140+10*cos(10*\x)},{10+10*sin(10*\x)});
\draw [domain=8.5:39.5] plot ({148.5+2*cos(10*\x)},{15-2*sin(10*\x)});\draw [domain=8.5:39.5] plot ({148.5+2*cos(10*\x)},{5+2*sin(10*\x)});
\end{tikzpicture}$$

The movie sequence $(f,f_1,f_2,f_3,f_4,\overline f)$ represents a bicolored surface $S$ which is not zero in $\Omega$. On the other hand, among the movies $f$,
$\overline f$ and the $f_i$, only $f_2$ and $f_3$ contains some type III Reidemeister moves. So we have:
$$s(S,\tau)=s(f_2,\tau)s(f_3,\tau)=(-1)(-1)=1$$
Then the map $\varphi$ is trivial on $\Omega$ and the result follows.\cqfd
\vskip 12pt
Now we are able to prove theorem C. Let $S$ be a closed oriented surface in $\R^4$. Denote by $S_i$ the $i$-th component of $S$ and by $p_i$ the genus of $S_i$. Denote
also by $u_i$ the element $(\delta/\pi_i)^{p_i-1}$. A $S$-state is characterized by the signs $e_i=\tau(S_i)$. So we have:
$$\Psi(S)=\sum_{\tau}\widehat S(\tau)=\sum_{e_*}\prod_i\Bigl(\Bigl({\delta\over\pi}\Bigr)^{p_i-1}\Bigr)^{(e_i)}=\sum_{e_*}\prod_i u_i^{(e_i)}=\prod_i(u_i+\overline u_i)$$
and the desired result follows from the obvious relation:
$$\forall u \in R,\ \ \ u+\overline u=\e(\delta u)\eqno{\cqfd}$$
\vskip 12pt
\noi{\bf 6.14 Remark:} We can extend the functor $\Psi$ to the category ${\rde L}''$ of decorated cobordisms of links. Consider a closed oriented surface $S$ in $\R^4$
decorated by $u$. This decoration sends each component $S_i$ of $S$ to an element $u_i\in R$. In this case we have:
$$\Psi(S,u)=\prod_i\e(u_i\delta^{p_i}\pi^{1-p_i})$$
where $p_i$ is the genus of $S_i$. In any case $\Psi(S,u)$ doesn't depend on the embedding $S\subset\R^4$.
\vskip 24pt
\noi{\bf References: }

\begin{list}{}{\leftmargin 45pt \labelsep 10pt \labelwidth 40pt \itemsep 0pt}
\item[{[Bl]}] Christian Blanchet -- {\sl An oriented model for Khovanov homology}, Journal of Knot Theory and its Ramifications Vol 19, N$^\circ$ 2 (2010), 291--312,
math.GT/1405.7246.
\item[{[BM]}] Dror Bar Natan, Scott Morrison -- {\sl The Karoubi Envelope and Lee's Degeneration of Khovanov Homology}, Alg. Geom. Topol. {\bf 6} (2006), 1459--1469,
math.GT/0606542.
\item[{[BN1]}] Dror Bar Natan -- {\sl On Khovanov's categorification of the Jones polynomial}, Alg. Geom. Topol. {\bf 2} (2002), 337--370, math.QA/0201043.
\item[{[BN2]}] Dror Bar Natan -- {\sl Khovanov's homology for tangles and cobordisms}, Geom. Topol. {\bf 9} (2005), 1443--1499, math.GT/0410495.
\item[{[CC]}] Carmen Livia Caprau -- {\sl sl(2) tangle homology with a parameter and singular cobordisms}, Alg. Geom. Topol. {\bf 8} (2008), 729--756.
\item[{[CKSS]}] J. Scott Carter, Seiichi Kamada, Masahico Saito, Shin Satoh -- {\sl A theorem of Sanderson on link bordisms in dimension $4$}, Alg. Geom. Topol. {\bf 1} 
(2001), 299--310.
\item[{[CMW]}] Davis Clark, Scott Morrison, Kevin Walker -- {\sl Fixing the functoriality of Khovanov homology}, Geom. Topol. {\bf 13} (2009), 1499--1582, 
math.GT/0701339.
\item[{[CS]}] J. Scott Carter, Masahico Saito -- {\sl Reidemeister moves for surface isotopies and their interpretation as moves to movies}, Journal of Knot Theory and 
its Ramifications Vol 2, N$^\circ$ 3 (1993), 251--284.
\item[{[Ja]}] Magnus Jacobson -- {\sl An invariant of link cobordisms from Khovanov homology}, Alg. Geom. Topol. {\bf 4} (2004), 1211--1254.
\item[{[Kh1]}] Mikhail Khovanov -- {\sl A categorification of the Jones polynomial}, Duke Math. J. Vol 101, N$^\circ$ 3 (2000), 359--426, math.QA/9908171.
\item[{[Kh2]}] Mikhail Khovanov -- {\sl A functor--valued invariant of tangles}, Alg. Geom. Topol. {\bf 2} (2002), 665--741, math.QA/0103190.
\item[{[Kh3]}] Mikhail Khovanov -- {\sl An invariant of tangle cobordism}, Trans. Amer. Math. Soc. 358 (2006), 315--327.
\item[{[Kh4]}] Mikhail Khovanov -- {\sl Link homology and Frobenius extensions}, Fund. Math. 190 (2006), 179--190, math.QA/0411447.
\item[{[Ko]}] Joachim Kock -- {\sl Frobenius algebras and 2D Topological Quantum Field Theories}, LMS Student Texts 59 (2003), Cambridge University Press.
\item[{[Le]}] Eun Soo Lee -- {\sl Khovanov's invariants for alternating links}, math.GT/0210213 (2002).
\item[{[Ra1]}] Jacob Rasmussen -- {\sl Khovanov homology and the slice genus}, Invent. Math. 182 (2010), 419--447, math.GT/0402131.
\item[{[Ra2]}] Jacob Rasmussen -- {\sl Khovanov's invariants for closed surfaces}, math.GT/0502527 (2005).
\item[{[Ta]}] Kokoro Tanaka -- {\sl Khovanov-Jacobsson numbers and invariants of surface-knots derived from Bar-Natan's theory}, Proc. Amer. Math. Soc. 134 no. 12 
(2006), 3685--3689.                                                                                               
\end{list}
\end{document}